\documentclass[reqno,10pt]{amsart}
\usepackage{amsmath,amsfonts,amssymb,amsxtra,latexsym,amscd,enumerate,enumitem,amsthm,verbatim,here}

\usepackage[dvipdfmx]{graphicx}
\usepackage{cancel}

\usepackage{mathtools}

\usepackage[dvipsnames]{xcolor}

\usepackage{hyperref}
\hypersetup{colorlinks=true, pdfstartview=FitV,linkcolor=BrickRed,citecolor=BrickRed, urlcolor=BrickRed}
\definecolor{labelkey}{rgb}{0.6,0,0}

\usepackage[margin=1.2in]{geometry}
\numberwithin{equation}{section}

\def\R{\mathbb{R}}

\newtheorem{theorem}{Theorem}[section]
\newtheorem{lemma}[theorem]{Lemma}
\newtheorem{corollary}[theorem]{Corollary}

\newtheorem{proposition}[theorem]{Proposition}
\newtheorem{definition}[theorem]{Definition}
\theoremstyle{remark}
\newtheorem{remark}[theorem]{Remark}

\title[VP with a wall]{The Vlasov-Poisson system with a perfectly conducting wall:
Convex domains.}

\author{Wenrui Huang}
\address{Brown University}
\email{wenrui\_huang@brown.edu}

\author{Benoit Pausader}
\address{Brown University}
\email{benoit\_pausader@brown.edu}

\author{Masahiro Suzuki}
\address{Nagoya Institute of Technology}
\email{masahiro@nitech.ac.jp}

\begin{document}

\date{\today}

\begin{abstract}
We consider the Vlasov--Poisson system in a $C^3$ convex domain $D$ with a perfectly conducting wall. We introduce the asymptotic domain $D_{\infty}$ for the domain $D$. Then under acceptable assumptions on $D$, we show that for localized initial data,  the velocity of particles is asymptotically supported in the (closure of the) asymptotic domain $\overline{D_{\infty}}$ and the solutions exhibit the asymptotics of modified scattering.
\end{abstract}

\maketitle

\section{Introduction}\label{S1}

The interaction of plasma with boundary is a central question in plasma physics, 
since it occurs in many applications of plasma such as neon lamps, microfabrication of semiconductor, and nuclear fusion.
In this work, we investigate the asymptotic behavior of solutions to the Vlasov-Poisson system in the simplest case of a perfectly conducting wall for a rarefied plasma.

\subsection{The Vlasov-Poisson system with a perfectly conducting wall}

We consider the Vlasov-Poisson system inside a convex domain $D$ with a perfectly conducting wall $\partial D$:
\begin{equation}\label{VPW}
\begin{split}
\partial_t\widetilde{\mu}+\{\widetilde{\mu},\mathcal{H}\}=0,&\,\,\mathcal{H}:=\frac{\vert v\vert^2}{2}+\lambda\phi\,\hbox{ in }\,D\times\mathbb{R}^{3}_{v},
\\
\widetilde{\mu}=\widetilde{\mu}_{0}\hbox{ at }t=0,&\quad \widetilde{\mu}=0 \hbox{ when }x\in \partial D\hbox{ and }{\bf n}\cdot v \ge 0,
\end{split}
\end{equation}
where ${\bf n}$ is the inward normal vector,  $\widetilde{\mu}(x,v,t):D\times\mathbb{R}^3\times\mathbb{R}_+\to\mathbb{R}$ is the particle distribution function, the electrostatic potential $\phi(x,t): D \times \mathbb{R}_{+} \to \mathbb{R}$ is defined by
\begin{equation}\label{VPPhi}
\Delta_x\phi=\int_{\mathbb{R}^3}\widetilde{\mu}^2(x,v,t)dv\hbox{ in }D,\quad\phi=0\hbox{ on }\partial D,\quad \phi \to0 \,\,\hbox{ as }\,\,\vert x\vert\to\infty,
\end{equation}
and we use the Poisson bracket
\begin{equation*}
\begin{split}
\{f,g\}=\nabla_xf\cdot\nabla_vg-\nabla_vf\cdot\nabla_xg.
\end{split}
\end{equation*}
We focus on the {\it plasma} case when $\lambda<0$. Although we are considering solutions close to vacuum, this distinction will be important in the analysis.

{We assume that $D$ is an admissible convex domain with a $C^3$ boundary that contains the origin. 
Roughly speaking, $D$ is unbounded and its asymptotically approaches an asymptomatic domain $D_\infty$ at large distances from the origin as Figure \ref{Fig1} below. 
The boundary $\partial D_\infty$ is the largest cone with its vertex at the origin and an infinite height, which is contained in $\overline{D}$.
Specifically, the distance between $\partial D$ and $\partial D_\infty$ becomes almost constant at large distances from the origin. 
The admissible convex domain also satisfies some other conditions, and its precise definition will be given in Definition \ref{AcceptableDomains} below. A half-space $\mathbb {R}^{3}_{+}$ is an example of an admissible convex domain with $D_\infty=\mathbb {R}^{3}_{+}$. Another example is a hyperbolic cylinder as Figure \ref{Fig1} (see Appendix \ref{hyper}).}

\begin{figure}[H]
	\centering
	\includegraphics[scale=0.7]{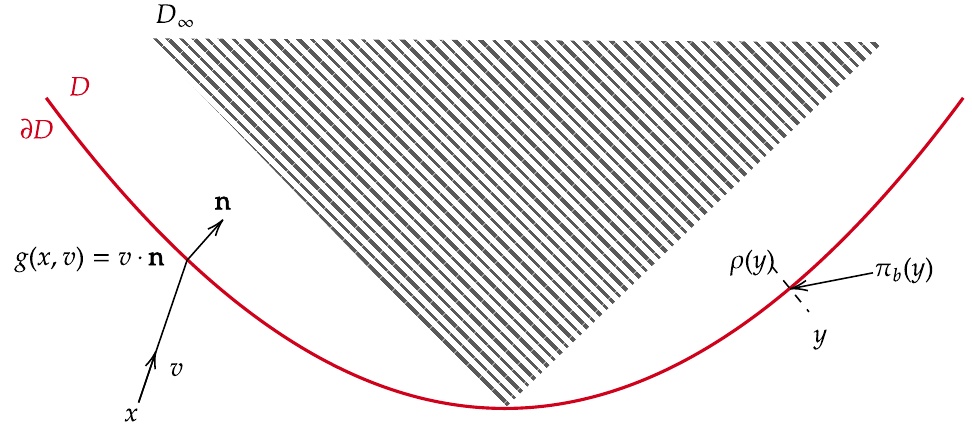}
	\caption{An illustration of the domain geometry}
	\label{Fig1}
\end{figure}

Our main result is the following, {where $\langle t\rangle:=(1+t^{2})^{1/2}$ and  $\langle x \rangle:=(1+|x|^{2})^{1/2}$ for $t \in \mathbb R$ and $x \in \mathbb R^{3}$}:

\begin{theorem}\label{MainThm}
Let $D\subset\mathbb{R}^3$ be an acceptable convex domain (see Definition \ref{AcceptableDomains} below). There exists $\varepsilon>0$ such that if the initial data $\widetilde{\mu}_0\in C^1_{x,v}(\overline{D}\times\mathbb{R}^3)$ satisfies the compatibility conditions of the zeroth and first orders as well as
\begin{equation*}
\begin{split}
\Vert\langle x\rangle^{4}\langle v\rangle^{4}\widetilde{\mu}_0\Vert_{L^\infty_{x,v}}+\Vert\nabla_{x,v}\widetilde{\mu}_0\Vert_{L^\infty_{x,v}}\le\varepsilon,
\end{split}
\end{equation*}
then there exists a unique global classical solution $\widetilde{\mu}\in C^1_{x,v,t}(D\times\mathbb{R}^3\times \mathbb{R}_{+})$ of the problem \eqref{VPW}. Besides the electrostatic potential $\phi$ decays at the optimal rate
\begin{equation*}
\begin{split}
\sqrt{1+\vert x\vert^2+t^2}\vert\phi(x,t)\vert+\left[1+\vert x\vert^2+t^2\right]\vert \nabla_x\phi(x,t)\vert\lesssim\varepsilon^2 \quad \hbox{for } x\in D, \ t>0.
\end{split}
\end{equation*}
In addition, we have modified scattering in the following sense: there {exists} an asymptotic domain $D_{\infty}\subset D$  and a limiting distribution
$\mu_\infty\in L^2_{x,v}\cap L^\infty_{x,v}\cap C^{0}_{x,v}(\mathbb{R}^3\times D_\infty)$ such that\footnote{The first line means that for any compact set $K\subset D\times D_\infty$, $\Vert \widetilde{\mu}(x,v,t)-\mu_\infty(x-tv-\lambda(\ln(t))E_\infty(v),v)\Vert_{L^\infty(K)}\to0$.}
\begin{gather*}
\widetilde{\mu}(x,v,t)\rightarrow \mu_\infty(x-tv-\lambda(\ln t) E_\infty(v),v) \quad a.e \ (x,v) \in D \times D_{\infty} \hbox{ as } t \to \infty,
\\
\Vert \widetilde{\mu}\Vert_{L^2_{x,v}(D\times D_\infty^c)}\le\varepsilon \langle t\rangle^{-1/4} \quad \hbox{for } t>0,
\end{gather*}
where $E_\infty$ is defined by the distribution $\mu_\infty$ and  the Dirichlet Green function $G_\infty$ of the domain  $D_\infty$ as
\begin{equation*}
E_\infty^{j}(a):=\iint_{\mathbb{R}^3 \times D_\infty} \partial_{a_{j}}G_\infty(a,w)\mu_\infty^{2}(y,w)dydw \quad \hbox{for }j=1,2,3.
\end{equation*}
\end{theorem}

\subsection{Prior works and general comments}

\subsubsection{Prior works} The literature of broadly related results is too vast to be surveyed here, so we will focus on a selection of more directly related results.
\begin{itemize}
	\item In the whole space case, we refer to earlier works \cite{BD1985, IoPaWaWiVacuum, Smu, Wang} for references on global wellposedness and dispersion, and to recent works \cite{FlOuPaWi,PaWi,PWY2022} for references on the symplectic approach to the study of asymptotics. We also refer to \cite{Big2022, HuKw, IaRoWi} for references on the study of asymptotics of some other Vlasov equations.
	\item In the presence of  boundary, the work \cite{Guo2} studied the simple case of a half line where the boundary is just a point. It constructed singular and regular solutions for different boundary conditions. The method there was then adapted in \cite{Guo,Hwa,HwaVel1,HwaVel2} to study three space dimensions. We also refer to \cite{CeIa,FerRea} for more recent works on global existence of solutions, also in a more general context. However, none of these studied the asymptotics.     \item In slightly different contexts, {the studies \cite{Bad2024,ST2,STZ1}} focus on plasma sheaths by considering different conditions at infinity and on the boundary. When the plasma contacts with a wall, there appears a non-neutral potential region between the wall and plasma, and a nontrivial steady state is achieved. This non-neutral region is referred to as a {\it sheath}. 
For more detailed discussion on the plasma sheath, see the references in \cite{ST2,STZ1}. The study \cite{ST2} classifies for all possible cases whether or not the steady state exists over the half line. In addition, \cite{STZ1} investigates the conditional stability of the steady states. {The study \cite{Bad2024} also analyzed the stability of the steady states in a bounded internal.}
\end{itemize}
As a consequence of these work, one knows that solutions are global in large generality, and the main novelty of Theorem \ref{MainThm} concerns the study of asymptotics.  Our goal here is to introduce robust methods and exhibit the important structural objects that are relevant in this study, and we hope in other generalizations of this problem such as the problems we mention at the end of the introduction. These include the asymptotic domain $D_\infty$ and the grazing function $g$ in Section \ref{SecDInfty}.

\subsubsection{Assumptions on the domain}

A half-space is certainly the most important example of a convex domain. In view of applications, it is also natural to consider general smooth convex domains (see Section \ref{S22}). 
Analytically, considering a convex domain has a number of technical advantages, including the fact that the Dirichlet Green function is slightly better behaved (see \cite{Fr}), that we have a natural foliation of the exterior of the domain (related to the function $b$ from \eqref{DefB}), which is useful to extend the domain, and that there is a simple blow-down analysis in the sense that rescaling of the domains $D_t$ form a monotone family \eqref{DInfty}.

The case of a general convex domain would be interesting, especially if the boundary has limited regularity. 
As a first step, it is natural to consider the case of a corner $\{x \in \mathbb R^{3} \, : \, x^{1}>0, \ x^{2}>0 \}$, since we can write explicitly the Dirichlet Green function. 
However, in this case, {it is easy to show that the crucial estimate \eqref{T1Hyp} is false and our method leads to worse estimates on the derivative.} 

Considering non-convex domains would also be interesting. In this case, the blow-down analysis giving the asymptotic domain can be quite rich, and it is not clear which general class of domains would be most natural. We expect our analysis to be stable under appropriate small perturbations of smooth convex domains.

\subsubsection{Boundary conditions}

We note an interesting consequence of our analysis which holds for general convex domains: if the asymptotic domain $D_\infty$ {is empty}, the asymptotic behavior of the plasma is trivial in the sense that $\iint_{D\times\mathbb{R}^3}\mu^2(x,v,t)dxdv\to 0$ as $t\to\infty$ (see Corollary \ref{CorSupp}). This is intuitively clear: almost every particles will hit the boundary in finite time, at which point it will be absorbed into the wall and disappear.

This is a clear difference with other natural boundary conditions for the Vlasov field, such as specular reflexion, in which case the billiard dynamic is much more interesting. Such problems also occur in application and would be very interesting to consider. 

\subsection{A problem on the whole space}
{In this section, we reduce the problem \eqref{VPW} to a Cauchy problem.
The problem \eqref{VPW} has a characteristic boundary condition if $x\in \partial D$ and ${\bf n}\cdot v = 0$. We say a characteristic boundary condition if the normal vector to a boundary at a certain point is orthogonal to the direction vector of a characteristic curve. 
In general, handling characteristic boundary conditions can be difficult. Indeed, when $\lambda>0$, the solution of \eqref{VPW} may immediately become discontinuous even with $C^1$ class initial data. This does not happen when $\lambda<0$. Those fact can be observed by tracing the characteristic curves near the boundary with the aid of the following consequence of the Hopf maximum principle:
\begin{equation}\label{ExtPotAss2}
{\bf n}\cdot\nabla_x\phi(x)<0\hbox{ whenever }x\in\partial D.
\end{equation}
However, handling characteristic boundary conditions is cumbersome even when $\lambda<0$, and our approach for easier treatment is to replace the problem \eqref{VPW} with a Cauchy problem, where particles in $D^c:=\mathbb{R}^3\setminus D$ do not contribute. 
}
To this end, we use an extension map
\begin{equation*}
\begin{split}
\phi\in C_t(\mathbb R_{+};W^{2,\infty}_x(D)) \mapsto \psi\in C_t(\mathbb R_{+};W^{2,\infty}_x(\mathbb R^{3}))
\end{split}
\end{equation*}
such that the following properties (1)--(3) hold:
\begin{enumerate}
\item a linear relation in the domain
\begin{equation}\label{ExtPotAss1}
\begin{split}
\psi(x,t)=\phi(x,t):=\iint_{D\times\mathbb{R}^3} G_D(x,y)\mu^2(y,v,t)dydv\hbox{ for }x\in \overline{D},
\end{split}
\end{equation}
where $G_D$ denotes the Green function for the Dirichlet Laplacian on $D$;
\item an invariant domain condition in the exterior of the domain
\begin{equation}\label{ExtPotAss3}
\{v\cdot\nabla b(x),\mathcal{H}\}\le 0\hbox{ for }(x,v,t)\in D^c \times \mathbb{R}^{3} \times \mathbb{R}_{+} ,\qquad\hbox{where }\mathcal{H}:=\frac{\vert v\vert^2}{2}+\lambda\psi(x,t),
\end{equation}
where $\nabla b(x)$ is a suitable extension of the normal vector in $D^{c}$ such that $\nabla b(x)={\bf n}$ on $\partial D$ (for more details, see \eqref{DefB} and Lemma \ref{BProj}).
\item a bounded condition
\begin{align}
\|\nabla \psi(t)\|_{L^{\infty}_{x}(\mathbb{R}^3)} &\lesssim \|\nabla \phi(t)\|_{L^{\infty}_{x}(D)} + \varepsilon_1^{2} \langle t\rangle^{-2},
\label{ExtPotAss11}
\\
\|\nabla^{2} \psi(t)\|_{L^{\infty}_{x}(\mathbb{R}^3)} &\lesssim \|\nabla \phi(t)\|_{W^{1,\infty}_{x}(D)} +\varepsilon_1^{2} \langle t\rangle^{-2},
\label{ExtPotAss12}
\end{align}
where $\varepsilon_1 \in (0,1)$ is a constant to be determined according to the size of the initial data $\mu_{0}$. More precisely, $\varepsilon_1$ is the same as in \eqref{Boot1}, and eventually we choose $\varepsilon_{1}=10\varepsilon$ for $\varepsilon$ in Theorem~\ref{MainThm}.
\end{enumerate}
We claim that an extension map with \eqref{ExtPotAss1}--\eqref{ExtPotAss12} exists if $D$ is a $C^{3}$ convex domain with \eqref{QuantEG0}. We will prove the claim in Lemma \ref{ExtOperator} in Appendix \ref{S8}.
{In the proof of the claim, \eqref{ExtPotAss2} plays a key role}.
In the rest of the paper, we denote the force field by
\begin{gather}\label{DefF1}
F(x,t):=\nabla_x\psi(x,t), \quad E(x,t):=\nabla_x\phi(x,t),
\end{gather}
where $F(x,t)=E(x,t)$ holds when $x \in \overline{D}$.

From now on we focus on the following related Cauchy problem:
\begin{equation}\label{RVPW}
\begin{split}
\partial_t\mu+\{\mu,\mathcal{H}\}=0,\  \mathcal{H}:=\frac{\vert v\vert^2}{2}+\lambda\psi\hbox{ in }\mathbb{R}^{3}_{x}\times\mathbb{R}^{3}_{v}, 
\qquad {\mu}={\mu}_{0}\hbox{ at }t=0
\end{split}
\end{equation}
for $\mu:\mathbb{R}^3_x\times\mathbb{R}^3_v\times \mathbb R_{+} \to\mathbb{R}$. 
We note that $\phi$ satisfies \eqref{VPPhi} with $\widetilde{\mu}=\mathfrak{1}_{\{x\in D\}} \mu$.
Furthermore, the restriction $\mathfrak{1}_{\{x\in D\}}\mu$ of the solution $\mu$ of \eqref{RVPW} solves the original problem \eqref{VPW}. This is summarized in the next lemma.

\begin{lemma}\label{SuffCondBC}
Let $D$ be a convex domain with the extension map satisfiying \eqref{ExtPotAss1} and \eqref{ExtPotAss3}.
Assume that $\mu(x,v,t)$ is a continuous solution of the related Cauchy problem \eqref{RVPW} with the initial data $\mu_{0}$ such that $\mu_{0}=0$ on $\mathcal{S}_{in}:={\left\{(x,v) \, : \, x \in D^{c}, \ v\cdot \nabla b(x) \ge0  \right\}}$, then $\mathfrak{1}_{\mathcal{S}_{in}}\mu\equiv 0$ and $\mathfrak{1}_{\{x\in D\}}\mu(x,v,t)$ is  a solution of the original problem \eqref{VPW}.
\end{lemma}

Before proving Lemma \ref{SuffCondBC}, we state our main theorem for the related Cauchy problem \eqref{RVPW}.

\begin{theorem}\label{Thm2}
Let $D\subset\mathbb{R}^3$ be an acceptable convex domain (see Definition \ref{AcceptableDomains}). 
There exists  $\varepsilon_{0} > 0$ such that if the initial data $\mu_0\in C^1_{x,v}(\mathbb{R}^3\times\mathbb{R}^3)$ satisfies
\begin{equation}\label{IDD1}
\begin{split}
&\mu_{0}=0 \quad \hbox{on } \mathcal{S}_{in}={\left\{(x,v) \, : \, x \in D^{c}, \ v\cdot \nabla b(x) \ge0  \right\}},
\\
&\Vert \langle x\rangle^4\langle v\rangle^4\mu_0\Vert_{L^\infty_{x,v}}+\Vert \nabla_{x,v}\mu_0\Vert_{L^\infty_{x,v}}\le\varepsilon_0,
\end{split}
\end{equation}
then there exists a global strong solution { $\mu\in C_t([0,\infty);L^{2}_{x,v}(\mathbb{R}^3\times\mathbb{R}^3)) \cap L^{\infty}_t([0,\infty);W^{1,\infty}_{x,v}(\mathbb{R}^3\times\mathbb{R}^3))$} {with $ (1+|v|)^{-1} \partial_{t} \mu \in L^{\infty}_{t}([0,T];L^{\infty}_{x,v}(\mathbb{R}^3\times\mathbb{R}^3))$} of  the related Cauchy problem \eqref{RVPW}. It satisfies the bounds: for $t>0$,
\begin{equation}\label{GlobalBds2}
\begin{split}
\Vert \langle x-tv\rangle^4\langle v\rangle^4\mu(t)\Vert_{L^\infty_{x,v}}&\lesssim \varepsilon_0\langle \ln(2+t)\rangle^{8},\\
\Vert \nabla_{x}\mu(t)\Vert_{L^\infty_{x,v}}+\langle t\rangle^{-1}\Vert \nabla_{v}\mu(t)\Vert_{L^\infty_{x,v}}  &\lesssim \varepsilon_0\langle t\rangle^{\sqrt{\varepsilon_0}},\\
[1+t^2]\left[\Vert F(t)\Vert_{L^\infty_x}+\Vert\nabla_xF(t)\Vert\right]&\lesssim\varepsilon_0^2.
\end{split}
\end{equation}
Besides inside the domain, $\mu \in C^{1}_{x,v,t}(D \times \mathbb R^{3} \times \mathbb R_{+} )$ and also 
there hold the stronger bounds: for $x\in D$, $v\in\mathbb{R}^3$, $t>0$,
\begin{equation}\label{GlobalBds3}
\begin{split}
\vert \nabla_{x}\mu(x,v,t)\vert &\le 2\varepsilon_0,\\
\vert (t\nabla_{x}+ \nabla_{v})\mu(x,v,t)\vert&\le\varepsilon_0+C\varepsilon_0^3\langle \ln(2+t)\rangle^{10},\\
\left[1+\vert x\vert^2+t^2\right]\vert E(x,t)\vert&\le\varepsilon_0^2,\\
\left[1+t^2\right]^{3/2} \vert \nabla_{x} E(x,t)\vert&  \le\varepsilon_0^2 \langle \ln(2+t)\rangle^{9}.
\end{split}
\end{equation}
In addition, there {exists} an asymptotic domain $D_{\infty}$ and a limit distribution 
$\mu_\infty\in L^2_{x,v}\cap L^\infty_{x,v}\cap C^{0}_{x,v} (\mathbb{R}^3\times D_\infty)$ such that 
\begin{gather*}
\mu(x+tv+\lambda(\ln t)E_\infty(v),v,t)\rightarrow \mu_\infty(x,v) \quad a.e \ (x,v) \in \mathbb{R}^3 \times D_{\infty} \hbox{ as } t \to \infty,
\\
\Vert {\mu}\Vert_{L^2_{x,v}(D\times D_\infty^c)}\le\varepsilon \langle t\rangle^{-1/4} \hbox{ for } t>0,
\end{gather*}
where the asymptotic electric field $E_{\infty}$ is given by the distribution $\mu_\infty$ and  the Dirichlet Green function $G_\infty$ of the domain $D_\infty$ as
\begin{equation*}
\begin{split}
E_{\infty}^{j}(a)=\iint_{\mathbb{R}^3 \times D_{\infty}}\partial_{a^j} G_{\infty}(a,w)\mu^2_{\infty}(z,w)dzdw \quad \hbox{for }j=1,2,3.
\end{split}
\end{equation*}
\end{theorem}

\begin{remark}
With our extension, we only have $\nabla^2\psi\in L^\infty$, but it will be discontinuous. As a result, we are only able to obtain strong (Lipshitz) solutions of the related Cauchy problem \eqref{RVPW}. On the other hand, we can recover classical solutions inside the domain $D$ as stated in Theorem \ref{Thm2}.
\end{remark}

Theorems \ref{MainThm} and \ref{Thm2} will be proved in Section \ref{S7}. In order to illustrate the use of the Poisson bracket, we present the simple proof of Lemma \ref{SuffCondBC}.

\begin{proof}[Proof of Lemma \ref{SuffCondBC}]

We compute that for $\theta_1,\theta_2>0$ and $\chi_{>,<}$ defined in \eqref{Partition1},
\begin{equation*}
\begin{split}
&\frac{d}{dt}\iint \mu^2(x,v,t)\chi_<(\theta_1^{-1}b(x))\chi_{>}(\theta_2^{-1}v\cdot\nabla b(x))dxdv\\
&=\iint\mu^2\{\chi_<(\theta_1^{-1}b(x))\chi_{>}(\theta_2^{-1}v\cdot\nabla b(x)),\mathcal{H}\}dxdv\\
&=\iint\mu^2\Big[\theta_1^{-1}(v\cdot\nabla b(x))\chi^\prime_<(\theta_1^{-1}b(x))\chi_>(\theta_2^{-1}v\cdot\nabla b(x))\\
&\qquad +\theta_2^{-1}\chi_<(\theta_1^{-1}b(x))\chi^\prime_{>}(\theta_2^{-1}v\cdot\nabla b(x))\{v\cdot\nabla_xb(x),\mathcal{H}\}\Big]dxdv\le 0,
\end{split}
\end{equation*}
where we have used \eqref{ExtPotAss3}, $\chi^\prime_< \leq 0$, and $\chi_<, \chi_>, \chi^\prime_{>} \geq 0$ in deriving the inequality. Integrating this and letting $\theta_1,\theta_2\to0$, we see that $\mu(x,v,t)=0$ on $\hbox{Int} \mathcal{S}_{in}$ for any $t>0$. Since $\mu\in C^0_{x,v,t}$, this fact means that $\mu=0$ for all times on $\{x\in\partial D,\,\,v\cdot{\bf n}\ge 0\}$.  Thus $\mathfrak{1}_{\{x\in D\}}\mu $ is a solution of \eqref{VPW}.
\end{proof}

\subsection{Open problems}

{
We now turn our attention to several interesting open problems. We first outline three main research directions.
The first direction involves investigating the asymptotic behavior of solutions in greater detail using higher-order expansions. The second is to change the domain and the boundary conditions, and study its effect on the asymptotic behavior of solutions. 
Finally, we will look at applying our technique to plasma sheaths.


It is worth noting that, for non-convex domains, propagating the regularity, even the continuity along a grazing trajectory seems very challenging.

The open problems can be summarized in a rough order of expected difficulty:
}

\begin{enumerate}
\item {Find a more precise asymptotic behavior of modified scattering via the expansion method in \cite{BR2024,ST2024}.}

\item Extend the present analysis to the case where the boundary is perturbed by a small, rapidly decaying function, making the domain non-convex, but preserving the smoothness and the nature of the asymptotic domain $D_\infty$.

\item Consider the cases of exterior of smooth convex domains. 

\item Consider the cases of general convex domains, in particular a corner $\Omega:=\{x \in \mathbb R^{3} \, : \,  x_{1}>0, \ x_{2}>0 \}$. 
\item Consider the case of a general wall, curved and when $D$ is not convex.
\item Consider other forms of boundary conditions. A natural choice is the case of specular reflexion for the particles on the wall, and insulating condition for the electric field:
\begin{equation*}
\begin{split}
\mu(x,v,t)={\mu(x,v-2(v\cdot{\bf n})\cdot{\bf n},t)},\qquad {\bf n}\cdot\nabla_x\phi(x,t)=0\qquad \hbox{when }x\in \partial D.
\end{split}
\end{equation*}

\item 
Consider the stability of the plasma sheath. In \cite{STZ1}, the stability is conditionally shown.
It is assumed that the support of initial perturbations is located in $\{v^{1} <0 \}$, where $v_{1}$ is the vertical velocity for the boundary and the particles with $v_{1}<0$ goes to the boundary.
It would be interesting to study the stability for more general initial perturbations by applying the techniques in this paper.

\end{enumerate}

\section{Notations and hypotheses}\label{S2}

\subsection{Notation}

We define a partition of unity subordinate to the covers $\mathbb{R}=(-\infty,-1)\cup(-2,2)\cup(1,\infty)$, and $\mathbb{R}_+\setminus\{0\}=\cup_{k\in\mathbb{Z}}(2^{k-1},2^{k+1})$. Namely, let $\varphi\in C^\infty_c(1/4,4)$ such that $\sum_{k\in\mathbb{Z}}\varphi(2^{-k}x)\equiv 1$ a.e. We also let
\begin{gather*}
\varphi_k(x):=\varphi(2^{-k}x),\qquad\varphi_{\le 0}:=\sum_{k\le 0}\varphi_k(x),\qquad \varphi_{>0}(x):=\sum_{k>0}\varphi_k(x),
\\
{ 
\varphi_{\le L}(x):=\sum_{k \leq L}\varphi_k(x), \quad \varphi_{\geq L}(x):=\sum_{k \geq L}\varphi_k(x) \quad  \hbox{  for } L \in \mathbb R
}
\end{gather*}
and set 
\begin{equation}\label{Partition1}
\begin{split}
\chi_>(x):=\varphi_{>0}(x),\qquad\chi_{<}(x):=\chi_{>}(-x),\qquad\chi_0(x):=\varphi_{\le0}(\vert x\vert) \quad \hbox{for } x\in\mathbb{R}.
\end{split}
\end{equation}

\subsection{Domain hypothesis}\label{S22}

Let us list the assumptions on the domain and their consequences. 

\subsubsection{Domain geometry}\label{221}
Let $D$ be a convex domain (not necessarily  $C^{3}$). Without loss of generality, we assume that the origin is on the wall $0\in\partial D$, and the wall is differentiable with normal unit vector ${\bf n}=(0,0,1)$ at the origin. We define the extended distance
\begin{equation}\label{DefB}
\begin{split}
b(x)&:=\begin{cases}-d(x,D)&\quad\hbox{ when } x\in D^c\\
d(x,D^c)&\quad\hbox{ when }x\in D,
\end{cases}
\end{split}
\end{equation}
where $d(\cdot,\cdot)$ denotes a standard distance function. 
{The function $b$ satisfies $b(0)=0$, $\nabla b(0)={\bf e}_3$, and the properties in Lemma \ref{BProj}}.

Since $D$ is convex, we can define the projection $\pi_b:\mathbb{R}^3\to \overline{D}$, and we define the reflection of a point $x\in D^c$ across $\partial D$ as
\begin{equation}\label{DefRho}
\begin{split}
D^c\ni x\mapsto \rho(x):=2\pi_b(x)-x.
\end{split}
\end{equation}
For the case of acceptable convex domains, we can see more detailed properties on the domain geometry.
See Appendix \ref{A2}.

\subsubsection{Asymptotic domain}\label{SecDInfty}
We also let $x_0=(0,0,1)$ be a reference point in the domain: $b(x_0)>0$.  We define the gauge of $D$ at $x\in D$ to be
\begin{equation}\label{j1}
\begin{split}
j_x(v):=\inf\{r>0 \, : \, x+ r^{-1}v\in D\}.
\end{split}
\end{equation}
This is also (the inverse of) the time of escape of $D$ for $(x,v)$. Note that $j_x$ is subadditive and positively homogeneous. In the following, we fix an arbitrary point $x_0\in D$, and sometimes abbreviate $j_{x_0}(v)$ as $j(v)$.

The asymptotic domain can be defined in several different ways. We note that $j_{x}^{-1}(\{0\})$ is independent of $x\in D$, and we define $D_\infty=\hbox{Int}(j_{x_0}^{-1}(\{0\}))$ to be its interior. In the following, we will assume that $D_\infty\ne\emptyset$ (else the asymptotic analysis is trivial and we simply define $E_\infty=0$).

If $D_\infty\ne\emptyset$, we note that $\overline{D_\infty}=j_{x_0}^{-1}(\{0\})$ is the maximal (solid) convex cone contained in $\overline{D}$ with vertex at the origin (see also \eqref{DInfty} for another definition). Let $C_{\infty}:=\partial D_\infty$ be its boundary. In the rest of the paper, we call $D_\infty$ the {\it asymptotic domain} and $C_{\infty}$ the {\it asymptotic cone}.

The wall can be represented as a the graph of a function $C_{\infty}\ni x\mapsto {\bf w}_{C_{\infty}}(x)\in \partial D$, whose radial derivative is decreasing almost everywhere. We set
\begin{equation*}
\begin{split}
\mathfrak{W}(r):=\sup_{x\in C_{\infty},r/2\le \vert x\vert\le 2r}\vert{\bf w}_{C_{\infty}}(x)-x\vert,
\end{split}
\end{equation*}
and then assume that
\begin{equation}\label{UCAssumption}
\begin{split}
 \mathfrak{W}(r)\lesssim 1,
\end{split}
\end{equation}
{which means the distance between $\partial D$ and $C_\infty$ is almost constant at large distances from the origin. For instance, $C^{3}$ convex domains such as hyperbolic cylinders and hyperboloids, as shown in Figure \ref{Fig1}, satisfy this condition. Appendix \ref{hyper} will demonstrate that hyperbolic cylinders satisfy the condition \eqref{UCAssumption}.}
Since $D_\infty$ is convex, we can define $\pi_\infty$ the projection onto ${\overline{D_\infty}}$. 
 
\subsubsection{Green functions}\label{223}

For a convex domain, we have the following bounds for the Dirichlet Green function: for $x,y \in D$,
\begin{equation}\label{QuantEGC}
\begin{split}
\vert x-y\vert\cdot \vert G(x,y)\vert&\lesssim \min\{1,\vert b(x)\vert\vert x-y\vert^{-1}, b(y)\vert x-y\vert^{-1},b(x)b(y)\vert x-y\vert^{-2}\},\\
\vert x-y\vert^2\cdot\vert \nabla_xG(x,y)\vert&\lesssim \min\{1,b(y)\vert x-y\vert^{-1}\},\\
\vert x-y\vert^2\cdot\vert \nabla_yG(x,y)\vert&\lesssim \min\{1,\vert b(x)\vert\vert x-y\vert^{-1}\},\\
\vert x-y\vert^3\cdot\vert\nabla_{x}\nabla_yG(x,y)\vert&\lesssim 1.
\end{split}
\end{equation}
These are shown in \cite[Proposition 2.3]{FrThesis}. For $t>0$, we introduce the rescaled domain $D_t:=\{x \, : \, tx\in D\}$, which is a decreasing sequence:
\begin{equation}\label{DInfty}
 D_\infty=\cap_{t>0}D_t\subset D_t\subset D_0=\cup_{t>0}D_t=\{x^3>0\}.
 \end{equation}
We will also use rescalings of the Green function and we set for $t>0$
\begin{equation*}
\begin{split}
G_t(x,y):=tG(tx,ty).
\end{split}
\end{equation*}
We note that $G_t(x,y)$ is the Green function of $D_t$ when $x,y\in D_t$. Indeed, letting $\mathcal{G}_t$ be the Green function of $D_t$,  we see that the function $x\mapsto H_t(x,y):=G_t(x,y)-\mathcal{G}_t(x,y)$ is harmonic and vanishes on the boundary and thus $\mathcal{G}_t=G_{t}$ holds. 
The family of {Green} functions $G_t$ is increasing in the sense that if $0\le s\le t\le\infty$ and $x\ne y\in D_t$, then
\begin{equation}\label{CompGF0}
G_0(x,y)\le G_s(x,y)\le G_t(x,y)\le 0,
\end{equation}
since $x\mapsto H_{s,t}(x,y):=G_s(x,y)-G_t(x,y)$ is harmonic on $D_t$ and negative on the boundary. 


\subsubsection{Acceptable domain}\label{AD}
All the assumptions in subsections \ref{221}--\ref{223} are \eqref{UCAssumption} and the convexity of $D$.
We will need the following stronger analytical properties of the Green function. 
\begin{equation}\label{QuantEG0}
\vert x-y\vert^3\cdot[\vert\nabla_{x}^2G(x,y) \vert+\vert\nabla_y^2G(x,y)\vert]\lesssim 1.
\end{equation}
\begin{equation}\label{T1Hyp}
 \left| \int_{R=0}^{1}  \nabla^2u_R(x)  \frac{dR}{R} \right| \lesssim 1, \quad u_R(x):=\int_D\varphi(R^{-1}|x-y|)G(x,y)dy,
\end{equation}
\begin{equation}\label{AssumptionExtensionEF}
\begin{split}
\sup_{x\in D,\, b(x)\le 1}\Vert [\nabla_xb(x)\cdot \nabla_xG(x,y)]_+\Vert_{L^1_y(D\cap\{\vert x-y\vert\le L\})}&\lesssim b(x) \ln (L+1),\\
 \sup_{x\in D,\, b(x)\le 1}\Vert [\nabla_xb(x)\cdot \nabla_xG(x,y)]_+\Vert_{L^\infty_y(D\cap\{\vert x-y\vert\ge L\})}&\lesssim b(x)L^{-3} \hbox{ for } L\ge1.
\end{split}
\end{equation}
\begin{definition}\label{AcceptableDomains}
We say that $D$ is an acceptable convex domain if $D$ is a $C^{3}$ convex domain with \eqref{UCAssumption} and  \eqref{QuantEG0}. 
\end{definition}

\begin{remark}
We note that the inequalities \eqref{QuantEG0}--\eqref{AssumptionExtensionEF} play a key role in the proof of our main theorem. 
Using an adapted modification of the double layer potential, we will show in Appendix \ref{SmoothDomain} that the Dirichlet Green function satisfies \eqref{T1Hyp} for all $C^{3}$ domains and \eqref{AssumptionExtensionEF} for all $C^{2,\alpha}$ domains with \eqref{QuantEG0}.
We also note that the half-space case is a fundamental example of the acceptable domain. In this case, it is clear that the assumption \eqref{UCAssumption} holds. We will confirm in Appendix \ref{ExFlatWall} that the inequalities \eqref{QuantEG0}--\eqref{AssumptionExtensionEF}  hold for the case of a half-space.

{
Hyperbolic cylinder domains of the form $\{(x_1,x_2,x_3)\, :\, x^1x^2>1, \, x^{1}>0\}$ are also admissible convex domains. In Appendix \ref{hyper}, we show that they satisfy \eqref{UCAssumption} and \eqref{QuantEG0}.
To the best of our knowledge, it remains an open question to determine the $C^3$ convex domains for which \eqref{QuantEG0} holds, although it is known to be valid for bounded domains.
}

\end{remark}

\subsection{{Gauge and grazing functions}}\label{JG}

We discuss the properties of the gauge function $j$ defined in \eqref{j1} and the grazing function $g$ defined as follows.
Let $\mathcal{I}(x,v)=\sup\{t\in\mathbb{R} : b(x+tv)\le 0\}$ be the impact time of $(x,v)$. We employ the grazing function:
\begin{equation}\label{DefGrazing}
\begin{split}
g(x,v):=\inf_{s\le \mathcal{I}(x,v)}v\cdot\nabla b(x+sv).
\end{split}
\end{equation}
Since $D_\infty$ is open, we have that $\mathcal{I}(x,v)<+\infty$ whenever $v\in D_\infty$. In this case, using \eqref{BNondegenerate},  we can write $g(x,v)=v\cdot{\bf n}(y)$ for a certain $y=x+\mathcal{I}(x,v)v \in \partial D\cap (x+\mathbb{R}v)$ and in particular $g(x,v)>0$ (see Figure \ref{Fig1}).

For any $v\in\{j_{x_0}>0\}$, we see that $z:=x_0+(1/j_{x_0}(v))v\in \partial D$, and so we can define {$\nu(v)\in D b(z)$} to be an inward unit normal vector, as well as $\alpha(v)=z\cdot\nu$, $\theta(v)=v\cdot\nu(v)<0$, and $\mathcal{H}_{t,v}:=\{x \, : \, x\cdot\nu(v)>t\}$. 
It is clear that\footnote{This and \eqref{DInfty} illustrate the role of the flat wall as a model case.}
\begin{equation*}
\begin{split}
D=\cap_{\{v\, : \, j_{x_0}(v)>0\}}\mathcal{H}_{\alpha(v),v}.
\end{split}
\end{equation*}
This allows to control the trajectories which stay in $D$ for a long time. It is summarized in Lemma \ref{EscapeLemma}. 
%
%
%
In addition, the following quantitative bound on the vanishing of $j$ and $g$ hold {for all $\alpha\in(0,\alpha_0)$ with some $\alpha_0 >0$}:
\begin{equation}\label{QuanfVanishingG}
\begin{split}
{\sup_{z\in\mathbb{R}^3}\vert J_\alpha(z) \vert\lesssim \alpha},
\end{split}
\end{equation}
where
\begin{equation*}
\begin{split}
J_\alpha(z):= \{ {w\in \mathbb{S}^{2}}  \, : \, z\in D^c,\, {0<}g(z,w)<\alpha\}&\cup\{{w\in \mathbb{S}^{2}} \, : \, z\in D,\, {0<j(w)<\alpha}\}.
\end{split}
\end{equation*}
{
This will be shown in Lemma \ref{GrazingLemma}.
The condition ensures that the measure of the velocity $w$ of particles that graze the boundary decreases as the grazing angle becomes smaller, and that the measure of the velocity $w$ of particles exiting the domain $D$ decreases over time if $w \not\in D_\infty$.}


\section{{Strategy and Plan}}\label{S3}

{Using a standard bootstrap argument, we will prove Theorem \ref{Thm2}, which combined with Lemma \ref{SuffCondBC} immediately yields Theorem \ref{MainThm}.
More precisely, we assume a bound on solutions over a certain time interval and then show a slightly improved bound. This enable us to establish the existence of global solutions and show the estimates \eqref{GlobalBds2}--\eqref{GlobalBds3} in Theorem \ref{Thm2}. The keys of our bootstrap argument are Propositions \ref{ProprBoot1} and \ref{ProprBoot2} below.
In subsections \ref{S3.1} and \ref{S3.2}, we discuss these propositions. 
Subsection \ref{S3.3} provides a renormalization of $\mu$, which is used throughout the proof.
Subsection \ref{S3.4} presents an outline of this paper.}

\subsection{Bootstrap 1: control of  moments}\label{S3.1}

First, assuming only {the bound on the moment $ \langle x-tv\rangle^4 \langle v\rangle^4\mu$}, we obtain global {bounds} of weak solutions and decay of the electric unknowns. Our first bootstrap assumption is as follows: for $0 \leq t \leq T$,
\begin{equation}\label{Boot1}
\begin{split}
\Vert \langle x-tv\rangle^4 \langle v\rangle^4\mu(t)\Vert_{L^\infty_{x,v}}&\le\varepsilon_1\langle \ln(2+t)\rangle^{10},\\
\left[1+t^2\right]\Vert F(t)\Vert_{L^{\infty}_{x}}&\le\varepsilon_1^2 \langle t\rangle^{\frac{1}{10}},
\end{split}
\end{equation}
where $F$ is defined in \eqref{DefF1}.
We claim that this bootstrap can be closed.

\begin{proposition}\label{ProprBoot1}
Let $D$ be a convex domain with \eqref{ExtPotAss3}, \eqref{ExtPotAss11}, \eqref{UCAssumption}, \eqref{QuanfVanishingG}, and \eqref{QuantEG0}.
Assume that $\mu$ solves \eqref{RVPW} in the weak sense for $0\le t\le T$, and also satisfies the bootstrap assumptions \eqref{Boot1} and
the initial bounds
\begin{equation}\label{ID1}
\begin{split}
&\Vert \langle x\rangle^4\langle v\rangle^4\mu(t=0)\Vert_{L^\infty_{x,v}}\le\varepsilon_0,
\\
&\mu(t=0)=0 \quad \hbox{on } \mathcal{S}_{in}. 
\end{split}
\end{equation}
Then $\mu$ satisfies the stronger bounds for $0 \leq t \leq T$:
\begin{equation}\label{Boot1CCl}
\begin{split}
\Vert \langle x-tv\rangle^4\langle v\rangle^4\mu(t)\Vert_{L^\infty_{x,v}}&\le\varepsilon_0+C\varepsilon_0\varepsilon_1\langle \ln(2+t)\rangle^{8},\\
\left[1+t^2\right]\Vert F(t)\Vert_{L^\infty_{x}}&\lesssim \varepsilon_1^2\langle \ln(2+t)\rangle.
\end{split}
\end{equation}
Furthermore, inside the domain $D$, there holds that for $x\in D$ and $0 \leq t \leq T$,
\begin{equation}\label{Boot1CCl2}
\begin{split}
\sqrt{1+t^2+\vert x\vert^2}\vert\phi(x,t)\vert&\lesssim \varepsilon_1^2,\\
\left[1+t^2+\vert x\vert^2\right]\vert E(x,t)\vert&\lesssim \varepsilon_1^2\langle \ln(2+t)\rangle.
\end{split}
\end{equation}

\end{proposition}

\begin{remark}

We observe that this limited bootstrap is enough to show, in a large generality (of convex domain and initial data) that weak solutions with finite moments have decreasing electrostatic potential and field. 
For instance, the corner $\Omega:=\{x \in \mathbb R^{3} \, : \,  x_{1}>0, \ x_{2}>0 \}$ satisfies all the conditions \eqref{ExtPotAss3}, \eqref{ExtPotAss11}, \eqref{UCAssumption}, \eqref{QuanfVanishingG}, and \eqref{QuantEG0}.
On the other hand, in order to upgrade this analysis to strong solutions and asymptotic behavior, we also need control of derivatives, and hence the stronger bootstrap \eqref{Boot2}. This comes at the cost of more stringent assumptions on the domain.

\end{remark}

Based only on control of moments, we can ascertain that the asymptotic data is supported on $v\in \overline{D_\infty}$.

\begin{corollary}\label{CorSupp}
Under the same assumptions as in Proposition \ref{ProprBoot1}, there holds that
\begin{equation*}
\begin{split}
\iint_{{D\times (\overline{D_{\infty}})^c }}\mu^2(x,v,t)dxdv\lesssim \varepsilon_{0 }^2\langle t\rangle^{-3/4}\langle \ln (2+t)\rangle^{20}.
\end{split}
\end{equation*}
\end{corollary}
\begin{proof}
We recall that $v\in (\overline{D_{\infty}})^c$ is equivalent to $j(v)>0$.
For any $\kappa>0$, we observe by using Lemma \ref{EscapeLemma} that for $t\ge C\kappa^{-1}$,
\begin{align*}
\iint_{\{x\in D,\,\, j(v)>\kappa\}}\mu^2(x,v,t)dxdv 
& \leq \iint_{\{x\in D,\,\, j(v)>\kappa\}}\langle x-tv\rangle^{-4} \langle v\rangle^{-4}  dxdv \cdot\Vert \langle x-tv\rangle^2 \langle v\rangle^2 \mu(t)\Vert_{L^\infty_{x,v}}^{2}
\\
& \lesssim \varepsilon_0^2\iint_{\{\vert y\vert\gtrsim \kappa t \}} \langle y\rangle^{-4}  \langle v\rangle^{-4} dxdv \langle \ln (2+t)\rangle^{20} 
\lesssim \varepsilon_{0}^2 \langle \kappa t\rangle^{-3}\langle \ln (2+t)\rangle^{20}.
\end{align*}
On the other hand, passing to polar coordinates in $v$, we get
\begin{align*}
\iint_{\{x\in D,\,\, 0<j(v)<\kappa\}}\mu^2(x,v,t)dxdv 
& \leq \iint_{\{0<j(v)<\kappa\}} \langle v\rangle^{-4} \langle x-tv\rangle^{-4} dxdv\cdot \Vert \langle v\rangle^2 \langle x-tv\rangle^2\mu(t)\Vert_{L^\infty_{x,v}}^{2}
\\
&\lesssim \varepsilon_{0}^2\cdot \vert\{\vert v\vert=1,\,\, 0<j(v)<\kappa\}\vert \langle \ln (2+t)\rangle^{20} 
\lesssim \varepsilon_0^2\kappa\langle \ln (2+t)\rangle^{20},
\end{align*}
where we have used \eqref{QuanfVanishingG}. Choosing $\kappa=\langle t\rangle^{-3/4}$, we obtain the desired result.
\end{proof}

\subsection{Bootstrap 2: control of derivatives and asymptotic behavior}\label{S3.2}
{Assuming {the bound on the derivatives $\nabla_{x}\mu$ and $(t\nabla_{x} + \nabla_{v})\mu$} as well, we obtain global {bounds} of strong solutions, and improved decay of $E$.} 
Our second bootstrap assumption is as follows: for $x,v \in \mathbb R^{3}$ and $0 \leq t \leq T$,
\begin{equation}\label{Boot2}
\begin{split}
\vert \nabla_{x} \mu(x,v,t)\vert&\le\varepsilon_1[1+\mathfrak{1}_{\{b(x)\le0\}}\langle t\rangle^{\delta_{1}}],\\
\vert (t\nabla_{x} + \nabla_{v})\mu(x,v,t)\vert&\le\varepsilon_1 \langle t\rangle^{\delta_{1}}[1+\langle t\rangle\mathfrak{1}_{\{b(x)\le0\}}],\\
\left[1+t^2\right]^{3/2} \Vert \nabla_{x} E(t)\Vert_{L^\infty_x(D)}&  \lesssim \varepsilon_1^2\langle t \rangle^{\frac{1}{10}},
\end{split}
\end{equation}
where $\delta_{1} \in (0,\frac{1}{10})$. 

\begin{proposition}\label{ProprBoot2}
Let $D\subset\mathbb{R}^3$ be an acceptable convex domain (see Definition \ref{AcceptableDomains}). 
Assume that $\mu$ solves \eqref{RVPW} for $0\le t\le T$ in the strong sense, and satisfies the bootstrap assumptions \eqref{Boot1} and \eqref{Boot2} as well as the initial bounds \eqref{ID1} and 
\begin{equation}\label{IDD1}
\begin{split}
\Vert \nabla_{x,v}\mu(t=0)\Vert_{L^\infty_{x,v}}&\le\varepsilon_0.
\end{split}
\end{equation}
Then $\mu$ satisfies the stronger bounds for $0 \leq t \leq T$:
\begin{equation}\label{Boot2CCl}
\begin{split}
\Vert \nabla_{x}\mu(t)\Vert_{L^\infty_{x,v}}&\le\varepsilon_0+C\delta_{1}^{-1}\varepsilon_1^2\langle t\rangle^{\delta_{1}},\\
\Vert (t\nabla_{x} + \nabla_{v})\mu(t)\Vert_{L^\infty_{x,v}} &\le\varepsilon_0+C\varepsilon_1^2\langle t\rangle^{\delta_{1}+1},\\
\left[1+t^2\right]\cdot \Vert  \nabla F(t)\Vert_{L^{\infty}_x}&\lesssim \varepsilon_1^2.
\end{split}
\end{equation}
Furthermore, inside the domain $D$, there holds that for $x\in D$, $v \in \mathbb R^{3}$, and $0 \leq t \leq T$,
\begin{equation}\label{Boot2CClD}
\begin{split}
\vert\nabla_{x}\mu(x,v,t)\vert &\le\varepsilon_0+C\varepsilon_0\varepsilon_1^2,\\
\vert (t\nabla_{x} \mu + \nabla_{v}\mu)(x,v,t)\vert&\le\varepsilon_0+C\varepsilon_0\varepsilon_1^2\langle \ln(2+t)\rangle^{12},\\
\left[1+t^2+\vert x\vert^2\right]\vert E(x,t)\vert&\lesssim\varepsilon_1^2,\\
\left[1+t^2\right]^{3/2} \vert \nabla E(x,t)\vert&  \lesssim \varepsilon_1^2 \langle \ln(2+t)\rangle^{11}.
\end{split}
\end{equation}
\end{proposition}

\subsection{Renormalization}\label{S3.3}

In order to obtain a fixed measure, we conjugate by the linear flow and define
\begin{equation*}
\begin{split}
\gamma(x,v,t):=\mu(x+tv,v,t),
\end{split}
\end{equation*}
which satisfies
\begin{equation}\label{RenormalizedVP}
\begin{split}
\partial_t\gamma+\{\gamma,\Phi\}=0,\qquad\Phi(x,v,t):=\lambda \psi(x+tv,t).
\end{split}
\end{equation}
This is the equation we will work on for most of the analysis.

\subsection{{Outline of this paper}}\label{S3.4}

This paper is organized as follows. In Section \ref{S4}, we show the decay of the electric field $E$ in \eqref{Boot1CCl2} and \eqref{Boot2CClD}. The decay of $F$ in \eqref{Boot1CCl} and \eqref{Boot2CCl} immediately follows from \eqref{ExtPotAss11}--\eqref{ExtPotAss12}, \eqref{DefF1}, and the decay of $E$. Section \ref{S5} is devoted to studying the bounds on $\mu$. Specifically, Subsections \ref{S5.1} and \ref{S5.2} provide the proofs of the bounds of the moments $\langle x-tv\rangle^4 \langle v\rangle^4\mu$ and the derivatives $\nabla_{x}\mu$, $(t\nabla_{x} + \nabla_{v})\mu$ in \eqref{Boot1CCl} and \eqref{Boot2CCl}--\eqref{Boot2CClD}, respectively.  Subsection \ref{S5.3} completes the proofs of Propositions \ref{ProprBoot1} and \ref{ProprBoot2}. In Section \ref{S6}, we analyze the asymptotic behavior of solutions by finding $\mu_\infty$ in Theorem \ref{Thm2}. It is worth noting that the decay of the $L^2$-norm of $\mu$ in Theorem \ref{Thm2} follows from Corollary \ref{CorSupp}. In Section \ref{S7}, we complete the proof of Theorem \ref{Thm2}. First, we show that $\mu$ is $C^1$ inside $D$. Then, we establish the existence of a global solution combining the time-local solvability and the bootstrap argument with Propositions \ref{ProprBoot1} and \ref{ProprBoot2}. Furthermore, we show that Theorem \ref{MainThm} is an immediate consequence of Theorem \ref{Thm2}.

In Appendix \ref{SD}, we construct the extension from $\phi$ to $\psi$. 
Appendix \ref{A2} details the geometric properties of the domain $D$. 
In Appendix \ref{ExFlatWall}, we demonstrate that a half-space is an admissible convex domain.
In Appendix \ref{hyper}, we show that hyperbolic cylinders satisfy \eqref{UCAssumption} and \eqref{QuantEG0}.
In Appendix \ref{ExCorner1}, we briefly discuss a counterexample of acceptable convex domain.
Finally, Appendix \ref{SmoothDomain} analyzes the properties of the Green function, which are then used in Appendix \ref{A.7} to prove \eqref{T1Hyp} and \eqref{AssumptionExtensionEF}.

\begin{remark}
If $D$ is a half space $\mathbb R^3_+:=\{x^3>0\}$, the proof can be simplified. 
For the reader's convenience, we explain where these simplifications occur.
First the admissible domain $D_\infty$ and the functions $b$, $j$, and $g$ can be explicitly written as
\[
D_\infty=\mathbb R^3_+, \quad b(x)=x^{3}, \quad j(v)=-v^{3}_{-}, \quad g(x,v)=v^{3}.
\]
Additionally, we can omit the estimates \eqref{VanishAwayFromDInfty}--\eqref{E0E1} in the proof of Proposition \ref{ControlDiffEProp}, since $D=D_{t}=D_{\infty}$, $G=G_{t}=G_{\infty}$, and $E^0=E^1$ hold.
The rest of the proof  is largely the same, with the exception that $b$, $j$, and $g$ have explicit forms.
Indeed, the proofs in Sections \ref{S4}--\ref{S7} does not change much. 
In Appendix \ref{SD}, it becomes easier a bit to follow the proof of Lemma \ref{ExtOperator} if we use $b(x)=x^{3}$.
We can omit the discussion in Appendix \ref{A2}.
In Appendix \ref{ExFlatWall}, we show not only the fact that the half-space $\mathbb R^3_+$ is an admissible convex domain, but also the assertions in Lemmas \ref{EscapeLemma}--\ref{GrazingLemma} in Appendix \ref{A2} and the estimates \eqref{T1Hyp} and \eqref{AssumptionExtensionEF} for this case.
Appendices \ref{hyper}--\ref{A.7} are not needed.
\end{remark}

\section{Control of the Electric field}\label{S4}

We now propagate the control of the electric field $E$. This will require controls on the density $\gamma$. We set
\begin{equation*}
\begin{split}
N_1&:=
\Vert \langle z\rangle^4\langle w\rangle^4\gamma\Vert_{L^\infty_{z,w}}^2.
\end{split}
\end{equation*}

\subsection{Bounds on electric unknowns}

\subsubsection{Decomposition along scales}

We decompose the electric potential and field along scales instead of using the level sets of the convolution kernel:

\begin{align}
{\phi}(x,t)&=\int_{R=0}^\infty {\phi}_{R}(x,t)\frac{dR}{R^2}, \ \
{\phi}_{R}(x,t):= R\iint  \mathfrak{1}_{\{b(y)>0\}} \varphi(R^{-1}|x-y|)G(x,y)\mu^2(y,w,t)dydw,
\label{DecPhiScale}\\
{E}(x,t)&=\int_{R=0}^\infty {E}_{R}(x,t)\frac{dR}{R^3}, \ \
{E}^j_{R}(x,t):= R^2\iint \mathfrak{1}_{\{b(y)>0\}} \partial_{x^j}\left\{\varphi(R^{-1}|x-y|)G(x,y)\right\}\mu^2(y,w,t)dydw.
\label{DecEScale}
\end{align}

Using the linear flow factorization, we can rewrite ${\phi}_{R}$ and its approximation $\phi^0_{r}(a,t)\simeq \phi_{tr}(ta,t)$:
\begin{equation*}
\begin{split}
 {\phi}_{R}(x,t)&= R\iint  \mathfrak{1}_{\{b(z+tw)>0\}} \varphi(R^{-1}|x-z-tw|)G(x,z+tw)\gamma^2(z,w,t)dzdw,\\
  \phi^0_{r}(a,t)&:= r\iint  \mathfrak{1}_{\{w\in D_\infty\}} \varphi(r^{-1}|a-w|)G_\infty(a,w)\gamma^2(z,w,t)dzdw.
 \end{split}
\end{equation*}
Similarly,
\begin{equation*}
\begin{split}
{E}^j_{R}(x,t)&= R^2\iint \mathfrak{1}_{\{b(z+tw)>0\}}  \partial_{x^j}\left[\varphi(R^{-1}|x-z-tw|)G(x,z+tw)\right]\gamma^2(z,w,t)dzdw,\\
E^{0,j}_{r}(a,t)&:= r^2\iint \mathfrak{1}_{\{w\in D_\infty\}}  \partial_{a^j}\left[\varphi(r^{-1}|a-w|)G_\infty(a,w)\right]\gamma^2(z,w,t)dzdw.
\end{split}
\end{equation*}
We also define the effective unknowns $\phi^0(a,t)$ and $E^{0,j}(a,t)$ by
\begin{equation*}
\begin{split}
\phi^{0}(a,t):=\frac{1}{t}\int_{t=0}^{\infty} \phi^{0}_{r}(a,t) \frac{dr}{r^{2}}
,\quad
E^{0,j}(a,t):=\frac{1}{t^2}\int_{t=0}^{\infty} E^{0,j}_{r}(a,t) \frac{dr}{r^{3}},
\end{split}
\end{equation*}
{which represent the leading-order asymptotics of $\phi$ and $E$ as $t \to \infty$. }
W remark that the decompositions along scales are relevant to prove almost optimal bounds for $\phi$ and $E$ as well as improved bounds for $\vert \phi-\phi^0\vert$ and $\vert E-E^0\vert$.

\subsubsection{{Crude estimates}}

We derive crude estimates for $\phi_{R}$, $E_{R}$, $\phi^{0}_{r}$, and $E^{0}_{r}$.

\begin{lemma}\label{crude}
Under the same assumption as in Proposition \ref{ProprBoot1}, there holds that
\begin{equation}\label{ScaleEstDiffPhi}
\begin{split}
\vert {\phi}_{R}(x,t)\vert+\vert {E}_{R}(x,t)\vert&
\lesssim \min\{1, (R/\langle t\rangle)^3 \} N_{1},\\
\vert\phi^0_{r}(x,t)\vert+\vert E^{0}_{r}(x,t)\vert&
\lesssim \min \{1,r^3\}N_{1}
\end{split}
\end{equation}
and
\begin{equation}\label{ScaleEstDiffPhi2}
\begin{split}
\mathfrak{1}_{\{r\le |a|/10\}}\left[\langle a\rangle\vert\phi^0_{r}(a,t)\vert+\langle a\rangle^{2}\vert E^{0}_{r}(a,t)\vert\right]&
\lesssim r^3N_{1}.
\end{split}
\end{equation}
\end{lemma}
\begin{proof}[Proof of Lemma \ref{crude}]
We {treat} only the bounds on $E$, the bounds on $\phi$ are similar (easier).
We start with the bounds \eqref{ScaleEstDiffPhi} and \eqref{ScaleEstDiffPhi2}. 
Using \eqref{QuantEGC} and \eqref{QuantEG0}, we observe that, on the support of $\varphi(R^{-1}|x-y|)$,
\begin{equation}\label{DecayGDiffE}
\begin{split}
\vert G(x,y)\vert+R\vert \partial_xG(x,y)\vert+R^2\vert \partial_{x^j}\partial_{x^k}G(x,y)\vert&\lesssim R^{-1}, \quad
\vert G(x,y)\vert+R\vert \partial_xG(x,y)\vert  \lesssim R^{-2}|y|.
\end{split}
\end{equation}
For $p=0,1,2$, {by \eqref{DecayGDiffE}}, a crude estimate gives the simple bounds
\begin{equation}\label{L2}
\begin{split}
\vert E_{R}(x,t)\vert+\vert E^0_{r}(x,t)\vert &\lesssim \Vert \gamma \Vert_{L^2_{z,w}}^2 \lesssim N_{1},\\
\vert E_{R}(x,t)\vert &\lesssim R^{-1} [ \Vert \langle z\rangle^{\frac12} \gamma \Vert_{L^2_{z,w}}^2 +  \langle t \rangle \Vert \langle w\rangle^{\frac12} \gamma \Vert_{L^2_{z,w}}^2] \lesssim R^{-1} \langle t \rangle N_{1},\\
\mathfrak{1}_{\{R\le \vert x\vert/10\}}\langle x\rangle^p\vert E_{R}(x,t)\vert&\lesssim  \Vert \langle z\rangle^{\frac{p}{2}}\gamma \Vert_{L^2_{z,w}}^2+\langle t\rangle^{p}\Vert \langle w\rangle^{\frac{p}{2}}\gamma \Vert_{L^2_{z,w}}^2 \lesssim \langle t\rangle^{p}  N_{1},\\
\mathfrak{1}_{\{r\le\vert a\vert/10\}}\langle a\rangle^p \vert E^0_{r}(a,t)\vert& \lesssim \Vert \langle w\rangle^{\frac{p}{2}}\gamma \Vert_{L^2_{z,w}}^2 \lesssim N_{1}.
\end{split}
\end{equation}
Independently, another crude estimate gives the simple bounds
\begin{equation}\label{Linfty}
\begin{split}
\vert E_{R}(x,t)\vert &\lesssim \min\{R^{3}\Vert \langle w\rangle^2\gamma\Vert_{L^\infty_{z,w}}^2, (R/t)^3\Vert \langle z\rangle^2\gamma\Vert_{L^\infty_{z,w}}^2\} \lesssim (R/\langle t\rangle)^3 N_{1} ,\\
\mathfrak{1}_{\{R\le \vert x\vert/10\}}\langle x\rangle^p \vert E_{R}(x,t)\vert &\lesssim  (R/\langle t\rangle)^3[\Vert \langle z\rangle^{2+\frac{p}{2}}\langle w\rangle^2\gamma \Vert_{L^\infty_{z,w}}^2+\langle t\rangle^p\Vert \langle z\rangle^2\langle w\rangle^{2+\frac{p}{2}}\gamma \Vert_{L^\infty_{z,w}}^2] 
\\
&\lesssim (R/\langle t\rangle)^3\langle t\rangle^p N_{1},\\
\vert E^{0}_{r}(x,t)\vert&\lesssim r^3\Vert \langle z\rangle^2\gamma\Vert_{L^\infty_{z,w}}^2 \lesssim r^{3} N_{1},\\
\mathfrak{1}_{\{r\le \vert a\vert/10\}}\langle a\rangle^p \vert E_{r}^{0}(a,t)\vert &\lesssim r^3\Vert \langle z\rangle^2 \langle w\rangle^{\frac{p}{2}} \gamma\Vert_{L^\infty_{z,w}}^2 \lesssim r^{3} N_{1}.
\end{split}
\end{equation}

{
For the reader's convenience, we demonstrate several estimates of $E_{R}$ in  \eqref{L2} and \eqref{Linfty}.
Using \eqref{DecayGDiffE}, we can obtain the first line of \eqref{L2} as 
\begin{align*}
		|E_R(x,t)|\lesssim R^2\iint   R^{-2} \gamma^{2}(z,w,t) {dz dw} = \Vert \gamma \Vert_{L^2_{z,w}}^2  \leq \Vert \langle z\rangle^2\langle w\rangle^2\gamma\Vert_{L^\infty_{z,w}}^2 \iint \frac{dzdw}{\langle z\rangle^4 \langle w\rangle^4 } \lesssim N_{1}.
\end{align*}
We note that $|x| \leq 4R + |z-tw|$ holds on the support of $\varphi(R^{-1}|x-(z-tw)|)$, and this with $R\le \vert x\vert/10$ yields $|x| \lesssim |z-tw|$. Then the third line of \eqref{L2} can be shown as
\begin{align*}
\mathfrak{1}_{\{R\le \vert x\vert/10\}}\langle x\rangle^p\vert E_{R}(x,t)\vert
&\lesssim R^2\iint   R^{-2} [\langle z\rangle^{p}+\langle t\rangle^{p} \langle w\rangle^{p}] \gamma^{2}(z,w,t) {dz dw}  \lesssim \langle t\rangle^{p}  N_{1}.
\end{align*}
We now turn to \eqref{Linfty}.
From \eqref{DecayGDiffE}, we can derive the first line of \eqref{Linfty} as follows:
\begin{align*}
		&|E_R(x,t)|\lesssim \Vert \langle w\rangle^2\gamma\Vert_{L^\infty_{z,w}}^2 R^2\int_{w\in \R^3}\int_{|z-(x+tw)|\leq 10R }  R^{-2} \frac{dz dw}{\langle w\rangle^{4}}\lesssim R^3 N_{1}, \\ 
		&|E_R(x,t)|\lesssim  \Vert \langle z\rangle^2\gamma\Vert_{L^\infty_{z,w}}^2 \int_{z\in \R^3}\int_{|w-(z-x)/t|\leq 10R/t } \frac{dzdw}{\langle z\rangle^4 }\lesssim (R/t)^3 N_{1}. 
\end{align*}
The second line can be shown as
\begin{align*}
		&\mathfrak{1}_{\{R\le \vert x\vert/10\}}\langle x\rangle^p |E_R(x,t)|
		\\
		&\lesssim  \left[\Vert \langle z\rangle^{\frac{p}{2}} \langle w\rangle^2\gamma\Vert_{L^\infty_{z,w}}^2 +\langle t\rangle^{p} \Vert \langle w\rangle^{2+\frac{p}{2}}\gamma\Vert_{L^\infty_{z,w}}^{2} \right] R^2\int_{w\in \R^3}\int_{|z-(x+tw)|\leq 10R }  R^{-2} \frac{dz dw}{\langle w\rangle^{4}}\lesssim R^3 N_{1}, \\ 
		&\mathfrak{1}_{\{R\le \vert x\vert/10\}}\langle x\rangle^p |E_R(x,t)|
		\\
		&\lesssim \left[\Vert \langle z\rangle^{2+\frac{p}{2}} \gamma\Vert_{L^\infty_{z,w}}^2 +\langle t\rangle^{p} \Vert \langle z\rangle^{2} \langle w\rangle^{\frac{p}{2}}\gamma\Vert_{L^\infty_{z,w}}^{2} \right] \int_{z\in \R^3}\int_{|w-(z-x)/t|\leq 10R/t } \frac{dzdw}{\langle z\rangle^4 }\lesssim (R/t)^3 N_{1}. 
\end{align*}
The others follow similarly.
}
%
From \eqref{L2} and \eqref{Linfty}, we deduce the estimate on $E$ in \eqref{ScaleEstDiffPhi} and \eqref{ScaleEstDiffPhi2}.
\end{proof}

\subsubsection{{Leading-order asymptotics}}

{
To show that $E^{0}$ gives the leading-order asymptotics of $E$, 
we prove that the difference $E - E^{0}$ decays faster than $t^{-2}$, 
which is the optimal decay of $E$ stated in \eqref{Boot2CClD}.
}

\begin{proposition}\label{ControlDiffEProp}
Under the same assumption as in Proposition \ref{ProprBoot1}, there holds  that for $x\in D$ and $t \geq 1$,
\begin{equation}\label{ControlDiff}
\begin{split}
\sqrt{1+t^2+\vert x\vert^2}\vert {\phi}(x,t)-\phi^0(\pi_\infty(x/t),t)\vert&\lesssim \langle t\rangle^{-\frac{1}{20}}N_1,\\
\left[1+t^2+\vert x\vert^2\right]\vert{E}(x,t)-E^0(\pi_\infty(x/t),t)\vert&\lesssim\langle t\rangle^{-\frac{1}{20}}N_1.
\end{split}
\end{equation}
\end{proposition}

\begin{proof}[Proof of Proposition \ref{ControlDiffEProp}]
We prove \eqref{ControlDiff} only for $E$. We first show that the electric field is well approximated by its value in $\overline{D_\infty}$:
\begin{equation}\label{VanishAwayFromDInfty}
\begin{split}
\sup_{x \in D \backslash  \overline{D_\infty}}\left[1+t^2+\vert x\vert^2\right]\vert E(x,t)-E(\pi_\infty(x),t)\vert&\lesssim \langle t\rangle^{-\frac1{10}}N_1.
\end{split}
\end{equation}
We see from \eqref{DecayGDiffE} that
\begin{equation*}
\begin{split}
[1+\vert x \vert^{2} \mathfrak{1}_{\{R\le\vert x\vert/10\}}]\cdot\vert E_R(x,t)-E_R(\pi_\infty(x),t)\vert&\lesssim (\vert x-\pi_\infty(x)\vert/R)\cdot\min\{1,(R/\langle t\rangle)^3\}N_1.
\end{split}
\end{equation*}
{It is also seen from \eqref{UCAssumption} that $\vert x-\pi_\infty(x) \vert \lesssim 1$.}
In the case $\vert x\vert\le\langle t\rangle^{\frac{6}{5}}$, we use the simple bounds \eqref{L2}--\eqref{Linfty} for the small and large scales to obtain
\begin{equation*}
\begin{split}
\vert E(x,t)-E(\pi_\infty(x),t)\vert&\lesssim N_1\left[\int_{R=0}^{\langle t\rangle^{\frac12}}\langle t\rangle^{-3}R^3\cdot \frac{dR}{R^3}+\int_{R=\langle t\rangle^2}^\infty\frac{dR}{R^3}+\int_{R=\langle t\rangle^{\frac12}}^{\langle t\rangle^2}\langle t\rangle^{-3}R^2\vert x-\pi_\infty(x)\vert\frac{dR}{R^3}\right]\\
&\lesssim N_1\langle t\rangle^{-2}\left[\langle t\rangle^{-\frac12}+\langle t\rangle^{-1}\ln \langle t\rangle\cdot\vert x-\pi_\infty(x)\vert\right].
\end{split}
\end{equation*}
In another case $\vert x\vert\ge\langle t\rangle^{\frac65}$, we compute 
\begin{equation*}
\begin{split}
&\vert x\vert^2\vert E(x,t)-E(\pi_\infty(x),t)\vert 
\\
&\lesssim N_1\left[\int_{R=0}^{\langle t\rangle^{\frac12}}{\langle t\rangle^{-3}}R^3\cdot \frac{dR}{R^3}+\int_{R=\langle t\rangle^{\frac65}/10}^{\vert x\vert/10}\frac{dR}{R^3}+\int_{R=\langle t\rangle^{\frac12}}^{\langle t\rangle^{\frac65}/10}{\langle t\rangle^{-3}}R^2\vert x-\pi_\infty(x)\vert\frac{dR}{R^3}\right]\\
& \qquad +\vert x\vert^2 \langle t \rangle N_1\int_{R=\vert x\vert/10}^\infty \frac{dR}{R^{4}}\\
&\lesssim N_1{\langle t\rangle^{-\frac{1}{5}}}.
\end{split}
\end{equation*}
These give the bound \eqref{VanishAwayFromDInfty}.

We set $D_{\infty,t}:=\{{\frac{x}{t}} \in D_{\infty} \, : \, d(\frac{x}{t},C_{\infty}) \geq \langle t\rangle^{-\frac14} \}$
and also denote by $\pi_{\infty,t}$ the projection onto $D_{\infty,t}$. 
For $\frac{x}{t} \in \overline{D_{\infty}} \backslash D_{\infty,t}$,
there exists $\frac{y}{t} \in \{\frac{x}{t} \in D_{\infty}  \, : \, d(\frac{x}{t},C_{\infty}) = \langle t\rangle^{-\frac14} \}$ such that $\pi_{\infty,t}(\frac{x}{t})=\frac{y}{t}$ and $|\frac{x}{t}-\pi_{\infty,t}(\frac{x}{t})| \leq \langle t\rangle^{-\frac 14}$.
We claim that
\begin{equation}\label{VanishAwayFromDInfty1}
\begin{split}
\sup_{\frac{x}{t} \in \overline{D_{\infty}} \setminus D_{\infty,t}}\left[1+t^2+\vert x\vert^2\right]\vert E(x,t)-E(t\pi_{\infty,t}({x}/{t}),t)\vert&\lesssim \langle t\rangle^{-\frac1{20}}N_1.
\end{split}
\end{equation}
We can show this claim as follows. It is seen from \eqref{DecayGDiffE} and $y=t\pi_{\infty,t}(\frac{x}{t})$ that
\begin{equation*}
\begin{split}
[1+\vert x \vert^2\mathfrak{1}_{\{R\le\vert x\vert/10\}}]\cdot\vert E_R(x,t)-E_R(y,t)\vert
&\lesssim (t\vert x-y\vert/t R)\cdot\min\{1,(R/\langle t\rangle)^3\}N_1 
\\
&\lesssim (\langle t\rangle^{\frac34} /R)\cdot\min\{1,(R/\langle t\rangle)^3\}N_1.
\end{split}
\end{equation*}
In the case $\vert x\vert\le\langle t\rangle^{\frac{17}{16}}$, we use the simple bounds \eqref{L2}--\eqref{Linfty} for the small and large scales to obtain
\begin{equation*}
\begin{split}
\vert E(x,t)-E(y,t)\vert&\lesssim N_1\left[\int_{R=0}^{\langle t\rangle^{-1}}\langle t\rangle^{-3}R^3\cdot \frac{dR}{R^3}+\int_{R=\langle t\rangle^2}^\infty\frac{dR}{R^3}+\int_{R=\langle t\rangle^{-1}}^{\langle t\rangle^2}\langle t\rangle^{-3}R^2 \langle t\rangle^{\frac34} \frac{dR}{R^3}\right]\\
&\lesssim N_1\langle t\rangle^{-2}\left[\langle t\rangle^{-2}+\langle t\rangle^{-\frac14}\ln \langle t\rangle \right].
\end{split}
\end{equation*}
In another case $\vert x\vert\ge\langle t\rangle^{\frac{17}{16}}$, we compute 
\begin{equation*}
\begin{split}
&\vert x\vert^2\vert E(x,t)-E(y,t)\vert 
\\
&\lesssim N_1\left[\int_{R=0}^{\langle t\rangle^{-1}}{\langle t\rangle^{-3}}R^3\cdot \frac{dR}{R^3}+\int_{R=\langle t\rangle^{\frac{17}{16}}/10}^{\vert x\vert/10}\frac{dR}{R^3}+\int_{R=\langle t\rangle^{-1}}^{\langle t\rangle^{\frac{17}{16}}/10}{\langle t\rangle^{-3}}R^2\langle t\rangle^{\frac{3}{4}} \frac{dR}{R^3}\right]\\
& \qquad +\vert x\vert^2 \langle t \rangle N_1\int_{R=\vert x\vert/10}^\infty \frac{dR}{R^{4}}\\
&\lesssim N_1{\langle t\rangle^{-\frac{1}{16}}}.
\end{split}
\end{equation*}
These give \eqref{VanishAwayFromDInfty1}.

We can now assume that $x \in \overline{D_{\infty}}$ and $\frac{x}{t} \in D_{\infty,t}$. Let us replace $E^{0,j}$ by
\begin{equation*}
\begin{split}
E^{1,j}(a,t)&:= \frac{1}{t^{2}} \int_{r=0}^{\infty} E^{1,j}_{r}(a,t) \frac{dr}{r^{3}}, \\
E^{1,j}_{r}(a,t)&:= r^2\iint \mathfrak{1}_{\{w\in D_\infty\}}  \partial_{a^j}\left[\varphi(r^{-1}|a-w|)G_t(a,w)\right]\gamma^2(z,w,t)dzdw.
\end{split}
\end{equation*}
More precisely, we show that
\begin{gather}\label{E0E1}
\left[1+t^2+\vert x\vert^2\right]  \left|E^{0,j}(\frac{x}{t},t)-E^{1,j}(\frac{x}{t},t)\right| \lesssim \langle t\rangle^{-\frac{1}{20}}N_1  \quad \hbox{for } \frac{x}{t} \in D_{\infty,t}, \ t\geq 1.
\end{gather}
We set $H(a,w):=|\partial_{a_{j}} [\varphi(r^{-1}|a-w|)(G_\infty(a,w)- G_t(a,w))]|$ and observe that
\begin{equation}
\begin{split}
\left[1+\langle a \rangle^2\mathfrak{1}_{\{r\le \langle a \rangle/10\}}\right] |E^{0,j}_r(a,t)-E^{1,j}_r(a,t)|
 \leq r^{2}\iint  \mathfrak{1}_{\{w\in D_\infty\}} H(a,w)\langle w\rangle^2\gamma^2(z,w,t)dzdw.
\end{split}
\label{withlem2.3}
\end{equation}

Let us first estimate $H(a,w)$.
For $a\in \overline{D_\infty}$ and $y\in C_{\infty}$, using the mean value theorem and the {{third}} line in \eqref{QuantEGC}, we have
{{\begin{equation*}
\begin{split}
|G_\infty(a,y)-G_t(a,y)|
=|tG(ta,ty)|=\left| t[G(ta,{\bf w}_{C_{\infty}}(ty))-G(ta,ty)] \right|&\lesssim  \sup_{z\in[ty,{\bf w}_{C_\infty} (ty)]} t \frac{|ty-{\bf{w}}_{C^\infty}(ty)|}{|ta-z|^2} \\ &\leq t \frac{\vert ty-{\bf w}_{C_{\infty}}(ty)\vert}{\vert t a-t y\vert^2},
\end{split}
\end{equation*}
where the last inequality follows from $\pi_\infty(z)=ty$ and 
$\langle z-ty, ta-ty \rangle\leq 0. $ }}
If $a \in D_{\infty,t}$, using \eqref{UCAssumption} and the maximum principle,  we deduce that
\begin{equation*}
\begin{split}
|G_\infty(a,w)-G_t(a,w)| \lesssim \frac{1}{td(a,C_{\infty})^2} \lesssim \langle t\rangle^{-\frac{1}{2}}.
\end{split}
\end{equation*}
Similarly, it follows from the second line in \eqref{QuantEGC} that
\begin{equation*}
\begin{split}
|\nabla_{a}[ G_\infty(a,w)-G_t(a,w)]| \lesssim \frac{1}{td(a,C_{\infty})^3} \lesssim \langle t\rangle^{-\frac{1}{4}}.
\end{split}
\end{equation*}
These lead to
\begin{gather*}
H(a,w) \lesssim r^{-1} \langle t\rangle^{-\frac{1}{2}} + \langle t\rangle^{-\frac{1}{4}}.
\end{gather*}
On the other hand, we observe on the support of $H$ that
\begin{gather*}
H(a,w) \lesssim  r^{-3} |w|.
\end{gather*}
From \eqref{withlem2.3} and the estimate of $H(a,w)$, it holds that
\begin{gather*}
\left[1+\langle a \rangle^2\mathfrak{1}_{\{r\le \langle a \rangle/10\}}\right] |E^{0,j}_r(a,t)-E^{1,j}_r(a,t)| \lesssim \min\{r \langle t\rangle^{-\frac{1}{2}} + r^{2} \langle t\rangle^{-\frac{1}{4}}, r^{-1}\} N_{1}.
\end{gather*}
Integrating this in $r$, we deduce that
\begin{align*}
& t^{2}\langle a \rangle^2   |E^{0,j}(a,t)-E^{1,j}(a,t)| 
\\
&\lesssim N_1 \left[1+t^{\frac16} \mathfrak{1}_{\{\langle a \rangle < t^{\frac1{12}} \}} \right] \cdot \left[\int_{r=0}^{\langle t\rangle^{-\frac14}} r^{3} \frac{dr}{r^{3}} +\int_{r=\langle t\rangle^{-\frac14}}^{\langle t\rangle^{\frac14}} \left[ r \langle t\rangle^{-\frac{1}{2}} + r^{2} \langle t\rangle^{-\frac{1}{4}} \right] \frac{dr}{r^3}+\int_{r=\langle t\rangle^{\frac14}}^{\infty} \frac{dr}{r^3}\right]
\\
& \quad + N_{1} \langle a \rangle^2 \mathfrak{1}_{\{\langle a \rangle \geq t^{\frac1{12}}\}}  \int_{r= \langle a \rangle/10}^{\infty}\frac{dr}{r^4}
\\
&\lesssim N_1 \langle t\rangle^{-\frac{1}{20}}.
\end{align*}
Using this and letting $x=ta$,  we arrive at the bound \eqref{E0E1}.

We complete the estimate on $E$ in \eqref{ControlDiff} by estimating $E-E^1$ when $x\in \overline{D_\infty}$.  
Using Lemma \ref{EscapeLemma}, we see that 
\begin{equation}\label{DiffSet}
\begin{split}
\vert\mathfrak{1}_{\{b(z+tw)>0\}}-\mathfrak{1}_{\{w\in D_\infty\}}\vert&=\mathfrak{1}_{\{b(z+tw)>0\}}\mathfrak{1}_{\overline{\{j(w)>0\}}}+\mathfrak{1}_{\{b(z+tw)\le 0\}}\mathfrak{1}_{\{w\in D_\infty\}}
\lesssim \mathfrak{1}_{J_{\alpha}}+\left(\frac{\langle z\rangle}{\alpha\vert w\vert t}\right)^2,
\end{split}
\end{equation}
where $\alpha \in (0,\frac{1}{10})$ are constants to be determined later.
We observe from $G_t(x,y)=tG(tx,ty)$ that $\partial_{x^j} G_t(x,y)=t^2 (\partial_{x^j}G)(tx,ty)$. Using \eqref{DiffSet} and letting $\widetilde{H}(a,b)=\partial_{a^j}[\varphi(r^{-1}|a-b|)G_t(a,b)]$, we have
\begin{equation*}
\begin{split}
\mathfrak{d}_r^{j}&:=[1+\langle a\rangle^2\mathfrak{1}_{\{r\le\langle a\rangle/10\}}]\vert E^{j}_{rt}(ta,t)-E^{1,j}_{r}(a,t)\vert\\
&\lesssim r^2\iint \left[\langle w\rangle+\frac{\langle z\rangle}{t}\right]^2\gamma^2(z,w,t)\left[\frac{\langle z\rangle^2}{\alpha^2 t^2\vert w\vert^2}+\frac{\langle z\rangle^2}{r^2t^2}\right]\cdot\left[ \vert \widetilde{H} (a,w)\vert+\vert \widetilde{H} (a,w+\frac{z}{t})\vert\right]dwdz\\
&\qquad +r^2\iint_{\{\langle z\rangle\le rt/10\}} \langle w\rangle^2\gamma^2(z,w,t)\cdot\mathfrak{1}_{J_\alpha}\cdot\vert \widetilde{H} (a,w)\vert dwdz\\
&\qquad+r^2\iint_{\{\langle z\rangle\le rt/10\}} \langle w\rangle^2\gamma^2(z,w,t)\left\vert \widetilde{H} (a,w)-\widetilde{H} (a,w+\frac{z}{t})\right\vert dwdz.\\
\end{split}
\end{equation*}
Let us estimate the terms on the right hand side one by one.
The first term is estimated by using \eqref{QuantEGC}  as
\begin{equation*}
\begin{split}
\hbox{(1st term)} &\lesssim (\alpha t)^{-2}\min\{1,r\}\left[ \Vert \langle z\rangle^4\langle w\rangle^3\gamma\Vert_{L^\infty_{z,w}}^2 + \Vert \langle z\rangle^2\langle w\rangle\gamma\Vert_{L^2_{z,w}}^2 \right]
\lesssim (\alpha t)^{-2}\min\{1,r\} N_{1}.
\end{split}
\end{equation*}
Using \eqref{QuanfVanishingG}, we can estimate the second term as
\begin{equation*}
\begin{split}
\hbox{(2nd term)}&\leq r^2\iint_{J_\alpha} \langle w\rangle^2\gamma^2(z,w,t)\cdot\vert \widetilde{H} (a,w)\vert dwdz
\\
&\le\Vert \langle z\rangle^2\langle w\rangle^4\gamma\Vert_{L^\infty_{z,w}}^2\iint_{{J_\alpha}}\frac{dw}{\langle w\rangle^4} \frac{dz}{\langle z\rangle^4}
\lesssim  {\alpha} N_{1}.
\end{split}
\end{equation*}
The third term can be handled by the mean value theorem and \eqref{QuantEGC} as
\begin{equation*}
\begin{split}
\hbox{(3rd term)} &\lesssim 
\frac{1}{rt}\iint_{\{\vert a-w\vert\le 10r\}} \langle z\rangle\langle w\rangle^2\gamma^2(z,w,t)dwdz
\lesssim \min\{r^{-1},r^2\}t^{-1}N_{1}.
\end{split}
\end{equation*}
Choosing 
$\alpha=\frac{1}{10}\langle t\rangle ^{-\frac{2}{3}}$, we obtain 
\begin{equation*}
\begin{split}
\mathfrak{d}_r^{j}\lesssim {t^{-\frac{2}{3}}}N_1.
\end{split}
\end{equation*}
Combining with the simple bounds $\mathfrak{d}_r^{j}\lesssim \min\{1,r^{3}\}N_1$
and plugging in \eqref{DecEScale} with $x=ta$, we deduce that
\begin{align*}
&\left[1+t^2+\vert x\vert^2\right]\left\vert {E}(x,t)-E^{1}(\frac{x}{t},t)\right\vert
\\
&\lesssim N_{1}\left[\int_{r=0}^{{\langle t\rangle^{-\frac{2}{9}}}} r^{3}\frac{dr}{r^3}
+ \int_{r={\langle t\rangle^{-\frac{2}{9}}}}^{\langle y\rangle/10} {\langle t\rangle^{-\frac{2}{3}}} \frac{dr}{r^3} 
+\langle y\rangle^2\int_{r=\langle y\rangle/10}^{\infty}\langle t\rangle^{-\frac{2}{3}}\frac{dr}{r^3}\right]
\lesssim N_1 {\langle t\rangle^{-\frac{2}{9}}}.
\end{align*}
Thus we arrive at \eqref{ControlDiff}.
\end{proof}

\subsubsection{Refined estimate}

{ From Proposition \ref{ControlDiffEProp}, we see that in order to obtain the optimal bound for $E$ in \eqref{Boot2CClD}, 
it suffices to refine the estimate of $E^{0}$. 
We therefore proceed to estimate $E^{0}_{r}$ by appealing to the evolution equation of $\gamma$. }

\begin{proposition}\label{ControlE0Prop}
Under the same assumption as in Proposition \ref{ProprBoot1}, there holds that for $a \in D_{\infty}$ and $0\leq t \leq T$,
\begin{equation}\label{ScaleEstDtPhi}
\begin{split}
\langle a\rangle \vert \partial_{t} \phi_{r}^{0}(a,t) \vert
&\lesssim \min\{1,r^2\} \| F(t) \|_{L^{\infty}_{x}} N_{1}, \\ 
\langle a\rangle^{2}\vert \partial_{t} E_{r}^{0}(a,t) \vert &\lesssim r \min\{1,r\} \| F(t) \|_{L^{\infty}_{x}}N_{1}
\end{split}
\end{equation}
and in particular, we have the uniform in time bound:
\begin{equation}\label{UniformBoundsPhi0}
\begin{split}
\langle a\rangle\vert\phi^0_{r}(a,t)\vert&\lesssim\varepsilon_1^2\min \{1,r^{2}\},
\\
\langle a\rangle^2\vert E^0_{r}(a,t)\vert&\lesssim\varepsilon_1^2 r\min \{1,r\}.
\end{split}
\end{equation}
In addition, if the same assumption in Proposition \ref{ProprBoot2} holds, then for $a \in D_{\infty}$, $0\leq t \leq T$, and $r\leq 1$,
\begin{equation}\label{ScaleEstDtE01}
\begin{split}
\langle a\rangle^{2}\vert \partial_{t} E_{r}^{0,j}(a,t) \vert &\lesssim r^{3} \left[ \| E(t) \|_{L^{\infty}_{x}(D)}  + t \| \nabla E(t) \|_{L^{\infty}_{x}(D)} \right] \left[N_{1}^{1/2}\Vert \mathfrak{1}_{\{b(z+tw)>0\}} \nabla_{w} \gamma\Vert_{L^\infty_{z,w}}+N_{1}\right],
\end{split}
\end{equation}
and in particular, we have the uniform in time bound:
\begin{equation}\label{UniformBoundsE01}
\begin{split}
\langle a\rangle^2\vert E^0_{r}(a,t)\vert& \lesssim\varepsilon_1^2 r^{3}.
\end{split}
\end{equation}

\end{proposition}

\begin{proof}[Proof of Proposition \ref{ControlE0Prop}]
We only detail the bounds on $E$, since the bounds on $\phi$ follow similarly. We start with the bound \eqref{ScaleEstDtPhi}.
From \eqref{QuantEGC}, we observe on the support of $\varphi(R^{-1}|a-w|)$ that
\begin{equation*}
\begin{split}
\vert G_{\infty}(a,w)\vert+R\vert \partial_aG_{\infty}(a,w)\vert+R^2\vert \partial_{a^j}\partial_{w^k}G_{\infty}(a,w)\vert&\lesssim R^{-1}.
\end{split}
\end{equation*}
Since $G_\infty(a,w)$ and $\partial_{a^j} G_\infty(a,w)$ vanish when $w\in C_{\infty}$, we see that
\begin{equation}\label{DEDT}
\begin{split}
\partial_t E^{0,j}_{r}(a,t)
&=r^2 \iint  \mathfrak{1}_{\{w\in D_\infty\}} \partial_{a^{j}}[\varphi(r^{-1}|a-w|)G_\infty(a,w)] \partial_t(\gamma^2)dzdw\\
&=r^2 \iint  \mathfrak{1}_{\{w\in D_\infty\}} \partial_{a^{j}}[\varphi(r^{-1}|a-w|)G_\infty(a,w)] \{\gamma^2,\Phi\} dzdw\\
&=r^2 \iint \mathfrak{1}_{\{w\in D_\infty\}}  \gamma^{2} \{\Phi,\partial_{a^{j}}[\varphi(r^{-1}|a-w|)G_\infty(a,w)]\} dzdw\\
&=r^2 \iint \mathfrak{1}_{\{w\in D_\infty\}}  \gamma^{2} \lambda F(z+tw,t)\cdot \nabla_w\partial_{a^{j}}[\varphi(r^{-1}|a-w|)G_\infty(a,w)] dzdw.
\end{split}
\end{equation}
From this, we arrive at \eqref{ScaleEstDtPhi}. We now deduce \eqref{UniformBoundsPhi0} by integrating \eqref{ScaleEstDtPhi} and the bootstrap assumptions \eqref{Boot1} to bound $\Vert F(t)\Vert_{L^\infty_x}$.

In order to prove \eqref{ScaleEstDtE01}, we need to use the derivative control on $\gamma$, but it is used only in the domain $D=\{b>0\}$.
Proceeding as in \eqref{DEDT}, we get
\begin{align*}
\partial_t E^{0,j}_{r}(a,t)&=r^2 \iint  \mathfrak{1}_{\{w\in D_\infty\}}\partial_{a^{j}}[\varphi(r^{-1}|a-w|)G_\infty(a,w)] \{\gamma^2,\Phi\}dzdw\\
&=-2\lambda r^2  \iint  \mathfrak{1}_{\{w\in D_\infty\}} \partial_{a^{j}}[\varphi(r^{-1}|a-w|)G_\infty(a,w)] \gamma F\cdot \nabla_w\gamma dzdw\\
&\quad+\lambda tr^2 \iint  \mathfrak{1}_{\{w\in D_\infty\}} \partial_{a^{j}}[\varphi(r^{-1}|a-w|)G_\infty(a,w)] F\cdot \nabla_z\gamma^2 dzdw\\
&=-2\lambda r^2 \iint  \mathfrak{1}_{\{w\in D_\infty\}}\mathfrak{1}_{\{b(z+tw)>0\}} \partial_{a^{j}}[\varphi(r^{-1}|a-w|)G_\infty(a,w)] \gamma E \cdot \nabla_w\gamma dzdw\\
&\quad-\lambda r^2 \iint  \mathfrak{1}_{\{w\in D_\infty\}}\mathfrak{1}_{\{b(z+tw)>0\}}  \partial_{a^{j}}[\varphi(r^{-1}|a-w|)G_\infty(a,w)] \hbox{div}_z(tE)\cdot \gamma^2 dzdw,
\end{align*}
where we have used Lemma \ref{SuffCondBC} to insert the additional characteristic function\footnote{Since on the support of $\varphi$, when $x=z+tw\in D^c$, then $w\cdot\nabla_xb(x)\ge g(x,w)>0$.} and also replaced $F$ by $E$.
From this, we arrive at \eqref{ScaleEstDtE01}. We now deduce \eqref{UniformBoundsE01} by integrating \eqref{ScaleEstDtE01} and
using the bootstrap assumptions \eqref{Boot1} and \eqref{Boot2}.
\end{proof}

{Propositions \ref{ControlDiffEProp} and \ref{ControlE0Prop} immediately gives the optimal bound for $E$ as follows:}

\begin{corollary}\label{ControlPhiE}
Under the same assumption as in Proposition \ref{ProprBoot1}, there holds that for $x\in D$ and $0\le t\le T$,
\begin{equation}\label{DecayE0}
\begin{split}
\sqrt{1+t^2+\vert x\vert^2}\vert {\phi}(x,t)\vert&\lesssim\varepsilon_1^2,\\
\left[1+t^2+\vert x\vert^2\right]\vert E(x,t)\vert&\lesssim\varepsilon_1^2\langle \ln (2+t)\rangle,
\end{split}
\end{equation}
In addition, if the same assumption in Proposition \ref{ProprBoot2} holds, then the optimal bound holds:
\begin{equation}\label{DecayE1}
\begin{split}
\left[1+t^2+\vert x\vert^2\right]\vert E(x,t)\vert&\lesssim\varepsilon_1^2.
\end{split}
\end{equation}
\end{corollary}

\begin{proof}[Proof of Proposition \ref{ControlPhiE}]
We {treat} only the bounds on $E$, since it is easier to obtain  the bounds on $\phi$.
This is clear if $\vert x\vert+t\le 1$. We suppose that $\vert x\vert+t\ge 1$.
Owing to \eqref{ControlDiff}, it suffices to get control on $E^0$ in the case $x\in D_\infty$. Using \eqref{UniformBoundsPhi0} for small scales and \eqref{ScaleEstDiffPhi2} for the other scale, we have
\begin{equation*}
\begin{split}
[t^2+\vert x\vert^2]\vert E^0(x/t,t)\vert &=\int_{r=0}^\infty [1+\vert x/t\vert^2]\vert E^0_{r}((x/t),t)\vert\frac{dr}{r^3}\\
&\lesssim\varepsilon_1^2\int_{r=\rho}^\infty \min\{1,r\}\frac{dr}{r^2}+\varepsilon_1^2\langle\ln(2+t)\rangle^{10}\int_{r=0}^\rho dr.
 \end{split}
\end{equation*}
Choosing $\rho=\langle\ln(2+t)\rangle^{-100}$ gives \eqref{DecayE0}.
Proceeding similarly with \eqref{UniformBoundsE01}, we deduce the optimal bound \eqref{DecayE1}.
\end{proof}

\subsection{Derivative of the Electric field}
In the case that $D\subset\mathbb{R}^3$ is an acceptable convex domain,  we can control the derivative of the electrostatic potential. 
As in the previous sections, we set
\begin{equation}\label{DefMR}
\begin{split}
&\partial_{x_{k}}E^{j}(x,t)=:{M}^{jk}(x,t)=\int_{R=0}^\infty {M}^{jk}_{R}(x,t)\frac{dR}{R^4}, 
\\
&{M}^{jk}_{R}(x,t):= R^3\iint \mathfrak{1}_{\{b(y)>0\}} \partial_{x^j}\partial_{x^k}\left[\varphi(R^{-1}|x-y|)G(x,y)\right]\mu^2(y,w,t)dydw.
\end{split}
\end{equation}
The thing to note here is that we use $\mu$ not $\gamma$, because it is easy to use the coordinate so that the boundary $\partial D$ is fixed in the proof of \eqref{ScaleEstDiffPhi3}, where the distance $b(x)$ from $\partial D$ appear.
\begin{proposition}\label{ControlDiffEProp*}
Under the same assumption as in Proposition \ref{ProprBoot2}, there holds that for $x\in D$ and $0\le t\le T$,
\begin{equation}\label{ScaleEstDiffPhi3}
\begin{split}
\vert {M}_{R}(x,t)\vert&\lesssim \min\{ 1,(R/\langle t\rangle)^{3}\} N_{1}.\\ 
\end{split}
\end{equation}
As a consequence, we have 
\begin{equation*}
\begin{split}
\langle t\rangle^{3} |\nabla E(x,t)| \lesssim  \varepsilon_{1}^{2} \langle\ln(2+t)\rangle^{11}.
\end{split}
\end{equation*}

\end{proposition}

\begin{proof}[Proof of Proposition \ref{ControlDiffEProp}]
It is straightforward  to get  \eqref{ScaleEstDiffPhi3} by using \eqref{QuantEGC} and \eqref{QuantEG0}.
We use the formula \eqref{DefMR} and the definition in \eqref{T1Hyp}, and then remove the constant to get
\begin{equation*}
\begin{split}
M_R^{jk}(x,t)&=R^3\iint\mathfrak{1}_{\{b>0\}}\partial_{x^j}\partial_{x^k}\left[\varphi(R^{-1}|x-y|)G(x,y)\right][\mu^2(y,v,t)-\mu^2(x,v,t)]dydv\\
&\quad+R^3n(x,t)\partial_{x^j}\partial_{x^k}u_R(x),
\end{split}
\end{equation*}
where 
\begin{gather*}
n (x,t) :=\int\mathfrak{1}_{\{b>0\}} \mu^2(x,v,t)dv. 
\end{gather*}
The first term can be estimated by using the mean-value theorem, \eqref{QuantEGC}, and \eqref{QuantEG0} as
\begin{equation*}
\begin{split}
\hbox{(1st term)}
\lesssim &R^4N_{1}^{1/2}\Vert\mathfrak{1}_{\{b>0\}}\nabla_x\mu(t)\Vert_{L^\infty_{x,v}}.
\end{split}
\end{equation*}
For the bound for $\nabla E(x,t)$, using \eqref{T1Hyp} we compute 
\begin{align*}
\langle t\rangle^{3} |\nabla E(x,t)| &\lesssim  
N_{1}^{1/2}\Vert\mathfrak{1}_{\{b>0\}}\nabla_x\mu(t)\Vert_{L^\infty_{x,v}} \langle t\rangle^{3}  \int_{R=0}^{\langle t\rangle^{-3}} R^{4}\frac{dR}{R^4}
\\
&\quad +N_{1}\left[\left|\int_{R=0}^{\langle t\rangle^{-3}}  \nabla^2_x u_{R} \frac{dR}{R} \right|
+ \int_{R=\langle t\rangle^{-3}}^{\langle t\rangle^{3}} R^{3}  \frac{dR}{R^4} 
+\langle t\rangle^{3} \int_{R=\langle t\rangle^{3}}^{\infty}\frac{dR}{R^{4}}\right]
\\
&\lesssim N_{1}^{1/2}\Vert\mathfrak{1}_{\{b>0\}}\nabla_x\mu(t)\Vert_{L^\infty_{x,v}} + N_{1} \langle\ln(2+t)\rangle.
\end{align*}
This with the bootstrap assumptions \eqref{Boot1} and \eqref{Boot2} gives the desired estimate.
\end{proof}

%
%
%
%
%
%

\section{Moments and derivative estimates, {and proofs of the key propositions}}\label{S5}

{In this section, we mainly study the bounds of $\mu$.
Subsections \ref{S5.1} and \ref{S5.2} deal with the control of the moment and derivatives, respectively.
Then subsection \ref{S5.3} completes the proofs of Propositions \ref{ProprBoot1} and \ref{ProprBoot2}.
}

\subsection{Control of the moments}\label{S5.1}
Assuming good control on the electric field $E$, we propagate control of the moments $\langle x-tv\rangle^4 \langle v\rangle^4\mu$. 
Given a weight $ \langle x-tv\rangle^k\langle v\rangle^l$ and a solution $\mu(t)$, we define the corresponding moment:
\begin{equation*}
\begin{split}
S_{k,l}(t):=\Vert \langle x-tv\rangle^k\langle v\rangle^l\mu  (x,v,t)\Vert_{L^\infty_{x,v}} \quad \hbox{for } k, l \in \mathbb N \cup \{0\}. 
\end{split}
\end{equation*}

\begin{proposition}\label{ControlOfMoments}
Suppose that the same assumption as in Proposition \ref{ProprBoot1} holds, and the force field $F$ satisfy the decay condition
\begin{equation}\label{BootwE0}
\begin{split}
\left[1+t^2 \right]\vert F(x,t)\vert&\lesssim \varepsilon_1^2\langle \ln(2+t)\rangle.\\
\end{split}
\end{equation}
Then there holds that for $k,l=0,1,\ldots,4$ and $0 \leq t \leq T$,
\begin{equation}\label{generalweight}
	S_{k,l}(t) \leq \varepsilon_0+C\varepsilon_0 \varepsilon_1 \langle \ln(t+2)\rangle^{2k}.
\end{equation}
\end{proposition}

\begin{proof}[Proof of Proposition \ref{ControlOfMoments}]
It is seen from \eqref{RVPW} that
\begin{equation*}
\begin{split}
\partial_t(\omega\mu )+\{(\omega\mu ),\mathcal{H}\}&=\mu [\partial_t\omega+\{\omega,\mathcal{H}\}]=\mu [\partial_t\omega+v\cdot\nabla_x \omega+\lambda \nabla_{x} \psi \cdot\nabla_v\omega],
\end{split}
\end{equation*}
where $\omega= \langle x-tv\rangle^k\langle v\rangle^l$.
Integrating this along the flow, we have
\begin{gather}\label{esSkl}
S_{k,l}(t) \leq S_{k,l}(0) + \int_{s=0}^t\Vert F(s)\Vert_{L^\infty_{x}}\left[skS_{k-1,l}+lS_{k,l-1}\right]ds.
\end{gather}
We note that $S_{0,0}$ is conserved.
{
By induction with the aid of \eqref{BootwE0} and \eqref{esSkl}, 
we obtain \eqref{generalweight} with $S_{k,0}$ for $k=1,\ldots,4$ and \eqref{generalweight} with $S_{0,l}$ for $l=1,\ldots,4$.
From these, \eqref{BootwE0}, and \eqref{esSkl}, it is seen that \eqref{generalweight} with $S_{1,1}$ holds. 
Using induction again, we have \eqref{generalweight} with $S_{k,1}$ for $k=2,\ldots,4$ and \eqref{generalweight} with $S_{1,l}$ for $l=2,\ldots,4$.
Then \eqref{generalweight} with $S_{2,2}$ also follows.
Similarly, for $i=2,3,4$, we conclude \eqref{generalweight} with  $S_{k,i}$ for $k=i,\ldots,4$ and \eqref{generalweight} with  $S_{i,l}$ for $l=i,\ldots,4$. Thus \eqref{generalweight} holds.
}
\end{proof}

\subsection{Control of the derivatives}\label{S5.2}

We recall that 
\begin{align*}
(\nabla_{z} \gamma)(z,w,t)&= (\nabla_{x} \mu)(z+tw,w,t), \\
(\nabla_{w} \gamma)(z,w,t)&=  (t \nabla_{x} \mu)(z+tw,w,t) + (\nabla_{v} \mu)(z+tw,w,t). 
\end{align*}
We estimate the supremum of the following $Z$ and $W$: 
\begin{gather*}
Z(x,v,t):=(\nabla_{x} \mu)(x,v,t), \quad W(x,v,t):=(t\nabla_x\mu+\nabla_v\mu)(x,v,t).
\end{gather*}
It is straightforward to see from \eqref{RVPW} that 
\begin{equation}\label{EvolDer}
\begin{split}
&\partial_{t} Z +\{Z,\mathcal{H}\}+ \lambda \nabla_{x} F (W-tZ)=0, \\
&\partial_{t} W +\{W,\mathcal{H}\}+ \lambda t \nabla_{x} F (W-tZ)=0.
\end{split}
\end{equation}


\subsubsection{Control of derivatives for the related Cauchy problem}
{We estimate $Z$ and $W$} in a whole space. 
\begin{proposition}\label{DerBoot1}
Suppose that the same assumption as in Proposition \ref{ProprBoot2} holds, and the force field $F$ satisfy the decay condition
\begin{equation}\label{BootwE}
\begin{split}
\left[1+t^2 \right]\vert \nabla F(x,t)\vert&\lesssim \varepsilon_1^2 \quad \hbox{for } x\in\mathbb R^{3}, t>0.\\
\end{split}
\end{equation}
%
Then there holds that for $0 \leq t \leq T$,
\begin{equation}\label{BootDerExt12}
\begin{split}
\Vert Z (t)\Vert_{L^\infty_{x,v}}&\le\varepsilon_0+C\delta_{1}^{-1}\varepsilon_1^2 \langle t\rangle^{\delta_{1}},\qquad
\Vert W(t)\Vert_{L^\infty_{x,v}}\le\varepsilon_0+C\varepsilon_1^2\langle t\rangle^{1+\delta_{1}}.
\end{split}
\end{equation}
\end{proposition}

\begin{proof}[Proof of Proposition \ref{DerBoot1}]
Integrating \eqref{EvolDer} along the flow, we see that
\begin{equation*}
\begin{split}
\Vert Z(t)\Vert_{L^\infty_{x,v}}&\le\Vert Z(0) \Vert_{L^\infty_{x,v}}+\int_{s=0}^t\Vert \nabla_xF(s)\Vert_{L^\infty_{x}}\left[s\Vert Z(s)\Vert_{L^\infty_{x,v}}+\Vert W(s)\Vert_{L^\infty_{x,v}}\right]ds,\\
\Vert W(t)\Vert_{L^\infty_{x,v}}&\le\Vert W(0)\Vert_{L^\infty_{x,v}}+\int_{s=0}^t\Vert \nabla_xF(s)\Vert_{L^\infty_{x}}\left[s^2\Vert Z (s)\Vert_{L^\infty_{x,v}}+s\Vert W(s)\Vert_{L^\infty_{x,v}}\right]ds.
\end{split}
\end{equation*}
These with the assumptions \eqref{Boot2} and \eqref{BootwE} immediately give \eqref{BootDerExt12}.
\end{proof}

While these bounds are strong enough to ensure the existence of a global strong solution, we will show refined estimates in the domain $D$.

\subsubsection{Refined control inside the domain}
We refine the estimate of $Z$ and $W$ inside the domain $D$. To do so, we set
\begin{gather*}
S_Z(t):=\Vert \mathfrak{1}_{\{b(x)>0\}}  Z (x,v,t)\Vert_{L^\infty_{x,v}}, \quad S_W(t):=\Vert \mathfrak{1}_{\{b(x)>0\}}  W (x,v,t)\Vert_{L^\infty_{x,v}}.
\end{gather*}

\begin{proposition}\label{ControlOfDerivatives}
Suppose that the same assumption as in Proposition \ref{ProprBoot2} holds, and {the force field $E$} satisfy the decay condition
\begin{equation}\label{BootwE2}
\begin{split}
\left[1+t^2 \right]^{3/2}\vert \nabla E(x,t)\vert& \lesssim \varepsilon_1^2\langle \ln(2+t)\rangle^{p},
\end{split}
\end{equation}
where $p \geq 1$.
Assume as well that the derivatives are bounded initially:
\begin{equation*}
\begin{split}
S_Z(0)+S_W(0)\le\varepsilon_0.
\end{split}
\end{equation*}
Then there holds that for $0\le t\le T$,
\begin{equation}\label{BootBounds2}
\begin{split}
S_Z(t) &\le \varepsilon_0+C\varepsilon_0\varepsilon_1^{2},\\
S_W(t) &\le\varepsilon_0 +C\varepsilon_0 \varepsilon_1^{2} \langle \ln(2+t)\rangle^{p+1}.
\end{split}
\end{equation}
\end{proposition}

\begin{proof}[Proof of Proposition \ref{ControlOfDerivatives}]

The proof is a simple variation on Gronwall lemma adapted to use the fact that trajectory just outside the domain are only outgoing. We give the details for the convenience of the reader.

For given nonnegative Lipshitz functions $A=A(t)$ and $\varphi=\varphi(\xi):=(\xi)_{+}$ as well as a given constant $B$, we define
\begin{equation*}
\begin{split}
0\le \widetilde{Z}(x,v,t)&=\langle (x-tv)\rangle^{-4}\langle v\rangle^{-4} \sum_{j=1}^{3}\varphi(Z^{j}(x,v,t)-B-A(t))\in C^0_tL^2_{x,v},
\end{split}
\end{equation*}
which satisfies
\begin{gather*}
\partial_t\widetilde{Z}+\{\widetilde{Z},\mathcal{H}\}=R, \\
R:=\langle (x-tv)\rangle^{-4}\langle v\rangle^{-4}  \sum_{j=1}^{3} \varphi' \cdot \left[-(\lambda \nabla_{x}F(W-tZ))^{j}-A'(t)\right]
+\lambda\frac{F \cdot \nabla_{v}\{\langle (x-tv)\rangle^{-4}\langle v\rangle^{-4}\}}{\langle (x-tv)\rangle^{-4}\langle v\rangle^{-4}}\widetilde{Z}.
\end{gather*}
We can integrate this in the region
$D_a:=\{b(x)>a\}$ for $-1\le a< 0$ to get
\begin{equation*}
\begin{split}
\iint_{D_a} \widetilde{Z}(t) dxdv-\int_{s=0}^t\int_{\{b(x)=a\}}\widetilde{Z} v\cdot{\bf n}d\sigma dvds&=\iint_{D_a} \widetilde{Z}(0) dxdv+\int_{s=0}^t\iint_{D_a}Rdxdvds.
\end{split}
\end{equation*}
Now, in the boundary $\{{b(x)}=a\}$, we have that $v\cdot{\bf n}=v\cdot\nabla b$ is nonpositive on the support of $\mu$ by Lemma \ref{SuffCondBC} and we deduce that
\begin{equation*}
\begin{split}
\iint_{D_a} \widetilde{Z} (t) dxdv&\le \iint_{D_a} \widetilde{Z}(0) dxdv+\int_{s=0}^t\iint_{D_a}Rdxdvds.
\end{split}
\end{equation*}
Letting $a \to 0$, we have
\begin{gather*}
\iint_{D} \widetilde{Z} (t) dxdv \le \iint_{D} \widetilde{Z}(0) dxdv + \int_{s=0}^t\iint_{D}Rdxdvds.
\end{gather*}
Now choosing $A$ and $B$ such that
\begin{equation*}\label{DefAB}
\begin{split}
A(t)=\int_{s=0}^t\Vert \nabla_xF(s)\Vert_{L^\infty_x(D)}[ S_W(s)+sS_Z(s)]ds, \quad 
B=\varepsilon_0 \geq S_Z(0),
\end{split}
\end{equation*}
we get for any fixed $t$,
\begin{equation*}
\begin{split}
\iint_{D} \widetilde{Z} (t) dxdv& \lesssim \int_{s=0}^t\Vert \langle s \rangle F(s)\Vert_{L^\infty_x(D)}\iint_{D}\widetilde{Z}  dxdvds.
\end{split}
\end{equation*}
Applying Gronwall's inequality, we see that $\widetilde{Z}(x,v,t)=0$ and so $(Z^{j}-\varepsilon_0-A(t) )_{+} = 0$ holds if $x \in D$. Thus
\begin{equation*}
\begin{split}
\Vert Z\Vert_{L^\infty_{x,v}(\{b>0\})}\le \varepsilon_0 + \int_{s=0}^t\Vert \nabla_xF(s)\Vert_{L^\infty_x(D)}[ S_W(s)+sS_Z(s)]ds.
\end{split}
\end{equation*}
Proceeding similarly for $S_W$, we arrive at
\begin{equation*}
\begin{split}
S_Z(t)&\le \varepsilon_0 + \int_{s=0}^t\Vert \nabla_xF(s)\Vert_{L^\infty_x(D)}[ S_W(s)+sS_Z(s)]ds, \\
S_W(t)&\le \varepsilon_0 + \int_{s=0}^t s\Vert \nabla_xF(s)\Vert_{L^\infty_x(D)}[ S_W(s)+sS_Z(s)]ds.
\end{split}
\end{equation*}
From these and \eqref{BootwE2}, we conclude \eqref{BootBounds2}.
\end{proof}

\subsection{{Proofs of the key propositions}}\label{S5.3}

We are now in a position to complete the proofs of Propositions \ref{ProprBoot1} and \ref{ProprBoot2}, which are keys for our bootstrap argument.

\begin{proof}[Proof of Proposition \ref{ProprBoot1}]
The bounds \eqref{Boot1CCl2} follow from Corollary \ref{ControlPhiE}; The bound for $F$ in \eqref{Boot1CCl} follows from \eqref{ExtPotAss11} and \eqref{Boot1CCl2}; The bound for the moment in \eqref{Boot1CCl} follows from Proposition \ref{ControlOfMoments}.
\end{proof}

\begin{proof}[Proof of Proposition \ref{ProprBoot2}]
The optimal bound for the electric field $E$ in \eqref{Boot2CClD} follows from Corollary \ref{ControlPhiE};
The bound for the derivative of $E$ in \eqref{Boot2CClD} follow from Proposition \ref{ControlDiffEProp*};
The bound for $\nabla F$ in \eqref{Boot2CCl} follows from \eqref{ExtPotAss12} and the bounds for $E$ in \eqref{Boot2CClD};
The bounds for the derivatives of $\mu$ in \eqref{Boot2CCl} follow from Proposition \ref{DerBoot1} with the bounds for $\nabla F$ in \eqref{Boot2CCl};
The precise bounds for the derivatives of $\mu$ in \eqref{Boot2CClD} follow from Proposition \ref{ControlOfDerivatives}.
\end{proof}

\section{Asymptotic Behavior}\label{S6}

We are now ready to extract the asymptotic behavior of solutions. We start by obtaining a limit for the acceleration:

\begin{lemma}\label{Convergenceofelectircfield}
Under the same assumption as in Proposition \ref{ProprBoot2}, 
there exists an asymptotic profile $E_\infty\in L^2\cap L^\infty\cap C^0(D_\infty)$ such that for $x\in D$ and sufficiently large $t \geq 1$, 
\begin{equation*}
	\begin{split}
	 \|\langle a\rangle^2(t^2 E^{0}(a,t)-E_{\infty}(a))\|_{L^{\infty}_a(D_\infty)}&\lesssim \varepsilon_1^2 t^{-\frac34},
\\
		\left[1+t^2+\vert x\vert^2\right] \vert E(x,t)-t^{-2}E_{\infty}(\pi_\infty(x/t))\vert&\lesssim \varepsilon_1^2  t^{-\frac{1}{50}}.
	\end{split} 
\end{equation*}	
\end{lemma}
	
\begin{proof}
It follows from \eqref{Boot1}, \eqref{Boot2}, \eqref{ScaleEstDtPhi}, and \eqref{ScaleEstDtE01} that
\begin{equation*}
\begin{split}
\langle a\rangle^2\vert\partial_tE^0_r(a,t)\vert\lesssim t^{-\frac74}\varepsilon_1^3\min\{r^3,r\}.
\end{split}
\end{equation*}
From this, there exists $E_\infty \in L^{\infty}_{a}$ such that
\begin{gather*}
t^{2} E^{0}(a, t) = \int_{s=1}^{t} \int_{r=0}^{\infty} {r^{-3}} \partial_{t} E^{0}_{r} (a,s) drds +E^{0} (a,t=1) \to E_{\infty}(a) \hbox{ as }t\to \infty.
\end{gather*}
This gives the first desired inequality and  $E_\infty\in L^2\cap L^\infty\cap C^0(D_\infty)$.
The second desired inequality also follows from \eqref{ControlDiff}.
\end{proof}

By Corollary \ref{CorSupp}, it makes sense to consider the asymptotic behavior only in the set $\{x\in\mathbb{R}^3,\,v\in D_\infty\}$. For $w\in D_\infty$, we define
\begin{equation}\label{DefSigma}
	\begin{split}
		\sigma(z,w,t):=\gamma(z+\lambda(\ln t)E_{\infty}(w),w,t).
	\end{split}
\end{equation}
We have pointwise convergence for  $\sigma(z,w,t)$.

\begin{proposition}\label{Convergenceofgamma}
Suppose that the same assumption as in Proposition \ref{ProprBoot2} holds. 
Then there exists a limit $\sigma_{\infty} \in L^2_{z,w}\cap L^\infty_{z,w} \cap C_{z,w}^{0}(\mathbb{R}^3\times D_\infty)$ such that for sufficiently large $t \geq 1$,
	\begin{equation}\label{UNIes1}
	\begin{split}
	\Vert \sigma(z,w,t)-\sigma_{\infty}(z,w)\Vert_{L^\infty_{z,w}({D^{c}}\times D_\infty\cap \{g(z,w)>t^{-\frac14}, \, |z| \leq t^{\frac{1}{100}}\})}\lesssim t^{-\frac{1}{100}},
	\\
	{\Vert \sigma(z,w,t)-\sigma_{\infty}(z,w)\Vert_{L^\infty_{z,w}(D\times D_\infty\cap \{d(w,C_{\infty})>t^{-\frac14}, \, |z| \leq t^{\frac{1}{100}}\})}\lesssim t^{-\frac{1}{100}}}.
	\end{split}
	\end{equation}
In addition, $\sigma_\infty$ is a distribution for $E_\infty$ in Lemma \ref{Convergenceofelectircfield} in the sense that
\begin{equation}\label{EInfty}
	\begin{split}
		E_{\infty}(a)=\iint_{\{w\in D_\infty\}}\partial_{a^j} G_\infty(a,w)\sigma^2_{\infty}(z,w)dzdw.
	\end{split}
\end{equation}
	
\end{proposition}

\begin{proof}
First we show \eqref{UNIes1}.
For $\kappa_{1},\kappa_{2}>0$, we define the set
\begin{equation*}
\begin{split}
\mathcal{S}_{\kappa_{1},\kappa_{2}}&:=\{ (z,w) \in {D^{c}}\times D_\infty \, : \, g(z,w)>\kappa_{1}^{-1}, \, |z| \leq \kappa_{2} \},
\\
{\mathcal{T}_{\kappa_{1},\kappa_{2}}}&{:=\{ (z,w) \in D \times D_\infty \, : \, d(w,C_{\infty})>\kappa_{1}^{-1}, \, |z| \leq \kappa_{2} \}},
\end{split}
\end{equation*}
where $g$ is defined in \eqref{DefGrazing}. Since $g(z,w)>0$ for all $(z,w)\in{D^{c}}\times D_\infty$, we see that $\cup_{\kappa_{1}>0,\kappa_{2}>0}\mathcal{S}_{\kappa_{1},\kappa_{2}}={D^{c}}\times D_\infty$. 
{It also holds that $\cup_{\kappa_{1}>0,\kappa_{2}>0}\mathcal{T}_{\kappa_{1},\kappa_{2}}=D\times D_\infty$}.
In addition, for given any $(z,w)\in \mathcal{S}_{\kappa_{1},\kappa_{2}} {\cup \mathcal{T}_{\kappa_{1},\kappa_{2}}}$, it is seen from Lemma \ref{EscapeLemma} (ii) that $z+tw\in D$ and $z+tw+\lambda(\ln t)E_\infty(w)\in D$ 
for sufficiently large $t$ with $t\ge 2(\kappa_{1}+\kappa_{2})^{3}+10\vert\lambda\vert|\ln t|\Vert E_\infty\Vert_{L^\infty}+1$
{and $t\ge 10\kappa_{1}^{3}\vert\lambda\vert|\ln t|\Vert E_\infty\Vert_{L^\infty}+1$.}
On the other hand, it holds that for sufficiently large $t$,
	\begin{equation*}
	\begin{split}
		\partial_t\sigma(z,w,t)&=\frac{\lambda}{t}[ E_{\infty}(w)-t^2 E(z+tw,t) ]\cdot\nabla_z\gamma(z+\lambda(\ln t)E_{\infty}(w),w,t)\\
		&\quad-t\lambda[E(z+tw+\lambda(\ln t)E_{\infty}(w),t)- E(z+tw,t) ]\cdot\nabla_z\gamma(z+\lambda(\ln t)E_{\infty}(w),w,t)\\
		&\quad-\lambda E(z+tw+\lambda(\ln t)E_{\infty}(w),t)\cdot\nabla_w\gamma(z+\lambda(\ln t) E_{\infty}(w),w,t).
	\end{split}
\end{equation*}
We note that
\begin{align*}
&E_{\infty}(w)-t^2 E(z+tw,t) 
\\
&= E_{\infty}(w)-t^{2}E^{0}(w,t) + t^{2}E^{0}(w,t) -  t^{2}E(tw,t) +  t^{2}E(tw,t) - t^2 E(z+tw,t).
\end{align*}
Using \eqref{Boot1}, \eqref{Boot2}, \eqref{ControlDiff}, Lemma \ref{Convergenceofelectircfield} and the mean value theorem, we estimate $\partial_t\sigma$ as 
\begin{equation*}
\begin{split}
\vert\partial_t\sigma\vert&\lesssim  t^{-1-\frac{1}{50} }\langle z\rangle
\Vert \mathfrak{1}_{\{b(z+tw)>0\}} \nabla_{z} \gamma\Vert_{L^\infty_{z,w}}
\\
&\quad +t(\ln t) \Vert \nabla E(t)\Vert_{L^\infty_x}\Vert \mathfrak{1}_{\{b(z+tw)>0\}} \nabla_{z} \gamma\Vert_{L^\infty_{z,w}}
+\Vert E(t)\Vert_{L^\infty_x}\Vert \mathfrak{1}_{\{b(z+tw)>0\}} \nabla_{w} \gamma\Vert_{L^\infty_{z,w}}
\\
&\lesssim  \langle z\rangle t^{-1-\frac{1}{50} }.
\end{split}
\end{equation*}
This allows to define $\sigma_\infty$ for any $(z,w)\in \mathcal{S}_{t^{\frac14},t^{\frac{1}{100}}} {\cup \mathcal{T}_{t^{\frac14},t^{\frac{1}{100}}}}$, and also gives \eqref{UNIes1}. 
In addition, Fatou's lemma ensures that $\sigma_\infty\in L^2_{z,w} \cap L^\infty_{z,w}$.
It also follows from $\sigma(t) \in W^{1,\infty}_{z,w} \subset C^{0}_{z,w}$ that $\sigma_\infty\in C^{0}_{z,w}$.

In order to prove \eqref{EInfty}, it suffices to show the convergence
\begin{equation*}
	\begin{split}
	        t^2 E^{0}(a,t)&=\iint_{\{w\in D_\infty\}}\partial_{a^j}G_\infty(a,w)\gamma^{2}(z+\lambda  (\ln t) E_{\infty}(w) ,w,t) dzdw
	        \\
	        &\to \iint_{\{w\in D_\infty\}}\partial_{a^j} G_\infty(a,w)\sigma^2_{\infty}(z,w)dzdw \hbox{ as } t\to\infty.
	\end{split}
\end{equation*}
We observe from $\{\vert z \vert\ge   t^{\frac{1}{1000}} \}{{\subset}}\{\vert z+\lambda(\ln t) E_\infty(w)\vert \geq t^{\frac{1}{2000}}\}$ that 
\begin{equation*}
\begin{split}
\left| \iint_{\{w\in D_\infty\}}\mathfrak{1}_{\{\vert z \vert\ge  t^{\frac{1}{1000}} \}}\partial_{a^j}G_\infty(a,w)\gamma^{2}(z+\lambda  (\ln t) E_{\infty}(w) ,w,t) dzdw \right| &\lesssim t^{-\frac{1}{1000}} N_{1}. 
\end{split}
\end{equation*}
Using \eqref{QuanfVanishingG} and passing to spherical coordinates in $w$, we have
\begin{equation*}
\begin{split}
\left| \iint_{\{{z\in D^{c},} \, w\in D_\infty\}}\mathfrak{1}_{\{ 0<g(z,w) \leq t^{-\frac{1}{4}}\}}\mathfrak{1}_{\{\vert z \vert\leq  t^{\frac{1}{1000}} \}}\partial_{a^j}G_\infty(a,w)\gamma^{2}(z+\lambda  (\ln t) E_{\infty}(w) ,w,t) dzdw \right| \lesssim t^{-\frac14}N_{1}. 
\end{split}
\end{equation*}
{It also holds that}
\begin{equation*}
\begin{split}
{\left| \iint_{\{ z\in D, \, w\in D_\infty\}}\mathfrak{1}_{\{d(w,C_{\infty}) \leq t^{-\frac{1}{4}}\}}\mathfrak{1}_{\{\vert z \vert\leq  t^{\frac{1}{1000}} \}}\partial_{a^j}G_\infty(a,w)\gamma^{2}(z+\lambda  (\ln t) E_{\infty}(w) ,w,t) dzdw \right| \lesssim t^{-\frac14}N_{1}. } 
\end{split}
\end{equation*}
Using these inequalities and \eqref{UNIes1}, we conclude the desired convergence.
Thus \eqref{EInfty} holds.
\end{proof}

\section{Regularity in the domain and proofs of the main theorems}\label{S7}

{In subsection \ref{S7.1}, we study the regularity of $\mu$ in the domain.
Subsection \ref{S7.2} complete the proofs of Theorems \ref{Thm2} and \ref{MainThm}.}

\subsection{Regularity in the domain}\label{S7.1}
We study the regularity in the domain $D$ assuming that the strong solution of \eqref{RVPW} exists.
By the uniqueness, the strong solution can be written as
\begin{equation*}
\begin{split}
\mu(x,v,t):=\mu_0(X(x,v,0,t),V(x,v,0,t)),
\end{split}
\end{equation*}
where the characteristics are defined by
\begin{equation*}
\begin{split}
X(x,v,s,t)&:=x-(t-s)v+\int_{\tau=s}^t(\tau-s)F(X(x,v,\tau,t),\tau)d\tau,\\
V(x,v,s,t)&:=v-\int_{\tau=s}^tF(X(x,v,\tau,t),\tau)d\tau.
\end{split}
\end{equation*}
It will be only Lipshitz\footnote{This is a priori best possible as can be seen by considering the solutions of $\dot{x}=\vert x\vert$.} in $(x,v)$ due to the limited regularity $F\in C^{0,1}_{x}$. This does not enable us to obtain a classical regularity of $\mu$ for all $(x,v,t) \in \mathbb R^{3}_{x} \times \mathbb R^{3}_{v} \times \mathbb R_{+} $.
On the other hand,  if $\mu_{0}=0$ holds on $\mathcal{S}_{in}={\left\{(x,v) \, : \, x \in D^{c}, \ v\cdot \nabla b(x)\geq 0  \right\}}$, we can expect a classical regularity of $\mu$ in the domain $D$ since $F=E\in C^{1}_x$ holds there.
To show this expectation, we start with the following lemma.

\begin{lemma}\label{LemTraj1}
Assume that $b(X(x,v,s,t))=0$ for some $ 0 \leq s \leq t$. Then the following hold:
\begin{enumerate}
\item If $V(x,v,s,t)\cdot{\bf n}(X(x,v,s,t))>0$ and ${\rm supp}(\mu_{0}) \subset \mathcal{S}_{in}^{c}$, then $\mu_0(X(x,v,0,t),V(x,v,0,t))\equiv0$ in a neighborhood of $(x,v)$;
\item If $V(x,v,s,t)\cdot {\bf n}(X(x,v,s,t))\le 0$, then $b(X(x,v,\sigma,t))<0$ for all $s < \sigma\le t$.
\end{enumerate}
\end{lemma}

\begin{proof}[Proof of Lemma \ref{LemTraj1}]
The proof is similar to that of Lemma \ref{SuffCondBC}. 
We note that
\begin{equation}\label{b1}
\begin{split}
\frac{d}{d\sigma}b(X(x,v,\sigma,t))=V(x,v,\sigma,t)\cdot \nabla b(X(x,v,\sigma,t)),
\end{split}
\end{equation}
which coincides with $V\cdot{\bf n}$ when $\sigma=s$, i.e. $b(X)=0$.
In addition, the case that $\{b(X(x,v,s,t))\le0\}$ and $b(X)\in W^{2,\infty}_\sigma$, we can bound the second derivative using \eqref{ExtPotAss3} as
\begin{equation}\label{b2}
\begin{split}
&\frac{d^2}{d\sigma^2}b(X(x,v,\sigma,t))=\{v\cdot\nabla b,\mathcal{H}\}(X,V) \leq 0.
\end{split}
\end{equation}
If $s>0$, $b(X(x,v,s,t))=0$, and $V(x,v,s,t)\cdot{\bf n}(X(x,v,s,t))>0$, we see from \eqref{b1} and \eqref{b2} that 
\begin{equation*}
\begin{split}
b(X(x,v,\sigma,t))<0,\quad V(X(x,v,\sigma,t))\cdot \nabla b(X(x,v,\sigma,t))>0 \quad \hbox{ for } 0\le \sigma < s.
\end{split}
\end{equation*}
Then plugging $\sigma=0$, we deduce from ${\rm supp}(\mu_{0}) \subset \mathcal{S}_{in}^{c}$ that $(X(x,v,0,t),V(x,v,0,t)) \in ({\rm supp}(\mu_{0}))^{c}$.
This holds in the case $s=0$ as well. Thus the assertion (1) holds.
If $V(x,v,s,t)\cdot{\bf n}(X(x,v,s,t))<0$, using \eqref{b1} and \eqref{b2}, we have $b(X(x,v,\sigma,t))<0$ for $s < \sigma\le t$.
Finally, if $V(x,v,s,t)\cdot{\bf n}(X(x,v,s,t))=0$,
it follows from \eqref{BNondegenerate} and \eqref{ExtPotAss2} that
\begin{equation*}
\begin{split}
\left. \frac{d^2}{d\sigma^2}b(X(x,v,\sigma,t)) \right|_{\sigma=s} \le \left.  E(X(x,v,\sigma,t))\cdot\nabla b(X(x,v,\sigma,t)) \right|_{\sigma=s} <0.
\end{split}
\end{equation*}
For $\sigma<s$ close to $s$, it holds that $b(X(x,v,\sigma,t))<0$ and $V(x,v,\sigma,t)\cdot\nabla_xb(X(x,v,\sigma,t))<0$.
Now we can proceed as the case $V(x,v,s,t)\cdot{\bf n}(X(x,v,s,t))<0$. 
Thus the assertion (2) holds.
\end{proof}

\begin{proposition}\label{LemTraj2}
Let $\mu$ be the strong solution of \eqref{RVPW}.
Assume that $E$ and $\nabla E$ belong to $L^{\infty} \cap C^{0}(D \times [0,\infty))$ as well as the initial data $\mu_{0}$ satisfies the condition $\mu_{0}=0$ on $\mathcal{S}_{in}$.
Then $\mathfrak{1}_{\{b>0\}}\mu\in C^1_{x,v,t}$.
\end{proposition}
\begin{proof}
First we consider the case ${\rm supp}(\mu_{0}) \subset \mathcal{S}_{in}^{c}$.
Owing to Lemma \ref{LemTraj1}, 
 if $x\in D$, the characteristic $(X(x,v,\sigma,t),V(x,v,\sigma,t))$ starting from $(x,v)$ at time $t$ satisfies either $b(X(x,v,s,t)) > 0$ or $(X(x,v,s,t),V(x,v,s,t)) \in ({\rm supp} (\mu_{0}))^{c}$ for all $0 \leq s \leq t$.
 Therefore we can use the derivative formula
{{\begin{align*}
&\partial_{x^j} \mu(x,v,t)
\\
&= (\partial_{x^k} \mu_0)(X(x,v,0,t),V(x,v,0,t)) \left[ \delta_{jk}+\int_{\tau=s}^t(\tau-s)\nabla E^{k}(X(x,v,\tau,t),\tau)\cdot\partial_{x_{j}} X(x,v,\tau,t)d\tau \right] \\
&\quad -(\partial_{v^k}\mu_{0})(X(x,v,0,t), V(x,v,0,t)) \int_{\tau=s}^t \nabla E^k(X(x,v,\tau,t),\tau)\cdot \partial_{x^j} X(x,v,\tau,t) d\tau.
\end{align*}}}
This leads to $\partial_{x^j} \mu \in C^{0}_{x,v,t}$.
Similarly, $\partial_{v^j} \mu \in C^{0}_{x,v,t}$ holds.

Next we consider the case $\mu_{0}=0$ on $\mathcal{S}_{in}$.
We can find a sequence $\mu_{0}^{j}$ such that ${\rm supp}\mu_{0}^{j} \subset \mathcal{S}_{in}^{c}$ and $\mu_{0}^{j} \to \mu_{0}$ in ${ C^{1}(K)}$ for any compact set $K \subset \mathcal{S}_{in}^{c}$. 
It is seen from the first case that $\mathfrak{1}_{\{b>0\}}\mu^{j} \in C^{1}_{x,v,t}$, where $\mu^{j}(x,v,t):=\mu_0^{j}(X(x,v,0,t),V(x,v,0,t))$. 
For any $K \subset D \times \mathbb R^{3} \subset \mathcal{S}_{in}^{c}$,  it is seen that  $\mu^{j} \to \mu$ in {$ C^{0}(K)$} and $\mu^{j}$ is a Cauchy sequence in {$ C^{1}(K)$}. From these, it follows that $\mu \in { C^{1}(K)}$. 
Consequently, $\mathfrak{1}_{\{b>0\}}\mu \in C^{1}_{x,v,t}$ holds.
\end{proof}

\subsection{Proofs of the main theorems}\label{S7.2}

We are now ready to prove Theorem \ref{Thm2}.

\begin{proof}[Proof of Theorem \ref{Thm2}]
Using a classical method {in \cite{Arse1975,UO1978}}, we have a time-local solution  {$\mu\in C_{t}([0,T];L^{2}_{x,v}(\mathbb{R}^3\times\mathbb{R}^3)) \cap L^{\infty}_{t}([0,T];W^{1,\infty}_{x,v}(\mathbb{R}^3\times\mathbb{R}^3))$} that satisfies {$ (1+|v|)^{-1} \partial_{t} \mu \in L^{\infty}_{t}([0,T];L^{\infty}_{x,v}(\mathbb{R}^3\times\mathbb{R}^3))$} and the bootstrap assumptions \eqref{Boot1} and \eqref{Boot2} on $[0,T]$. 
In addition, we have a blow-up criterion in the sense that if $\Vert F(x,t)\Vert_{L^1_t([0,T]:W^{1,\infty}_x)}<\infty$, the solution can be extended beyond $T$. 
We can now use Propositions \ref{ProprBoot1} and \ref{ProprBoot2} to obtain a global strong solution satisfying \eqref{GlobalBds2} and \eqref{GlobalBds3}.
Proposition \ref{LemTraj2} ensures that the global strong solution has a classical regularity in the domain $D$.
Finally, the statement on the asymptotic behavior follows from Proposition \ref{Convergenceofgamma}. 
We note that $\sigma_{\infty} \in C^{0}_{x,v}$ follows from \eqref{UNIes1} and  the classical regularity of $\sigma$ in the domain $D$.
%
%
%
%
%
\end{proof}

We are now ready to prove the main theorem.

\begin{proof}[Proof of Theorem \ref{MainThm}]
We extend the initial data $\widetilde{\mu}_0$ to $\mu_0\in C^1_{x,v}(\mathbb{R}^3\times\mathbb{R}^3)$, and then apply Theorem \ref{Thm2}. Using Lemma \ref{SuffCondBC}, we see that $\widetilde{\mu}:=\mathfrak{1}_{\{x\in D\}}\mu$ is the desired solution.
\end{proof}

\begin{appendix}

\section{}\label{S8}
\subsection{Extension of the electrostatic potential}\label{SD}

In this subsection, we suppose that $D$ is an acceptable convex domain and then find an extension map of  the electrostatic potential $\phi$, which satisfies \eqref{ExtPotAss1}--\eqref{ExtPotAss2}.
To this end, we first extend the Green function $G(x,y)$. If $D$ is a half-space, we can use an odd extension for the outside $D^{c}$, else we use a similar idea in \eqref{DefGext}. 

Let us consider the case that  $D$ is a $C^{3}$-domain.
Using the reflexion (see also Figure \ref{Fig1})
\begin{equation}\label{DefRho}
\begin{split}
D^c\ni x\mapsto \rho(x):=2\pi_b(x)-x,
\end{split}
\end{equation}
we define the extended domain
\begin{equation*}
\begin{split}
D^{ext}:=\{x\in D^c \, : \, \rho(x)\in D\}.
\end{split}
\end{equation*}
We set 
\begin{gather*}
\delta^{ext}:= \min\{1,d(D,(D^{ext}\cup \overline{D})^{c})\}>0.
\end{gather*}
The thing to note is that $\pi_{b}(x)=x-b(x)\nabla b(x)$ (for more details, see the proof of Lemma \ref{BProj}).
{From this and Lemma \ref{BProj} with $\pi=\pi_{b}$, it follows that $b$ is $C^{3}$ near $\partial D$, but only on $D^{c}$. For simplicity, we suppose that $b$ is $C^{3}$ on $\{-1<b<0\}$.}
Then it holds that 
\begin{gather}\label{RhoLipshitz}
\rho(x)=x-2b(x)\nabla b(x) \hbox{ for $x\in D^{ext}$}, \quad
\quad D\rho=Id-2\nabla b\otimes\nabla b-2b\nabla^2b.
\end{gather}
We extend the Green function $G(x,y)$ for $y\in D$ and $x\in D^{ext}$ using the reflexion \eqref{DefRho} as
\begin{equation}\label{DefGext}
\begin{split}
G(x,y):=-G(\rho(x),y).
\end{split}
\end{equation}
Now we can extend the electric potential as
\begin{equation}\label{DefExt}
\begin{split}
\phi^0(x,t)&=\iint G(x,y)\mu^2(y,v,t)dydv,\\
\psi(x,t)&:= \chi_{>}(b(x)+\delta_{*})\phi^0(x,t) + A_1(t) V_{*}(b(x)) \quad \hbox{for $x \in \mathbb R^{3}, \ t>0$},
\end{split}
\end{equation}
where the function $V_{*}$ and the constant $\delta_{*}$ are defined by
\begin{gather*}
V_{*}(a):=
\left\{
\begin{array}{ll}
0 & a \geq 0,
\\
\displaystyle  \int_{0}^{a} \frac{s}{1+s^{4}} ds & a < 0,
\end{array}
\right.
\quad
{\delta_{*}:= \delta^{ext}+2}.
\end{gather*}

The properties of the extension are summarized in the following lemma.
\begin{lemma}\label{ExtOperator}
Assume that $D$ is {a $C^{3}$ convex domain}, and $\mu$ satisfies the bootstrap assumptions \eqref{Boot1}.
Then there exists sufficiently large $C_{*}>0$ such that if $A_1(t) = C_{*} \langle t\rangle^{-2}$, then
the extension $\psi$ in \eqref{DefExt} satisfies \eqref{ExtPotAss1}--\eqref{ExtPotAss2}.
\end{lemma}

\begin{proof}[Proof of Lemma \ref{ExtOperator}]
It is clear that \eqref{ExtPotAss1} {holds}, since $\psi(x,t)=\phi^{0}(x,t)=\phi(x,t)$ for $x \in D$.
The Hopf maximum principle immediately gives \eqref{ExtPotAss2}.

To show \eqref{ExtPotAss3}, we use {\eqref{BNondegenerate} and \eqref{bproperty}}. We observe that
\begin{equation*}
\begin{split}
&\{v\cdot\nabla b(x),\mathcal{H}\}\\
&=v^jv^k\partial_{x^j}\partial_{x^k}b(x)-\chi_{>}(b(x)+\delta_{*})\lambda\nabla b(x)\cdot\nabla \phi^0(x)-\vert\nabla b(x)\vert^2[\chi_{>}^\prime(b(x)+\delta_{*})\phi^0+A_1V_\ast^\prime(b)]\\
&\le -\chi_{>}(b(x)+\delta_{*})\lambda\nabla b(x)\cdot\nabla_x\phi^0(x)-A_1 \vert b(x)\vert\left[1+\vert b(x)\vert^4\right]^{-1}.
\end{split}
\end{equation*}
We claim that the right hand side is nonpositive for $A_1(t) = C_{*} \langle t\rangle^{-2}$ with sufficiently large $C_{*}>0$. 
It is sufficient to show that 
\begin{equation}\label{ClaimHopfMaxPrinc}
\begin{split}
\nabla b(x)\cdot\nabla \phi^0(x,t)&\lesssim \vert b(x) \vert   \langle t\rangle^{-3}\langle \ln(2+t) \rangle^{40} 
\quad \hbox{for $x\in D^{ext}$}.
\end{split}
\end{equation}
From \eqref{RhoLipshitz}, \eqref{DefGext}, \eqref{bform1}, {and \eqref{bproperty}}, it follows that
\begin{equation*}
\begin{split}
\nabla b(x)\cdot\nabla G(x,y)=\nabla b(\rho(x))\cdot(\nabla_xG)(\rho(x),y) \quad \hbox{for $x\in D^{ext}$}.
\end{split}
\end{equation*}
Therefore it suffices to estimate $\nabla b(x)\cdot \nabla \phi (x,t)$ for $x \in D$. To this end, we decompose as
\begin{equation*}
\begin{split}
\nabla \phi (x,t)=(\nabla \phi)_{\le L}(x,t)+(\nabla \phi)_{\ge L}(x,t),\quad 
\end{split}
\end{equation*}
where
\begin{align*}
(\nabla \phi)_{\le L}(x,t)&:=\iint\nabla_xG(x,y)\varphi_{\le L}(|x-y|)\mu^2(y,v,t)dydv, \\
(\nabla \phi)_{\ge L}(x,t)&:=\iint\nabla_xG(x,y)\varphi_{\ge L}(|x-y|)\mu^2(y,v,t)dydv.
\end{align*}
Then we use  \eqref{AssumptionExtensionEF} with $L=\langle t\rangle^{3}$ and \eqref{Boot1} to estimate
\begin{equation*}
\begin{split}
\vert [\nabla b(x)\cdot (\nabla \phi)_{\le L}(x,t)]_+\vert&\lesssim \Vert n(x,t)\Vert_{L^\infty_{x}}\Vert [\nabla b(x)\cdot\nabla_xG(x,y)]_+\varphi_{\le L}(|x-y|)\Vert_{L^1_y}
\lesssim \vert b(x)\vert \langle t\rangle^{-3}\langle\ln(2+t)\rangle^{21},
\end{split}
\end{equation*}
where $n (x,t) :=\int\mathfrak{1}_{\{b>0\}} \mu^2(x,v,t)dv$.
On the other hand, a crude estimate gives 
\begin{equation*}
\begin{split}
\vert [\nabla b(x)\cdot (\nabla \phi)_{\ge L}(x,t)]_+\vert &\lesssim  \Vert \mu \Vert_{L^2_{x,v}}^{2}\Vert [\nabla_xb(x)\cdot\nabla_xG(x,y)]_+\varphi_{\ge L}(|x-y|)\Vert_{L^\infty_y}\lesssim \vert b(x)\vert \langle t\rangle^{-3}\langle\ln(2+t)\rangle^{20}.
\end{split}
\end{equation*}
From these, we deduce \eqref{ClaimHopfMaxPrinc} and hence \eqref{ExtPotAss3}.

Let us show the inequalities \eqref{ExtPotAss11} and \eqref{ExtPotAss12}.
It is seen by a direct computation that for $x \in D^{c} $,
\begin{equation*}
\begin{split}
\partial_{x^j}\psi(x,t)
&=-\chi_>(b(x)+\delta_{*})\partial_{x^j}\rho^p(x)\partial_{x_{p}}\phi(\rho(x),t)+\chi_{>}^\prime(b(x)+\delta_{*})\phi(\rho(x),t)\partial_{x^j}b(x)\\
&\quad  +A_1(t) V^\prime_{*}(b(x))\partial_{x^j} b(x)
\end{split}
\end{equation*}
and
\begin{equation*}
\begin{split}
\partial_{x^i}\partial_{x^j}\psi(x,t)&=\chi_{>}(b(x)+\delta_{*})[\partial_{x^j}\rho^p(x)\partial_{x^i}\rho^q(x)\partial_{x^q}\partial_{x_{p}}\phi(\rho(x),t)+\partial_{x^j}\partial_{x^i}\rho^p(x)\partial_{x_{p}}\phi(\rho(x),t)]\\
&\quad+[\chi_{>}^{\prime\prime}(b(x)+\delta_{*})\partial_{x^i}b(x)\partial_{x^j}b(x)+\chi_{>}^\prime(b(x)+\delta_{*})\partial_{x^i}\partial_{x^j}b(x)]\phi(\rho(x),t)\\
&\quad+\chi_{>}^\prime(b(x)+\delta_{*})\partial_{x^i}\rho^p(x)\partial_{x_{p}}\phi(\rho(x),t)\partial_{x^j}b(x)\\
&\quad+A_1(t) V^\prime_{*}(b(x))\partial_{x^j}\partial_{x^i}b(x)+A_1(t) V^{\prime\prime}_{*}(b(x))\partial_{x^i}b(x)\partial_{x^j}b(x).
\end{split}
\end{equation*}
{We note that $b$  is Lipschitz in $\mathbb R^{3}$ (see Lemma \ref{BProj}).}
From {these, \eqref{BNondegenerate},} and the fact $\phi=0$ on $\{b=0\}$, we arrive at \eqref{ExtPotAss11} and \eqref{ExtPotAss12}.
The proof is complete.
\end{proof}

\begin{remark}
To minimize the smoothness requirement on $b$, we could have chosen to extend the electric field in the simpler way
\begin{equation*}
\widetilde{F}^j(x,t)=E^j(\rho(x),t)+A_1(t) V^\prime_{*}(b(x))\nabla b(x).
\end{equation*}
However, $\widetilde{F}$ cannot, in general, be written by a gradient of a certain potential, so we lose the symplectic structure \eqref{RVPW}.
\end{remark}

\subsection{Domain geometry}\label{A2}
In this subsection, we consider the geometry of the domain. For a point $x$ on the boundary $\partial D$, we have the {inward normal vector} ${\bf n}(x)$ and then take an orthogonal frame $(\tau_1(x),\tau_2(x),{\bf n}(x))$. 
Using a translation and/or a rotation, we have a coordinate system so that $x=0$ and {${\bf n}(x)=(0,0,1)$}.
Let us fix it and then the domain is of the form $D=\{ y \, :\, y^3 > w_{x}(y^1,y^2)\}$,
where $w_{x}:\mathbb R^{2} \to \mathbb R$ is a $C^{3}$-function as well as $w_{x}(0)=0$ and $\nabla w_{x}(0)=0$ hold.
There exists an orthogonal matrix $P$ such that
\begin{gather*}
P(x) \tau_1(x) = {\bf e}_{1}, \quad P(x) \tau_2(x) = {\bf e}_{2}, \quad P(x) {\bf n} (x) = {{\bf e}_{3}},
\end{gather*}
where $\{{\bf e}_{1},{\bf e}_{2}, {\bf e}_{3}\}$ is the standard basis in the $y$-coordinate.

We define $II$ and $III$ by
\begin{equation*}
\begin{split}
II(x)=\nabla^2_{y}w_{x}(0),\qquad III(x)=\nabla^3_{y}w_{x}(0).
\end{split}
\end{equation*}

For a point $x_0 \in D$ close enough to the boundary $\partial D$, we can define uniquely its orthogonal projection $\pi(x)$ onto the boundary. We can then extend the orthogonal frames defined above by setting $\tau_j(x)=\tau_j(\pi(x))$ and ${\bf n}(x)={\bf n}(\pi(x))$.
We use the following $3 \times 3$ matrix $M$  by 
\begin{equation*}
\begin{split}
M(x):=
\begin{bmatrix}
[Id_{2}-b(x)\nabla^2_{y^{1},y^{2}}w_{\pi(x)}(0)]^{-1} & 0
\\
0 & 0
\end{bmatrix}.
\end{split}
\end{equation*}
Furthermore, for tangent vectors $\theta,\vartheta \in \mathbb R^{3}$ at $x$ with $\theta\cdot{\bf n}(x)=\vartheta\cdot{\bf n}(x)=0$, we define
\begin{equation*}
\begin{split}
II(x)[\theta,\vartheta]&:=^{t}\![M(x)P(\pi(x)) \theta ] \, II(\pi(x)) \, [M(x) P(\pi(x))\vartheta], \\
III(x)[\tau,\theta,\vartheta]&:=^{t}\![M(x) P(\pi(x)) \theta ] \, III(\pi(x))[P(\pi(x))\tau]] \, [M(x) P(\pi(x)) \vartheta]. 
\end{split}
\end{equation*}

\begin{lemma}\label{BProj}
Suppose that $D$ is a $C^{3}$ convex domain.
We have the following useful identities and properties for $\pi(x)$ and $b(x)$:
\begin{enumerate}
	\item The boundary projection $x\mapsto \pi(x)$ is twice differentiable and
\begin{equation}\label{DerPi}
\begin{split}
\qquad & (\tau\cdot\nabla_x) \pi (x)={}^{t} [ P(\pi(x)) \, \cdot \, ] M(x)[P(\pi(x)) \tau ], 
\quad({\bf n}\cdot\nabla_x)\pi(x)=0,
\\
\qquad &\tau^{c}\theta^a\vartheta^b\partial_{x^a}\partial_{x^b}\pi^c(x)=b(x) III(x)[\tau,\theta,\vartheta], \quad 
{\bf n}^{c}\theta^a\vartheta^b\partial_{x^a}\partial_{x^b}\pi^c(x)=II(x)[\theta,\vartheta],
\\
\qquad & \theta^c\tau^a{\bf n}^b\partial_{x^a}\partial_{x^b}\pi^c(x)=II(x)[\tau,\theta], \quad 
{\bf n}^{c}\tau^a{\bf n}^b\partial_{x^a}\partial_{x^b}\pi^c(x)=0,\quad
{\bf n}^a{\bf n}^b\partial_{x^a}\partial_{x^b}\pi(x)=0,
\end{split}
\end{equation}
where $a,b,c=1,2,3$ and $\tau, \theta, \vartheta$ are tangent vectors at $x$ with $\tau\cdot{\bf n}(x)=\theta\cdot{\bf n}(x)=\vartheta\cdot{\bf n}(x)=0$.
	\item  The function $b$  is Lipschitz in $\mathbb R^{3}$, and $C^1$ in $D^c$ with $\nabla_xb(x)={\bf n}(x)$ on $\partial D=\{b=0\}$ and 
\begin{equation}\label{bform1}
	\begin{split}
		\nabla b(x)=\nabla b(\rho(x)) \quad \mathrm{for} \quad  \rho(x)\in D,
	\end{split}
\end{equation}
where $\rho$ is defined in \eqref{DefRho}.
In addition, it is concave and nondegenerate in the sense that
\begin{equation}\label{BNondegenerate}
\begin{split}
0\le  b(x)\nabla^2b(x)\le Id\quad\hbox{ uniformly when }\quad b(x)\le 0.
\end{split}
\end{equation}
It also holds that  
\begin{gather}\label{bproperty}
{
\vert\nabla b\vert=1, \quad
\nabla^2b(x)\cdot\nabla b(x)=0.
}
\end{gather}
\end{enumerate}
\end{lemma}

\begin{proof}
First we show the assertion (1).
We fix a point $x \in D$ close enough to the boundary $\partial D$, and we choose a coordinate chart as above, where $x=(0,0,b(x))$ and $\pi(x)=0$. We can compute the (squared) distance $d$ between a point $\widetilde{x}=(p^1,p^2,q) \in D$ and a boundary point $(v^1,v^2,w(v)) \in \partial D$:
\begin{equation*}
\begin{split}
d&=\frac{1}{2}[\vert v-p\vert^2+(w(v)-q)^2],
\end{split}
\end{equation*}
where $w(v)=w_{\pi(x)}(v)$.
At a critical point, we obtain the implicit relation
\begin{equation*}
\begin{split}
p&=v+(w(v)-q)\nabla_{v} w(v)=v-q\nabla_{v} w(v)+\nabla_{v}\Bigl(\frac{w^2(v)}{2}\Bigr).
\end{split}
\end{equation*}
This gives
\begin{equation*}
\begin{split}
\frac{\partial v}{\partial p}&= \widetilde{M}:=\left[Id_{2}-q\nabla_{v}^2w(v)+\nabla_{v}^2\frac{w^2(v)}{2}\right]^{-1},\quad\partial_qv^j=\partial_jw(v).
\end{split}
\end{equation*}
Similarly we can compute the second derivative as follows:
\begin{equation*}
\begin{split}
\frac{\partial^2v}{\partial p^2}&=\widetilde{M}\cdot[[q\nabla_{v}^3w(v)-\nabla_{v}^3(w^2/2)(v)]\nabla_pv]\cdot \widetilde{M},\quad \frac{\partial^2v}{\partial q^2}=\nabla^2_{v}w(v)\nabla_{v} w(v), \\ 
\frac{\partial^2v}{\partial p\partial q}&=\widetilde{M}\cdot[\nabla_{v}^2w(v)+[q\nabla_{v}^3w(v)-\nabla_{v}^3(w^2/2)(v)]\nabla _{v}w(v)]\cdot \widetilde{M}.
\end{split}
\end{equation*}
Evaluating these at $p=0$ yields
\begin{equation*}
\begin{split}
 \widetilde{M}_0 &= [Id_{2}-b(x_0)\nabla_{v}^{2}\omega(0))]^{-1},\quad \frac{\partial v}{\partial p}=\widetilde{M}_0,\quad \partial_q v=0,\\
\frac{\partial^2 v}{\partial p^2}&=q\widetilde{M}_0 \nabla_{v}^{3}\omega(0) \widetilde{M}_0,\quad \frac{\partial^2 v}{\partial p\partial q}=\widetilde{M}_0 \nabla_{v}^{2}\omega(0) \widetilde{M}_0,\quad \frac{\partial^2 v}{\partial q^2}=0.
\end{split}
\end{equation*}
Similarly, evaluating at $x_0$, we find
\begin{equation*}
\begin{split}
\frac{\partial [w(v)]}{\partial (p,q)}=0,\quad\frac{\partial^2[w(v)]}{\partial p^2}=\widetilde{M}_0 \nabla_{v}^{2}\omega(0) \widetilde{M}_0,\quad \frac{\partial^2[w(v)]}{\partial p\partial q}=\frac{\partial^2[w(v)]}{\partial q^2}=0.
\end{split}
\end{equation*}
From these, we arrive at the assertion (1).

For the proof of the assertion (2),
it is easy to see that \eqref{bform1} holds. Indeed, it is clear that for $x \in D^c$,
\begin{equation*}
\begin{split}
b(x)=-\vert x-\pi_b(x)\vert,\quad\nabla b(x)=-\frac{x-\pi_b(x)}{\vert x-\pi_b(x)\vert}
={\bf n}(\pi_b(x)).
\end{split}
\end{equation*}
On the other hand, for $\rho(x)\in D$,
\begin{equation*}
\begin{split}
b(\rho(x))=\vert \rho(x)-\pi_b(x)\vert,\quad\nabla b(x)=-\frac{x-\pi_b(x)}{\vert x-\pi_b(x)\vert}= \frac{\rho(x)-\pi_{b}(x)}{|\rho(x)-\pi_{b}(x)|}=\nabla b(\rho(x)).
\end{split}
\end{equation*}
Thus \eqref{bform1} holds.
Furthermore, we can compute the Hessian of $b(x)$:
\begin{equation*}
	\begin{split}
b(x)\nabla^2b(x)=Id-\nabla\pi_{b}(x)-\bf{n}\otimes\bf{n}.
	\end{split}
\end{equation*}
By the first line of \eqref{DerPi}, where the equalities  also hold outside of the domain, we can easily see that \eqref{BNondegenerate} holds. We also see that $\vert\nabla b\vert=1$ and thus $\nabla b$ is in the null space of $\nabla^2b$.
\end{proof}

The next lemmas provide properties of $j_{x}(v)$ and $g(x,v)$ defined in \eqref{j1} and \eqref{DefGrazing}.

\begin{lemma}\label{EscapeLemma}
Suppose that $D$ is a convex domain. Let $x_0\in D$ be fixed.  \\
$(i)$ Let $v\in (\overline{D_{\infty}})^c $ and assume that $x+tv\in D$ for some $x\in\mathbb{R}^3$ and $t\ge0$. Then
\begin{equation*}
\begin{split}
\vert x-x_0\vert\ge \frac12 t j_{x_0}(v)\cdot d(x_0,\partial D) \quad \hbox{for }t\ge 2/j_{x_0}(v).
\end{split}
\end{equation*}
$(ii)$ Let $v\in D_{\infty}$ and assume that $x+tv\notin D$  for some $x\in\mathbb{R}^3$ and $t\ge0$. Then
\begin{equation*}
\begin{split}
\vert x-x_0\vert\ge tg(x,v).
\end{split}
\end{equation*}
\end{lemma}

\begin{proof}
$(i)$ It follows from $x+tv\in D$ that $(x+tv)\cdot\nu(v)>\alpha(v)=(x_0+(1/j_{x_0}(v))v)\cdot\nu(v)$. This gives
\begin{equation*}
\begin{split}
(x-x_0)\cdot\nu(v)> (t-1/j_{x_0}(v))(-v\cdot\nu(v))\ge (t/2)(-v\cdot\nu(v)).
\end{split}
\end{equation*}
By definition of $z$, we see that
\begin{equation*}
-v\cdot\nu(v)=j_{x_0}(v)\cdot[\alpha(v)-x_0\cdot\nu(v)]=j_{x_0}(v)\cdot[x_1\cdot\nu(v)-x_0\cdot\nu(v)]\ge j_{x_0}(v)\hbox{dist}(x_0,\partial D),
\end{equation*}
where $x_{1}$ is the orthogonal projection of $x_0$ onto $\{x\cdot\nu(v)=\alpha(v)\}\subset D^c$.
From these, we arrive at the desired inequality.

$(ii)$ We note that $x+sv\notin D$ for $0\le s\le t$ and
\begin{equation*}
\begin{split}
d(x,x_0)\ge-b(x)=-b(x+tv)+\int_{s=0}^tv\cdot\nabla b(x+sv)ds.
\end{split}
\end{equation*}
Using \eqref{DefGrazing}, we see that the integrand is lower bounded by $g(x,v)$. Thus the desired inequality holds.
\end{proof}

\begin{lemma}\label{GrazingLemma}
{
Suppose that $D$ is a $C^{1}$ convex domain with \eqref{UCAssumption}. Then the condition \eqref{QuanfVanishingG} holds.}
\end{lemma}
\begin{proof} 
We can suppose that $\partial D =\{ x^{3}=f_{\partial D}(x^{1},x^{2}) \} $
and $\partial D_{\infty} = C_{\infty} = \{ x^{3}=f_{C_{\infty}}(x^{1},x^{2}) \}$. 
It is clear that $r f_{C_{\infty}}(x^{1},x^{2}) = f_{C_{\infty}}(r x^{1},r x^{2})$ holds for any $r>0$.
Furthermore, due to \eqref{UCAssumption}, there exists $t_{0}>0$ independent of $(x^{1},x^{2})$ such that 
\begin{gather}\label{tbound1}
f_{C_{\infty}}(x^{1},x^{2}) - t_{0} \leq  f_{\partial D}(x^{1},x^{2}) \leq f_{C_{\infty}}(x^{1},x^{2}).
\end{gather}
In addition, it also follows from \eqref{UCAssumption} that
\begin{gather}
\sup_{(x^{1},x^{2}) \in \mathbb R^{2}} |\nabla f_{\partial D}(x^{1},x^{2})| < \infty.
\label{nbound1}
\end{gather}

First we show that $|\{{w\in\mathbb S^2: 0<j(w)<\alpha }\}| \lesssim \alpha$.
Let $x_0+r^{-1}w \in \partial D$, where $x_{0}=(0,0,1)$ and $r \in (0,\alpha)$.
This means that $ r^{-1}w^{3}= f_{\partial D}(r^{-1}w^{1},r^{-1}w^{2}) -1$.
It is easy to see from \eqref{tbound1} that 
\begin{gather*}
 f_{C_{\infty}}(r^{-1}w^{1},r^{-1}w^{2}) - (t_{0}+1)  \leq  r^{-1}w^{3}= f_{\partial D}(r^{-1}w^{1},r^{-1}w^{2}) -1 \leq  f_{C_{\infty}}(r^{-1}w^{1},r^{-1}w^{2}),
\end{gather*}
which gives
\begin{gather*}
  f_{C_{\infty}}(w^{1},w^{2}) - \alpha (t_{0}+1) \leq f_{C_{\infty}}(w^{1},w^{2}) - r (t_{0}+1)  \leq  w^{3}  \leq  f_{C_{\infty}}(w^{1},w^{2}).
\end{gather*}
Now we observe that
\begin{align*}
|\{{w\in\mathbb S^2: 0<j(w)<\alpha }\}| & \leq |\{{w\in\mathbb S^2:  f_{C_{\infty}}(w^{1},w^{2}) - \alpha (t_{0}+1)  \leq  w^{3}  \leq  f_{C_{\infty}}(w^{1},w^{2})}\}| 
\lesssim \alpha.
\end{align*}

Next let us show that $\sup_{z \in \mathbb R^{3}}\ |\{{w\in\mathbb S^2: 0<g(z,w)<\alpha }\}| \lesssim \alpha$.
Suppose the ray $z+tw$ hits the boundary $\partial D$ at time $t=t_0$ and also write $x=z+t_0w$.  
Then there holds that 
\begin{gather*}
0<w\cdot n(x)=\frac{1}{\sqrt{|\nabla f_{\partial D}|^2+1}} \left[w^{3}- (\partial_{x^{1}} f_{\partial D}) w^{1} - (\partial_{x^{2}} f_{\partial D}) w^{2} \right] < \alpha.
\end{gather*}
From \eqref{nbound1}, we have 
\begin{gather*}
0< \left[w^{3}- (\partial_{x^{1}} f_{\partial D}) w^{1} - (\partial_{x^{2}} f_{\partial D}) w^{2} \right] <  C \alpha,
\end{gather*}
where $C>1$ is independent of $x$. Now we observe that
\begin{align*}
&\sup_{z \in \mathbb R^{3}} |\{{w\in\mathbb S^2: 0<g(z,w)<\alpha }\}| 
\\
&\leq \sup_{z \in \mathbb R^{3}} |\{w\in\mathbb S^2:  (\partial_{x^{1}} f_{\partial D}) w^{1} + (\partial_{x^{2}} f_{\partial D}) w^{2}  \leq  w^{3}  \leq (\partial_{x^{1}} f_{\partial D}) w^{1} + (\partial_{x^{2}} f_{\partial D}) w^{2} + C \alpha \}| 
\lesssim \alpha.
\end{align*}
Thus the condition \eqref{QuanfVanishingG} holds.
\end{proof}

\subsection{Example: flat wall}\label{ExFlatWall}{In the case of a half-space  with $x_0=(0,0,1)$, $j(v)=-v^3_-$, and $g(x,v)=v^3$, we have the following inequality for cases $(i)$ and $(ii)$ in Lemma \ref{EscapeLemma}:
\begin{equation*}
\begin{split}
&\vert x-x_0\vert\ge  -t v^3 -1\geq -\frac{1}{2}tv^3, \quad \hbox{for}\,\, v^3<0,  x+tv\in D, t\geq -2/v^3, \\
&|x-x_0|\geq tv^3+1, \quad \hbox{for}\,\, v^3>0, \, x+tv\notin D.
\end{split}
\end{equation*}
Then \eqref{UCAssumption} and \eqref{QuanfVanishingG} follows immediately. }

We next show \eqref{QuantEG0}--\eqref{AssumptionExtensionEF}. 
The extended distance $b$ and the Green function $G$ are given by
\begin{equation*}
\begin{split}
b(x)&:=x^{3}, \quad 
G(x,y):=-\frac{1}{4\pi}\frac{1}{\vert x-y\vert} + H(x,y), \\
H(x,y)&:=\frac{1}{4\pi}\frac{1}{\vert \sigma(x)-y\vert}, \quad 
\sigma(x^1,x^2,x^3):=(x^1,x^2,-x^3).
\end{split}
\end{equation*}
Obviously \eqref{QuantEG0} holds.
In addition, a direct computation gives
\begin{gather}\label{halfGreen1}
\partial^2_{x^3} G(x,y)+\partial^2_{y^3}G(x,y)=0, \quad
\partial^2_{y^3}G(x,y)\lesssim \frac{b(x)}{|x-y|^4}.
\end{gather}
Let us  show \eqref{AssumptionExtensionEF}.
It is seen that when $x^3 ,y^3 \geq 0$,
\begin{equation*}
\begin{split}
\partial_{x^3}G(x,y)=\frac{1}{4\pi}\left(\frac{x^3-y^3}{\vert x-y\vert^3}\mathfrak{1}_{\{y^3\ge  x^3\}}-\frac{x^3+y^3}{\vert \sigma(x)-y\vert^3}\right)+\frac{1}{4\pi}\frac{x^3-y^3}{{\vert x-y\vert^3}}\mathfrak{1}_{\{0\le y^3\le  x^3\}}.
\end{split}
\end{equation*}
The first term on the right hand side is nonpositive and so negligible,
while the second term satisfies \eqref{AssumptionExtensionEF}.
Thus we have \eqref{AssumptionExtensionEF}.

Finally, we show \eqref{T1Hyp}. 
It is clear that $u_{R}(x)$ depends only on ${x^{3}}$, and so we need only to estimate of the second derivative in $x^{3}$.
We decompose as
\begin{gather*}
\int_{R=0}^{\infty} \partial_{x^{3}}^{2} u_{R}(x) \frac{dR}{R} 
=\int_{R=0}^{b(x)/10} \partial_{x^{3}}^{2} u_{R}(x) \frac{dR}{R} 
+\int_{R=b(x)/10}^{1} \partial_{x^{3}}^{2} u_{R}(x) \frac{dR}{R}.
\end{gather*}
To estimate the first term on the right hand side, we observe by using integration by parts that for $ 0< R < b(x)/10$,
\begin{align*}
& \partial_{x^{3}}^{2} u_{R}(x) 
\\
&= \int_D[\partial_{x^3}^2G(x,y)\varphi(R^{-1}|x-y|)+2\partial_{x^3}G(x,y)\partial_{x^3} (\varphi(R^{-1}|x-y|))+G(x,y)\partial_{x^3}^{2} (\varphi(R^{-1}|x-y|))]dy
\\
&= \int_D[\partial_{x^3}^2H(x,y)\varphi(R^{-1}|x-y|)+2\partial_{x^3}H(x,y)\partial_{x^3} (\varphi(R^{-1}|x-y|))+H(x,y)\partial_{x^3}^{2} (\varphi(R^{-1}|x-y|))]dy.
\end{align*}
On the other hand, a direct computation gives
\begin{gather*}
| \partial_{x^3}^k H(x,y) | \lesssim \frac{1}{\vert \sigma(x)-y\vert^{-1-k}} \hbox{ for } k=0,1,2.
\end{gather*}
From these, we can arrive at 
\begin{gather*}
\left| \int_{R=0}^{b(x)/10} \partial_{x^3}^2 u_{R}(x) \frac{dR}{R} \right| \leq 
\int_{R=0}^{b(x)} \frac{R}{b(x)} \frac{dR}{R} \lesssim 1.
\end{gather*}
What is left is  to show 
\begin{gather}\label{T1Hyp2}
\left| \int_{R=b(x)/10}^{1} \partial_{x^3}^2 u_{R}(x) \frac{dR}{R} \right| \lesssim 1.
\end{gather}
{{When the two derivatives both hit $G(x,y)$ or $\varphi(R^{-1}|x-y|)$, we can easily bound those terms. When one derivative hits $G(x,y)$ and the other hits $\varphi(R^{-1}|x-y|)$, we use the following identity:
\begin{equation*}
	\begin{split}
		\nabla_x [\varphi(R^{-1}|x-y|)]=-R\partial_R[\varphi(R^{-1}|x-y|)]\frac{x-y}{|x-y|^2}.
	\end{split}
\end{equation*}
By  using integration by parts, we can get:
\begin{equation*}
	\begin{split}
&\left\vert\int_{R=b(x)/10}^{1}\int_{y^3>0}\partial_{x^3}G(x,y)\partial_{x^3}[\varphi(R^{-1}|x-y|)]dy	\frac{dR}{R} \right\vert \\
 &\leq\left\vert\int_{y^3>0}\frac{x-y}{|x-y|^2}\partial_{x^3}G(x,y)\left[\varphi(|x-y|)-\varphi(10b(x)^{-1}|x-y|)\right]dy\right\vert \lesssim 1.
	\end{split}
\end{equation*}
Thus \eqref{T1Hyp} holds.}}

\subsection{{ Example: hyperbolic cylinder}}\label{hyper}
For $k>0$ and $l>0$, we define hyperbolic cylinders by
\begin{gather*}
D=\{x \in \mathbb R^{3} \, : \,  x^3>f(x^1):=\sqrt{l^{2}+|kx^{1}|^{2}}-l \}.
\end{gather*}
These are asymptotically like cones of angle $\arctan{|k|}$, since $f(x^{1})/|x^{1}|\rightarrow |k|$ as $x\rightarrow\infty$. Clearly, $D$ is $C^3$ convex domains. Moreover, $|| \nabla f||_{C^2}\leq C(k,l)$, where $C(k,l)$ is a constant which only depends on $k$ and $l$.  
We verify \eqref{UCAssumption} as follows:
\begin{equation*}
	\begin{split}
	\mathfrak{W}(r)&=\sup_{x\in C_{\infty},r/2\leq|x|\leq 2r}	|{\bf{w}}_{C_\infty}(x)-x |\lesssim \sup_{x\in \R^2}[ |kx^{1}|-(\sqrt{l^{2}+|kx^{1}|^2}-l)]\lesssim 1.
	\end{split}
\end{equation*}

We show that this hyperbolic cylinder $D$ also satisfies \eqref{QuantEG0} in the case $k=1$ and $l=\sqrt{2}$.
To this end, applying a suitable rotation and translation, we may assume that $D=\{x \in \mathbb R^{3} \, : \, x^{1}x^{2}>1, \ x^{1}>0 \}$.
Then we set 
\begin{gather*}
c(x):=\sqrt{(x^1)^2+(x^2)^2}.
\end{gather*}
Furthermore, we use the Green function $G_{\infty}$ for the asymptotic domain $D_{\infty}:=\{x^{1}>0,\ x^{2}>0\}$ of the hyperbolic cylinder $D$, namely
\begin{gather*}
G_{\infty}(x,y):=-\frac{1}{4\pi}\left[\frac{1}{\vert x-y\vert}-\frac{1}{\vert \sigma_1(x)-y\vert}-\frac{1}{\vert \sigma_2(x)-y\vert}+\frac{1}{\vert \sigma_1\sigma_2(x)-y\vert}\right],\quad \sigma_jx^p=x^p-2\delta_j^px^p,
\end{gather*}
and also the barrier function $B=B(x):=x^{1}x^{2}-1$, which satisfies
\begin{enumerate}
\item $\Delta B\equiv 0$, and $B>0$ in $D$,
\item $B\equiv 0$ on $\partial D$, and $\vert\nabla_xB(x)\vert\ge  c(x)$ for any $x \in D$,
\item $\vert\nabla^2_x B\vert\le \sqrt{2}$ in $D$.
\end{enumerate}
We divid the proof into the following two lemmas:

\begin{lemma}\label{lemHc1}
For any $\kappa \in  (0,1)$, if $c(x) \geq \kappa |x-y|$, then $|\nabla_{x}^{2} G(x,y)| \leq C(\kappa)|x-y|^{-3}$ holds, where $C(\kappa)>0$ is independent of $x,y$ but depends on $\kappa$.
\end{lemma}
\begin{proof}
We consider a ball centered at $x$ and of radius $r_\kappa:=\kappa\vert x-y\vert/100$. 
If the ball does not intersect the boundaries, 
the proof follows from standard interior convolution estimates for the fundamental solution of the Laplacian:
\begin{equation*}
\begin{split}
H_y(x):=H_y\ast\chi_{r_\kappa}(x),\quad \chi_{r_\kappa}(z):=r_\kappa^{-3}\chi(r_\kappa^{-1}z), \quad H_y(x):=G(x,y)-\frac{1}{4\pi}{\vert x-y\vert}
\end{split}
\end{equation*}
for $\chi$ a smooth, compactly supported radial approximation of unity.
If the ball intersects the boundary, it follows from flattening the boundary and applying standard Schauder boundary estimates. See also \cite[Chapter 2]{CW}.
\end{proof}

\begin{lemma}\label{lemHc2}
If $c(x) \leq |x-y|/10$, then $|\nabla_{x}^{2} G(x,y)| \leq C|x-y|^{-3}$ holds, where $C>0$ is independent of $x,y$.
\end{lemma}
\begin{proof}
It suffices to consider the case $x^3=0$. 
Indeed, let $\tau_a$ be the translation operator defined by $\tau_a (x^1,x^2,x^3):=(x^1,x^2,x^3+a)$.
Then we have $G(\tau_a x,\tau_a y)=G(x,y)$. 
Taking $a=-x^3$, we obtain
\begin{equation*}
	\begin{split}
	G(x,y)=G((x^1,x^2,0),(y^1,y^2,y^3-x^3)), \quad 
	|\nabla_x^2 G(x,y)|=|\nabla_x^2 G((x^1,x^2,0),(y^1,y^2,y^3-x^3)) |.	
	\end{split}
\end{equation*}
Note that $|x-y|=|(x^1,x^2,0)-(y^1,y^2,y^3-x^3)|$. 
Therefore, it suffices to show that
\begin{equation*}
	\begin{split}
	\sup_{x\in D, x^3=0, y\in D}| |x-y|^3 \nabla_x^2 G(x,y)| \leq C.
	\end{split}
\end{equation*}
Hereafter we may suppose that $x^{3}=0$ and $x^{1} \geq x^{2}$. Note that $c(x)=|x|$.

First, we show that if $\vert x\vert<\vert x-y\vert/3$, then
\begin{equation}\label{ComparisonWedge}
\begin{split}
0\le G(x,y)\le G_\infty(x,y) \lesssim  \frac{x^1x^2}{\vert x-y\vert^3}.
\end{split}
\end{equation}
The first two inequalities are clear. 
Moreover, we have $\vert\sigma_1(x)-y\vert^2-\vert x-y\vert^2=4x^1y^1$ and hence
\begin{equation*}
\begin{split}
G_\infty(x,y)\le\frac{4x^1x^2y^1y^2}{\pi\vert x-y\vert^5} \leq \frac{4x^1x^2 (|x|+|x-y|)^{2}}{\pi\vert x-y\vert^5} \lesssim \frac{x^1x^2}{\vert x-y\vert^3}.
\end{split}
\end{equation*}
Thus, \eqref{ComparisonWedge} follows.

We next consider a rescaling
\begin{equation*}
\begin{split}
R&:=|x-y|, \quad D_R:=\{x\, :\, Rx\in D\}, \\
\quad g_R(x)&:=RG(Rx,y),\quad b_R(x):=R^{-2}B(Rx), 
\quad q_R(x):=\frac{g_R(x)}{b_R(x)}=R^3\frac{G(Rx,y)}{B(Rx)},
\end{split}
\end{equation*}
and claim that 
\begin{equation}\label{claim1}
\begin{split}
Q(x):= R^{-3} q_R(R^{-1} x)\lesssim R^{-3}, \quad x\in D \cap B(0,R/4).
\end{split}
\end{equation}
To prove this, we show that
\begin{gather}\label{estdb}
{\rm dist}(x,\partial D_{R}) \lesssim \frac{b_{R}(x)}{|x|}, \quad x \in D_{R} \cap \partial B(0,1/4).
\end{gather}
Since $\nabla b_{R} \neq 0$ on $\partial D_{R}$, there exists a neighborhood $U$ of $\partial D_{R}$ such that every point $x \in U$ admits a unique nearest point $\pi(x) \in \partial D_{R}$. Moreover, the projection map $\pi : U \to \partial D_{R}$ is $C^1$, and we may write
\[
x = \pi(x) + t(x)\,\nu(\pi(x)), \quad 
\nu(\pi(x)) := \frac{\nabla_{x} b_{R}(\pi(x))}{|\nabla_{x} b_{R}(\pi(x))|}, \quad t(x) = |x-\pi(x)|={\rm dist}(x,\partial D_{R}), 
\]
Since $b_{R}(\pi(x))=0$, a Taylor expansion gives
\[
b_{R}(x) = \nabla_{x} b_{R}(\pi(x))\cdot (x-\pi(x)) + {\mathcal R}(x),
\]
where $|{\mathcal R}(x)| \le 2 |x-\pi(x)|^2$.
Substituting $x - \pi(x) = t(x)\,\nu(\pi(x))$ into this expansion and using $|\nabla_x b_{R}(x)| \geq |x|$, we conclude \eqref{estdb}. 
Applying the maximum principle to $g_{R}$ together with  \eqref{QuantEGC}  and \eqref{estdb} on $\partial B(0,1/4) \cap D_R$, and the fact $g_{R}=b_{R}=0$  on $\partial D_R\cap B(0,1/4)$, we obtain
\[
0 \leq q_{R}(x) \lesssim  b_{R}(x), \quad x \in D_{R} \cap B(0,1/4),
\]
which immediately gives \eqref{claim1}.

We also claim that 
\begin{equation}\label{claim2}
\begin{split}
\vert\nabla_xQ(x)\vert&\lesssim \left[(x^1)^2+(x^2)^2\right]^{-\frac{1}{2}} R^{-3}, \quad x\in \partial D \cap B(0,R/10).
\end{split}
\end{equation}
To prove this, fix a point $x = \mathbf{a} = (a,a^{-1},0) \in \partial D \cap B(0,R/10)$ with $a \ge 1$, and consider the ball $B(\mathbf{a}, a/10)$.
We introduce the following set and functions
\begin{equation*}
\begin{split}
S_a&:=\{z:\,\,{\bf a}+az\in  D\}=\left\{a^2(1+z^1)(\frac{1}{a^2}+z^2)\ge 1\right\}, \\
g_a(z)&:=G({\bf a}+az,y),\quad b_a(z)=a^{-2}B({\bf a}+az).
\end{split}
\end{equation*}
Note that these functions are defined for $z\in S_a\cap B(0,1/10)$ and
\begin{equation*}
\begin{split}
S_a&=\{z^2>X_a(z^1)\},\qquad X_a(s)=-\frac{1}{a^2}\frac{s}{1+s}, \\
b_a(z)&=z^2+\frac{1}{a^2}z^1+z^1z^2,\qquad \nabla b_a=\begin{pmatrix}a^{-2}+z^2 \\ 1+z^1\end{pmatrix}.
\end{split}
\end{equation*}
Since $\Vert X_a\Vert_{C^2(B(0,1/4))}\lesssim a^{-2}\le 1$ and $b_a(0,1/10,0)=1/10$, we can apply Theorem 1.1 in \cite{DDSOS} to obtain
\begin{equation*}
\begin{split}
\left\Vert \frac{g_{a}(z)}{b_a(z)}\right\Vert_{C^1(S_a\cap B(0,1/10))}\lesssim \Vert g_a\Vert_{L^\infty(S_a\cap B(0,1/5))}\lesssim \frac{a^2}{ R^3},
\end{split}
\end{equation*}
where we have also used \eqref{ComparisonWedge} and $|{\bf{a}}+az-y|\geq R-|x|-|{\bf{a}}+az|\geq R/10$ in deriving the last inequality.
We deduce that
\begin{equation*}
\begin{split}
(\nabla_xQ)({\bf a})=a^{-1}\nabla_z[Q({\bf a}+az)]  \Big|_{z=0}=a^{-3}\nabla_z\left[\frac{g_a(z)}{b_a(z)}\right] \Big|_{z=0}\lesssim \frac{1}{\vert{\bf a}\vert}\frac{1}{R^3}.
\end{split}
\end{equation*}
Thus, \eqref{claim2} follows.

Finally, we show $|\nabla_{x}^{2} G(x,y)| \leq C|x-y|^{-3}$. We observe that
\begin{equation*}
\begin{split}
\nabla_x^2G(x,y)&=B(x)\nabla^2_xQ(x)+\left[\nabla_xB(x)\otimes\nabla_xQ(x)+\nabla_xQ(x)\otimes\nabla_xB(x)\right]+Q(x)\nabla^2_xB(x).
\end{split}
\end{equation*}
By the maximum principle, it suffices to estimate these terms when $x\in \partial(D\cap B(0,R/4))$. 
On the boundary $\partial B(0,R/4)$, we see that $c(x)=|x| \leq R/4=|x-y|/4$. 
Hence, we apply Lemma \ref{lemHc1} to directly estimate the Hessian $\nabla_{x}^{2}G$. 
On the boundary $\partial D = \{B=0\}$, the first term vanishes, and it suffices to consider the last two terms on the right hand side. 
Using \eqref{claim1} and the boundedness of $\nabla_x^2 B$, the last term is easily controlled. 
For the middle term, we use \eqref{claim2} together with the estimate
\[
|\nabla B(x)| \lesssim \sqrt{(x^1)^2 + (x^2)^2}.
\]
This completes the proof.
\end{proof}

\smallskip

\subsection{Counterexample: corner}\label{ExCorner1}
{
The Green function $G_{c}$ for the corner $D=\{x^1>0, \ x^{2}>0\}$ satisfies \eqref{QuantEG0}, but it does not satisfy \eqref{T1Hyp}. In fact, a direct computation shows that \eqref{T1Hyp} fails in this case.
This observation suggests that it would be of interest to investigate domains whose boundaries are not smooth.
}

\subsection{Green function estimates for $C^{3}$ convex domains}\label{SmoothDomain}
The bound \eqref{QuantEGC} holds for any convex domain, see  \cite[Proposition 2.3]{FrThesis}. 
In the case of the $C^{3}$ convex domains, 
we show \eqref{T1Hyp} and \eqref{AssumptionExtensionEF} concerning {the Green function} $G$.
To this end, we need to employ the approximate Green function $G_W$ discussed below.


\subsubsection{Approximate Green Function}

Given $x \in D$ near the boundary $\partial D$, denote the orthogonal projection onto the boundary $\partial D$ and the reflection to the exterior $D^{c}$ by $\pi(x)$ and $\sigma(x):=2\pi(x)-x$, respectively.
We define the approximate Green function $G_{W}$ by
\begin{equation*}
\begin{split}
G_W(x,y)&:=-\frac{1}{4\pi}\left[\frac{1}{\vert x-y\vert}-\frac{1}{\vert \sigma(x)-y\vert}\right],
\end{split}
\end{equation*}
and the correction by 
\begin{equation}\label{DefHGGW}
\begin{split}
H(x,y):=G(x,y)-G_W(x,y).
\end{split}
\end{equation}
We note that $H$ is a harmonic function in $y$.
In order to estimate $H$, it is convenient to use the acceptable upper bound
\begin{equation}\label{DefinitionMajorant}
\begin{split}
M_x(y):=\frac{1}{\vert \sigma(x)-y\vert}+{\bf n}(x)\cdot F(\sigma(x)-y), \quad \hbox{where } 
F(p):=\frac{p}{\vert p\vert^3} \hbox{ for } p\in\mathbb R^{3},
\end{split}
\end{equation}
where ${\bf n}(x)$ is defined in subsection \ref{A2}.
This is a nonnegative harmonic function in $y$, and we can easily compute its boundary value. 
Indeed, it is seen for the case $x=(0,0,b(x))$ and $y=(v,w(v))$ that 
\begin{equation}\label{DefinitionMajorant2}
\begin{split}
M_x(y)=\frac{1}{\left[\vert v\vert^2+(b(x)+w(v))^2\right]^\frac{1}{2}}\left[1+\frac{b(x)+w(v)}{\vert v\vert^2+(b(x)+w(v))^2}\right].
\end{split}
\end{equation}
For the approximate Green function $G_{W}$, the following lemma holds:

\begin{lemma}\label{EstimatesGW}
If $\tau \cdot{\bf n}(x)=\theta\cdot {\bf n}(x)=\vartheta\cdot {\bf n}(x)=0$, $\vert y\vert\ge 2\vert x\vert$, and $b(x) \leq 10$,  then there holds that
\begin{equation*}
\begin{split}
\vert \tau \cdot(\nabla_x+\nabla_y)\nabla_xG_W(x,y)\vert&\lesssim \frac{1}{\vert x-y\vert^2}+\frac{b(x)}{\vert x-y\vert^3},
\end{split}
\end{equation*}
while
\begin{equation*}
\begin{split}
\vert {\bf n}^i(x){\bf n}^j(x)\partial_{x^i}\partial_{x^j}G_W(x,y)\vert \lesssim\frac{b(x)}{\vert x-y\vert^4}, \quad
\vert \theta^i\vartheta^j\partial_{x^i}\partial_{x^j}G_W(x,y) \vert \lesssim \frac{1}{\vert x-y\vert^2}+\frac{b(x)}{\vert x-y\vert^3}+\frac{b(x)}{\vert x-y\vert^4}.
\end{split}
\end{equation*}
\end{lemma}

\begin{proof}

We compute that
\begin{equation*}
\begin{split}
\partial_{x^i}G_W(x,y)&=\frac{1}{4\pi}\left[F_i(x-y)-\partial_{x^i}\sigma^k(x)F^{k}(\sigma(x)-y)\right],\\
\partial_{x^i}\partial_{x^j}G_W(x,y)&=\frac{1}{4\pi}\left[A_{ij}(x-y)-\partial_{x^i}\sigma^k(x)\partial_{x^j}\sigma^l(x)A_{kl}(\sigma(x)-y)-\partial_{x^i}\partial_{x^j}\sigma^k(x)F^{k}(\sigma(x)-y)\right],
\end{split}
\end{equation*}
where
\begin{equation*}
A_{ab}(x):=\partial_{x^a}F^{b}(x)=\frac{1}{\vert x\vert^3}\left[\delta_{ab}-3\frac{x^{a}x^{b}}{\vert x\vert^2}\right] \ \hbox{for} \ a,b=1,2,3.
\end{equation*}
From \eqref{DerPi} and the formula $\sigma(x)=2\pi(x)-x$, it holds that
\begin{equation*}
\begin{split}
{\bf n}(x)\cdot\nabla_x\sigma^k=-{\bf n}^k(x), \quad
{\bf n}^{a}(x){\bf n}^b(x)\partial_{x^a}\partial_{x^b}\sigma(x)=0, \quad
\tau\cdot\nabla_x (u^{k} \sigma^k)= u^{k}\tau^k+2   u^{a}  N_{ab}(x) \tau^{b} , 
\end{split}
\end{equation*}
where $u \in \mathbb R^{3}$ and $N(x):={}^{t}P(\pi(x))[M(x)-Id_{3}]P(\pi(x))$. 
Therefore we deduce that
\begin{align}
{\bf n}^i{\bf n}^j\partial_{x^i}\partial_{x^j}G_{W}(x,y)&=\frac{1}{4\pi}[A_{ab}(x-y)-A_{ab}(\sigma(x)-y)]{\bf n}^a{\bf n}^b,
\label{GWnn}\\
{\bf n}^i\tau^j\partial_{x^i}\partial_{x^j}G_W(x,y)&=\frac{1}{4\pi}[A_{ab}(x-y)+A_{ab}(\sigma(x)-y)]{\bf n}^a\tau^b  
\notag \\
&\quad -\frac{3}{2\pi}\frac{\langle \sigma(x)-y,N(x)\tau\rangle\langle \sigma(x)-y,{\bf n}\rangle}{\vert \sigma(x)-y\vert^5}
-\frac{1}{2\pi}II(x)[\tau,F(\sigma(x)-y)],
\label{GWnt}\\
\theta^i\vartheta^j\partial_{x^i}\partial_{x^j}G_W(x,y)&=\frac{1}{4\pi}[A_{ab}(x-y)-A_{ab}(\sigma(x)-y)]\theta^a\vartheta^b
\notag \\
&\quad  -\frac{1}{4\pi}A_{ab}(\sigma(x)-y)[((Id+2N)\theta)^a((Id+2N)\vartheta)^b-\theta^a\vartheta^b]
\notag \\
&\quad  -\frac{1}{4\pi}II(x)[\theta,\vartheta]{\bf n}^kF^{k}(\sigma(x)-y)-\frac{b(x)}{4\pi}III(x)[\theta,\vartheta,\tau]\tau^kF^{k}(\sigma(x)-y),
\label{GWtt}
\end{align}
where $\tau,\theta, \vartheta$ are tangent vectors at $x$ with $\tau\cdot{\bf n}(x)=\theta\cdot{\bf n}(x)=\vartheta\cdot{\bf n}(x)=0$.
These give the desired estimates.
\end{proof}

To estimate correction $H$, we often use the following maximum principle:
\begin{lemma}\label{LemMaxPrinciple}
Let $H_{p}=H_{p}(y)$ be a nonnegative harmonic function. 
Assume that $f=f(y)$ is a harmonic function in $D$ and  satisfies the following for $y\in \partial D$:
\begin{equation}\label{MajorizedByM}
\begin{split}
\vert f(y)\vert\le H_{p}(y).
\end{split}
\end{equation}
Then \eqref{MajorizedByM} holds for all $y\in D$.
\end{lemma}

For pure tangential and normal derivatives, we have the next lemma.

\begin{lemma}\label{EstimateH} If $\theta\cdot {\bf n}(x)=\vartheta\cdot {\bf n}(x)=0$ and $b(x) \leq 10$, then there holds that 
\begin{equation}\label{esM1}
	\begin{split}
	       |({\bf n}\cdot\nabla_x) H(x,y) | \lesssim \frac{1}{\vert \sigma(x)-y\vert}  &\hbox{ for } y \in \partial D, \\
		|{\bf n}^i(x){\bf n}^j(x)\partial_{x^i}\partial_{x^j} H(x,y)| + |\theta^i\vartheta^j \partial_{x^i}\partial_{x^j} H(x,y)|\lesssim M_{x}(y) &\hbox{ for } y \in \partial D.
	\end{split}
\end{equation}	
In particular, 
\begin{equation}\label{esM2}
	\begin{split}
	|({\bf n}\cdot\nabla_x) H(x,y) | \lesssim \frac{1}{\vert x-y\vert}  &\hbox{ for } y \in D, \\
		|{\bf n}^i(x){\bf n}^j(x)\partial_{x^i}\partial_{x^j} H(x,y)| + |\theta^i\vartheta^j \partial_{x^i}\partial_{x^j} H(x,y)| \lesssim \frac{1}{|x-y|}  +\frac{1}{|x-y|^{2}}  &\hbox{ for } y \in D.
	\end{split}
\end{equation}

\end{lemma}
\begin{proof}
First we show the bound \eqref{esM1}, but we treat only the second inequality, since the first one can be obtained similarly.
We note that $\partial_{x^i}\partial_{x^j} H|_{y\in\partial D}=-\partial_{x^i}\partial_{x^j} G_{W}|_{y\in\partial D}$ holds.
Using a translation and/or a rotation, we may suppose that $\pi(x)=0$, $x=(0,0,b(x))$, $y=(v,w(v))$, $\omega(0)=0$, and $\nabla \omega(0) =0$.
In the case $|x-y| \geq 1$, it is obvious that  the second inequality in \eqref{esM1} holds.
We treat another case $|x-y| \leq 1$, where $w(v)\lesssim |v|^2$ holds due to $\omega(0)=0$ and $\nabla \omega(0) =0$.
Using \eqref{GWnn}, we have 
\begin{equation*}
	\begin{split}
		\vert {\bf n}^i{\bf n}^j\partial_{x^i}\partial_{x^j}G_{W}(x,y)\vert&\lesssim \left[\frac{1}{\vert x-y\vert^3}-\frac{1}{\vert\sigma(x)-y\vert^3}\right]+\left[\frac{1}{\vert x-y\vert^5}-\frac{1}{\vert\sigma(x)-y\vert^5}\right]\cdot\vert (\sigma(x)-y)\cdot{\bf n}\vert^2\\
&\quad+\frac{\vert ((x-y)\cdot{\bf n})^2-((\sigma(x)-y)\cdot{\bf n})^2\vert}{\vert x-y\vert^5}\\
&\lesssim\frac{b(x)w(v)}{\vert x-y\vert^5} \lesssim M_{x}(y). 
	\end{split}
\end{equation*}
Similarly, it follows from \eqref{GWtt} that
\begin{equation*}
	\begin{split}
		|\theta^i\vartheta^j\partial_{x^i}\partial_{x^j}G_W(x,y)|
		&\lesssim \frac{b(x)w(v)}{|x-y|^5}+\frac{\vert v\vert^2b(x)w(v)}{|x-y|^7}
		 +\frac{b(x)}{\vert \sigma(x)-y\vert^3} + M_{x}(y) + \frac{b(x)}{\vert \sigma(x)-y\vert^2}\\
		 & \lesssim M_{x}(y).
	\end{split}
\end{equation*}
Thus \eqref{esM1} holds.

We next show \eqref{esM2}.
Recall that $M_{x}(y)$ is positive and harmonic in $y$ as well as $\partial_{x^i}\partial_{x^j} H$ is also harmonic in $y$.
Hence, the bound \eqref{esM2} immediately follows from \eqref{DefinitionMajorant} and \eqref{esM1} with the aid of Lemma \ref{LemMaxPrinciple} with $H_{p}=M_{x}, \, |\sigma(x)-\cdot|^{-1}$.
\end{proof}
We now estimate the normal and tangential derivatives of the correction $H(x,y)=G(x,y)-G_W(x,y)$ in \eqref{DefHGGW}. Given $x \in D$, a tangent vector $\tau$ at $x$ with $\tau\cdot{\bf n}(x)=0$, $y\in D$, and $(z,\beta(z))\in \partial D$, where $\beta(0)=\nabla\beta(0)=0$ holds, we define 
\begin{equation}\label{DefVqD1}
\begin{split}
{\bf V}_x(y)&:=\frac{1}{2\pi}II(x)[\tau,  \, \cdot \, ] \frac{1}{|\sigma (x)-y|} 
+ { \frac{{\bf n}^{3}(x)}{4\pi|\sigma (x)-y|}
\begin{bmatrix}
\nabla^{2}_{z} \beta(\pi^{1}(x),\pi^{2}(x)) & 0
\\
0 & 0
\end{bmatrix} \tau},   \\
q(x,y)&:=\frac{1}{2\pi}\frac{{\bf n}(x)\cdot(\pi(x)-y)}{\vert\sigma(x)-y\vert^3},\\
 \mathcal{D}_{\tau}[q]&:=  \tau^{1} \partial_{y^{1}} q(z,\beta(z)) + \tau^{2} \partial_{y^{2}} q(z,\beta(z)) + \tau^{1} [ \partial_{z^{1}} \beta]  [\partial_{y^{3}} q(z,\beta(z))]+ \tau^{2} [ \partial_{z^{2}} \beta][\partial_{y^{3}} q(z,\beta(z))]  \\
&\ = \left[\tau^i\partial_{y^i} q(y) + [\tau \cdot \nabla_{y}(\beta(y^{1},y^{2})-y^{3}) ] \partial_{y^{3}}q(y)\right]_{y=(z,\beta(z))}.
\end{split}
\end{equation}
We note that ${\bf V}_x^k$ and $\hbox{div}_y[{\bf V}_x]$ are harmonic in $y$. 
We decompose the normal and tangential derivatives of $H$ by using these functions.
Then we show in the following lemma that all the terms are acceptable except $\mathcal{D}_{\tau}[q]$.
On the other hand, $\mathcal{D}_{\tau}[q]$ will be handled later by using a harmonic extension defined in subsection \ref{normaltangential1}.

\begin{lemma}\label{StructureMixedDer}
If $b(x)\leq 10$, then there holds that for $y\in \partial D$,
\begin{equation}\label{Hes-tn1}
\begin{split}
\vert {\bf V}_x(y)\vert+
\left\vert \tau^i{\bf n}^{j}(x)\partial_{x^i}\partial_{x^j}H(x,y)-{\rm div}_y[{\bf V}_x]-\mathcal{D}_\tau[q]\right\vert&\lesssim M_x(y),\\
\vert {\bf n}_k(y){\bf V}^k_x(y)\vert&\lesssim 1.
\end{split}
\end{equation}
\end{lemma}

\begin{proof}[Proof of Lemma \ref{StructureMixedDer}]

First we estimate  ${\bf V}_{x}$ in \eqref{Hes-tn1} following the proof of Lemma \ref{EstimateH}.
It is clear that $\vert {\bf V}_x(y)\vert \lesssim M_x(y)$. 
Furthermore, we observe that when $x=(0,0,b(x))$ and $y=(v,w(v)) \in \partial D$, 
\begin{gather*}
\left| {\bf{n}}_{k}(y){\bf{V}}_x^k(y) \right| \lesssim \left| \frac{\nabla_{v}w(v)}{[|v|^2+(b(x)+w(v))^2]^{\frac{1}{2}}} \right| \lesssim 1.
\end{gather*}
Thus the second inequality in \eqref{Hes-tn1} also holds.

It remains to handle the second term on the left hand side of \eqref{Hes-tn1}.
From $\partial_{x^i}\partial_{x^j} H|_{y\in\partial D}=-\partial_{x^i}\partial_{x^j} G_{W}|_{y\in\partial D}$ and \eqref{GWnt}, it holds that for $y \in \partial D$,
\begin{equation}\label{Hes-tn3}
\begin{split}
&\tau^i{\bf n}^{j}(x)\partial_{x^i}\partial_{x^j}H(x,y) -{\rm div}_y[{\bf V}_x]-\mathcal{D}_\tau[q] \\
&=-\frac{1}{4\pi}\left[A_{ij}(x-y)+A_{ij}(\sigma(x)-y)\right]{\bf n}^j\tau^i-\tau^i\partial_{y^i} q 
\\
&\quad +\frac{3}{2\pi}\frac{\langle \sigma(x)-y,N(x)\tau\rangle\langle \sigma(x)-y,{\bf n}\rangle}{\vert \sigma(x)-y\vert^5}\\
&\quad - [\tau \cdot \nabla_{y}(\beta(y^{1},y^{2})-y^{3}) ] \partial_{y^{3}}q
-{\rm div}_y\left[ \frac{{\bf n}^{3}(x)}{4\pi|\sigma (x)-y|}
\begin{bmatrix}
\nabla^{2}_{z} \beta(\pi^{1}(x),\pi^{2}(x)) & 0
\\
0 & 0
\end{bmatrix} \tau \right].
\end{split}
\end{equation}
We can simplify the first two terms on the right hand side as
\begin{equation*}
\begin{split}
\text{(frist two terms)}=&-\frac{1}{4\pi}\left[A_{ij}(x-y)+A_{ij}(\sigma(x)-y)+6\frac{(\pi(x)-y)^j(\sigma(x)-y)^i}{\vert\sigma(x)-y\vert^5}\right]{\bf n}^j\tau^i\\
=&-3\tau^i(\sigma(x)-y)^i{\bf n}^j(x)\left[\frac{(x-y)^j}{\vert x-y\vert^5}+\frac{(\sigma(x)-y)^j}{\vert\sigma(x)-y\vert^5}-\frac{2(\pi(x)-y)^j}{\vert\sigma(x)-y\vert^5}\right].
\end{split}
\end{equation*}
The last two terms can be written by using $\tau \cdot \nabla_{y}(\beta(y^{1},y^{2})-y^{3})|_{y=\pi(x)}=0$ and $\sigma(x)=2\pi(x)-x$ as
\begin{equation*}
\begin{split}
&\text{(last two terms)} \\
&=-\tau \cdot\left[\nabla_{y} \beta(y^{1},y^{2})-\nabla_{y}\beta(\pi^{1}(x),\pi^{2}(x)) \right]\partial_{y^{3}}q 
+ \tau^{j} \partial_{y^{j}y^{i}} \beta (\pi^{1}(x),\pi^{2}(x)) [\sigma(x)-y]^{i}   \frac{{\bf n}^{3}(x)}{4\pi |\sigma(x)-y|^{3}}.
\\
&=- \tau \cdot [\nabla_{y} \beta(y^{1},y^{2})-\nabla_{y}\beta(\pi^{1}(x),\pi^{2}(x))] \frac{{\bf n}(x)\cdot(\pi(x)-y)}{2\pi} \partial_{y^{3}} \left[ \frac{1}{|\sigma(x)-y|^{3}} \right]\\
&\quad  + \tau^{j} \left[\partial_{y^{j}} \beta(y^{1},y^{2})-\partial_{y^{j}}\beta(\pi^{1}(x),\pi^{2}(x)) - \frac{1}{2}\partial_{y^{j}y^{i}} \beta (\pi^{1}(x),\pi^{2}(x)) [y-\pi(x)]^{i}  \right] \frac{{\bf n}^{3}(x)}{2\pi |\sigma(x)-y|^{3}} \\
&\quad  - \frac{1}{2}\tau^{j} \partial_{y^{j}y^{i}} \beta (\pi^{1}(x),\pi^{2}(x)) [x-\pi(x)]^{i}   \frac{{\bf n}^{3}(x)}{2\pi |\sigma(x)-y|^{3}}.
\end{split}
\end{equation*}
In the same way as in the proof of Lemma \ref{EstimateH}, 
we can estimate the right hand side of \eqref{Hes-tn3} by $M_{x}(y)$.
Thus \eqref{Hes-tn1} holds.
\end{proof}

\subsubsection{Harmonic extensions}\label{normaltangential1}

For any harmonic functions $f$, we can reconstruct $f$ from its boundary value.  
Given a function $f=f(x,y)$ defined on $\partial D$, we define $f^{H}=f^H(x,y)$ by solving the  Dirichlet problem
\begin{equation}\label{HarmonicExt}
\begin{split}
\Delta_yf^H=0,\quad y\in D,\qquad f^H=f\quad\hbox{ on }\partial D.
\end{split}
\end{equation}
We call the mapping from the boundary data $f$ to the solution $f^{H}$ of \eqref{HarmonicExt} {\it a harmonic extension}.
%

We define an approximate solution to the Dirichlet problem \eqref{HarmonicExt}. Fix $L \in (0,1)$ and $y_0\in \partial D$. We introduce a coordinate $y=(z,t)\in\mathbb{R}^2\times\mathbb{R}$ such that $y_0=0$ and $\partial D=\{(z,t):\, t=\beta(z)\}$, where $\beta$ satisfies $\beta(0)=\nabla\beta(0)=0$. Given a vector $\underline{\tau} \in \mathbb R^{2}$ and a function $g$ defined on $\mathbb R^{2}$, we define $\underline{f}(z)$ by
\begin{equation}\label{Deffb}
\begin{split}
\underline{f}(z):=\underline{\tau} \cdot\nabla_z [g(z)],
\end{split}
\end{equation}
where $g$ satisfies
\begin{equation}\label{Deffb2}
\begin{split}
\left[L+\vert z\vert\right]^{1+\vert \alpha\vert}\vert \partial_z^\alpha g(z)\vert\le 1, \quad 0\le\vert\alpha\vert \le 3.
\end{split}
\end{equation}
Now we define $f^\flat$ by solving another Dirichlet problem
\begin{equation}\label{FlatExt}
\begin{split}
\Delta_yf^\flat(z,t)=0,\quad (z,t) \in \mathbb R^{2} \times (-L^{2},\infty),\qquad f^\flat(z,-L^2)=\underline{f}(z)=\underline{\tau}\cdot\nabla_z[g(z)].
\end{split}
\end{equation}
Since the boundary is flat, thanks to the Poisson kernel, we can explicitly solve the Dirichlet problem \eqref{FlatExt}:
\begin{equation*}
\begin{split}
f^\flat(z,t)=\int_{\mathbb{R}^2}P_{L^2+t}(z-a)\underline{f}(a)da,
\quad P_{L^2+t}(z-a):= \frac{1}{2\pi} \frac{L^2+t}{\left[\vert z-a\vert^2+(L^2+t)^2\right]^\frac{3}{2}}.
\end{split}
\end{equation*}

\begin{proposition}\label{SolveDir}

Let $0<L\leq 1$, and consider the boundary data $\underline{f}$ defined in \eqref{Deffb}
and the solution $f^{H}$ of \eqref{HarmonicExt} with $f=\underline{f}$.
Then there holds that for $y=(z,t)\in\mathbb{R}^2\times\mathbb{R}_{+}$,
\begin{equation}\label{Estimfflat}
\begin{split}
\vert f^\flat(y)\vert&\lesssim \left[L+\vert  y \vert\right]^{-1}+\left[L+\vert y \vert\right]^{-2},\\
\vert f^H(y)-f^\flat(y)\vert&\lesssim \left[L+\vert y \vert\right]^{-1} + \left[L+\vert  y \vert\right]^{-2}.
\end{split}
\end{equation}
\end{proposition}

\begin{proof}
We start with estimating $f^\flat$. Given $\Lambda:= L+\vert z\vert$, we decompose
\begin{equation*}
\begin{split}
f^\flat(z,t)=A+B,\qquad
A&:= \int_{\mathbb{R}^2}P_{L^2+t}(z-a)\underline{f}(a)\chi(\Lambda^{-1}\vert z-a\vert)da,\\
B&:= \int_{\mathbb{R}^2}P_{L^2+t}(z-a)\underline{f}(a)[1-\chi(\Lambda^{-1}\vert z-a\vert)]da,
\end{split}
\end{equation*}
where $\chi \in C^{\infty}(\mathbb R)$ satisfies $\chi(s)=1$ if $s \leq \frac{1}{4}$, and $\chi(s)=0$ if $s \geq \frac{1}{2}$.
By a direct computation, it is seen that 
 \begin{equation*}
\begin{split}
\vert A\vert&\lesssim \left[L^2+\vert z\vert^2\right]^{-1}\int_{\mathbb{R}^2}\frac{L^2+t}{\left[\vert z-a\vert^2+(L^2+t)^2\right]^\frac{3}{2}}\chi(\Lambda^{-1}\vert z-a\vert)da\lesssim \left[L^2+\vert z\vert^2+t^2\right]^{-1}.
\end{split}
\end{equation*}
To handle $B$, we take advantage of the derivative structure as
 \begin{equation*}
\begin{split}
B&=-\frac{1}{2\pi}\int_{\mathbb{R}^2}g(a)\underline{\tau}\cdot\nabla_a\left[\frac{L^2+t}{\left[\vert z-a\vert^2+(L^2+t)^2\right]^\frac{3}{2}}[1-\chi(\Lambda^{-1}\vert z-a\vert)]\right]da.
\end{split}
\end{equation*}
Using this and dividing into two cases $t \geq |a|$ and  $t \leq |a|$, we have
\begin{equation*}
\begin{split}
\vert B\vert&\lesssim (L^2+t)\int_{\{\vert z-a\vert\ge \frac{1}{4}(L+\vert z\vert)\}}\frac{1}{L+\vert a\vert}\left[\frac{1}{\left[\vert z-a\vert^2+t^2\right]^2} +\frac{1}{\left[\vert z-a\vert^2+t^2\right]^{\frac{3}{2}}}\right]da \\
&\lesssim [L^2+\vert z\vert^2+t^2]^{-1} + [L^2+\vert z\vert^2+t^2]^{-\frac{1}{2}}.
\end{split}
\end{equation*}
These estimates lead to the first line of \eqref{Estimfflat}.
 
We now turn to the difference between $f^H$ and $f^\flat$. We begin with estimating it on the boundary $\{t=\beta(z)\}$ in the case $|z| \leq 1$.
Since the Poisson kernel is normalized and $\int_{\mathbb{R}^2}(z-a)h(|z-a|) da=0$ holds for any function $h$, we observe that
\begin{equation*}
\begin{split}
f^H(z,\beta(z))-f^\flat(z,\beta(z))&=\int_{\mathbb{R}^2}P_{L^2+\beta(z)}(z-a)[\underline{f}(z)-\underline{f}(a)]da= D_{1} + D_{2} +D_{3},
\end{split}
\end{equation*}
where $\Lambda:=L+\vert z\vert/10$ and
\begin{equation*}
\begin{split}
D_{1}&:=\int_{\mathbb{R}^2}P_{L^2+\beta(z)}(z-a)[\underline{f}(z)-\underline{f}(a)-(z-a)\cdot\nabla_z\underline{f}(z)]\chi(\Lambda^{-1}\vert z-a\vert))da,\\
D_{2}&:=\int_{\mathbb{R}^2}P_{L^2+\beta(z)}(z-a)\underline{f}(z)[1-\chi(\Lambda^{-1}\vert z-a\vert)]da,\\
D_{3}&:=-\int_{\mathbb{R}^2}P_{L^2+\beta(z)}(z-a)\underline{f}(a)[1-\chi(\Lambda^{-1}\vert z-a\vert)]da.\\
\end{split}
\end{equation*}

We can estimate $D_{1}$, $D_{2}$, and $D_{3}$ as follows.
From $\beta(0)=\nabla\beta(0)=0$ and $|z| \leq1$, it holds that $\beta(z) \lesssim |z|^{2}$.
We observe that
\begin{equation*}
\begin{split}
\vert D_{1} \vert&\lesssim \int_{\{\vert z-a\vert \leq \frac{1}{2}\Lambda\}} \frac{L^2+\beta(z)}{[\vert z-a\vert^2+(L^2+\beta(z))^2]^\frac{3}{2}}\frac{\vert z-a\vert^2}{[L^2+\vert z\vert^2]^2}da  \lesssim [L+\vert z\vert]^{-1}.
\end{split}
\end{equation*}
A direct computation yeilds
\begin{equation*}
\begin{split}
\vert D_{2} \vert&\lesssim |f(z,\beta(z))|\cdot [L^2+\beta(z)]/ \Lambda  \lesssim \Lambda^{-1} \lesssim [L+\vert z\vert]^{-1}.
\end{split}
\end{equation*}
Using the fact that $f(a,\beta(a))=\underline{\tau}\cdot\nabla_a [g(a)]$, we rewrite $D_{3}$ as
\begin{equation*}
\begin{split}
D_{3}&=\int_{\mathbb{R}^2}g(a)\underline{\tau}\cdot\nabla_a \left[[1-\chi(\Lambda^{-1}\vert z-a\vert)]\cdot P_{L^2+\beta(z)}(z-a)\right]da,\\
\end{split}
\end{equation*}
which gives
\begin{equation*}
\begin{split}
\vert D_{3}\vert \lesssim [L+\vert z\vert]^{-1}.
\end{split}
\end{equation*}
As a result we see that, for $y\in\partial D$,
\begin{equation*}
\begin{split}
\vert f^H(y)-f^\flat(y)\vert  \lesssim M_{x_0}(y),\quad x_0:=(0,0,-L).
\end{split}
\end{equation*}
This inequality is obviously true even in the case $|z| \geq 1$. 
Recall all the functions $f^{H}$, $f^\flat$, and $M_{x_0}$ are harmonic in $D$. 
So we can also extend this inequality inside by Lemma \ref{LemMaxPrinciple} with $H_{p}=M_{x_0}$. 
Thus we obtain the second inequality in \eqref{Estimfflat}.
\end{proof}


\subsection{Proofs of \eqref{T1Hyp} and \eqref{AssumptionExtensionEF}}\label{A.7}

Now we are ready to prove \eqref{T1Hyp} and \eqref{AssumptionExtensionEF}.

\begin{proposition}\label{T1Prop} 
Assume  the domain $D$ is $C^{3}$. Then $\eqref{T1Hyp} $ holds.
\end{proposition}
\begin{proof}
First we consider the interval $0 < R\leq {b(x)}/{10}$. 
We introduce a harmonic function
\begin{equation*}
	\begin{split}
		K(x,y):=G(x,y)+\frac{1}{4\pi}\frac{1}{|x-y|}.
	\end{split}
\end{equation*}
It is seen that 
\begin{equation*}
\begin{split}
\nabla_xu_R(x)&=\nabla_x\left(\int_{y\in D}\varphi(R^{-1}|x-y|)K(x,y)dy\right),\qquad
\vert \nabla_x^\alpha K(x,y)\vert\lesssim [b(x)]^{-1-\vert\alpha\vert},
\end{split}
\end{equation*}
where we have used the maximum principle in deriving the second inequality.
From these, we have
\begin{equation*}
	\begin{split}
	\int_0^{\frac{b(x)}{10}} |\nabla^2 u_R(x) | \frac{dR}{R} \lesssim  \int_0^{\frac{b(x)}{10}} \frac{R}{b(x)} \frac{dR}{R}  \lesssim 1.
	\end{split}
\end{equation*}

Next we consider the interval $b(x)/10 \leq R \leq 1$ and we show 
\begin{equation*}
	\begin{split}
	  \left| \int_{R=b(x)/10}^{1}  \nabla^2u_R(x)  \frac{dR}{R} \right| \lesssim 1.
	\end{split}
\end{equation*}	
Fixed $x \in D$, we can decompose Hessian into tangential and normal derivatives. 
Let $\theta$, $\tau$, and $\vartheta$ be tangent vectors such that
$\theta \cdot {\bf n}(x)=\vartheta \cdot {\bf n}(x)=\tau\cdot {\bf n}(x)=0$. 
We divide the proof into three cases as below.

\medskip

\noindent
Case I (two normal derivatives):
Taking two normal derivatives, we have
\begin{equation*}
	\begin{split}
{\bf n}^a {\bf n}^b \partial_{x^a}\partial_{x^b} u_R(x)&=\int_D  [{\bf n}^a {\bf n}^b \partial_{x^a}\partial_{x^b} G(x,y)]\varphi(R^{-1}|x-y|)dy\\
&\quad+2\int_D[({\bf n}\cdot\nabla_x)G(x,y)]({\bf n}\cdot\nabla_x)\varphi(R^{-1}|x-y|)dy\\
&\quad+\int_D G(x,y){\bf n}^a {\bf n}^b\partial_{x^a}\partial_{x^b}(\varphi(R^{-1}|x-y|))dy\\
&=: I^{1}_{R} + I^{2}_{R} + I^{3}_{R}.
	\end{split}
\end{equation*}
Using Lemmas \ref{EstimatesGW} and \ref{EstimateH}, we observe on the support of $\varphi(R^{-1}|x-y|)$ that 
\begin{equation*}
\begin{split}
\vert {\bf n}^a(x){\bf n}^b(x)\partial_{x^a}\partial_{x^b}G(x,y)\vert&\lesssim\frac{1}{R} +\frac{1}{R^{2}}+\frac{b(x)}{R^4}.
\end{split}
\end{equation*}
Then the first term $I^{1}_{R}$ can be estimated as
\begin{equation*}
	\begin{split}
	  \left| \int_{R=b(x)/10}^{1} I^{1}_{R}  \frac{dR}{R} \right| \lesssim 1.
	\end{split}
\end{equation*}	
To estimate the second term $I_{R}^{2}$,
we use the following equality as long as derivatives hit $\varphi(R^{-1}|x-y|)$:
\begin{equation*}
	\begin{split}
		\nabla_x [\varphi(R^{-1}|x-y|)]=-R\partial_R[\varphi(R^{-1}|x-y|)]\frac{x-y}{|x-y|^2}.
	\end{split}
\end{equation*}
Applying integration by parts in $R$, we find
\begin{equation*}
	\begin{split}
	\left| \int_{b(x)/10}^1 I_{R}^{2} \frac{dR}{R} \right|
	&=\left| \int_{D}\frac{{\bf n}\cdot(x-y)}{\vert x-y\vert^2}({\bf n}\cdot \nabla_x)G(x,y)\cdot [\varphi(|x-y|)-\varphi(10b(x)^{-1}|x-y|)] dy  \right| \lesssim 1.
	\end{split}
\end{equation*}
Similarly, one can handle the third term $I_{R}^{3}$. Thus the desired bound holds.

\medskip

\noindent
Case II (two tangential derivatives):
In the case of two tangential derivatives, we observe on the support of $\varphi(R^{-1}|x-y|)$ using Lemmas \ref{EstimatesGW} and \ref{EstimateH} that
\begin{equation*}
	\begin{split}
		|\theta^i \vartheta^j\partial_{x^i}\partial_{x^j}G(x,y) |\lesssim \frac{1}{R}+\frac{1}{R^{2}}+\frac{b(x)}{R^3}+\frac{b(x)}{R^4}.
	\end{split}
\end{equation*}
Then following the proof of Case I, we can get the desired bound. 

\medskip

\noindent
Case III (one tangential derivative and one normal derivative):  
Taking one tangential derivative and one normal derivative and using integration by parts, we have
\begin{equation*}
\begin{split}
({\bf n}\cdot\nabla_x)(\tau \cdot\nabla_x)u_R(x)&=\int_D \left[({\bf n}\cdot\nabla_x)[(\tau \cdot\nabla_x)+(\tau \cdot\nabla_y)]G(x,y)\right]\varphi(R^{-1}|x-y|)dy\\
&\quad+\int_D\left[[(\tau \cdot\nabla_x)+(\tau \cdot\nabla_y)]G(x,y)\right]({\bf n}\cdot\nabla_x)(\varphi(R^{-1}|x-y|))dy.
\end{split}
\end{equation*}
The second term can be handled by the same technique as in Case I. 
Let us estimate the first term. To this end, we decompose $G=G_{W}+H$.
Using Lemma \ref{EstimatesGW}, we observe on the support of $\varphi(R^{-1}|x-y|)$ that
\begin{equation*}
	\begin{split}
		|({\bf n}\cdot\nabla_x)[(\tau\cdot\nabla_x)+(\tau\cdot\nabla_y)]G_{W}(x,y)|\lesssim \frac{1}{R^2}+\frac{b(x)}{R^3}.
	\end{split}
\end{equation*}
Multiplying this by $R^{-1}\varphi(R^{-1}|x-y|)$ and integrating the result in $y$ and $R$, we have an acceptable bound. 

It remains to estimate the correction $H$.
We observe that 
\begin{equation}\label{Hes-tn11}
	\begin{split}
&\int_D \left[({\bf n}\cdot\nabla_x)[(\tau \cdot\nabla_x)+(\tau \cdot\nabla_y)]H(x,y) \right] \varphi(R^{-1}|x-y|)dy
\\
&=\int_D [({\bf n}\cdot\nabla_x)(\tau \cdot\nabla_x) H(x,y)]\varphi(R^{-1}|x-y|)dy
+ \int_D [({\bf n}\cdot\nabla_x)(\tau \cdot\nabla_y) H(x,y)]\varphi(R^{-1}|x-y|)dy.
	\end{split}
\end{equation}
To estimate the first term on the right hand side of \eqref{Hes-tn11}, 
we decompose the second-order derivatives in $x$ of  the correction $H$ into several components as
\begin{equation}\label{DecomH1}
	\begin{split}
		\tau^i{\bf n}^{j}(x)\partial_{x^i}\partial_{x^j}H(x,y)=
				[\tau^i{\bf n}^{j}(x)\partial_{x^i}\partial_{x^j}H(x,y)-\hbox{div}_y[{\bf V}_x]-(\mathcal{D}_{\tau}q)^{H}] + \mathrm{div}_{y}[{\bf V}_x] + (\mathcal{D}_{\tau}q)^{H},
	\end{split}
\end{equation}
where ${\bf V}_x$, $q$, and $\mathcal{D}_{\tau}$ are defined in \eqref{DefVqD1}, and 
$(\, \cdot \, )^{H}$ denotes the harmonic extension defined in \eqref{HarmonicExt}. 
Now we claim that 
\begin{equation}\label{Hes-tn10}
	\begin{split}
		|\tau^i{\bf n}(x)^j\partial_{x^i}\partial_{x^j}H(x,y)| \lesssim  \frac{1}{R}+\frac{1}{R^{2}}.	
	\end{split}
\end{equation}
Indeed, it is easy to estimate the first term on the right hand of \eqref{DecomH1} by  Lemmas \ref{LemMaxPrinciple} and \ref{StructureMixedDer}.
Furthermore, we can handle the second term by using integration by parts and Lemma \ref{StructureMixedDer}.
For the last term $(\mathcal{D}_{\tau}q)^{H}$, it follows from Proposition \ref{SolveDir} that on the support of $\varphi(R^{-1}|x-y|)$,
\begin{equation*}
	\begin{split}
		\vert (\mathcal{D}_{\tau}q)^{H} \vert\lesssim \frac{1}{R}+\frac{1}{R^2}.
	\end{split}
\end{equation*}
Thus \eqref{Hes-tn10} is valid. 
Multiplying \eqref{Hes-tn10} by $R^{-1}\varphi(R^{-1}|x-y|)$ and integrating the result in $y$ and $R$, we have an acceptable bound.

Finally we handle the second term on the right hand side of \eqref{Hes-tn11}.
It can be written as
\begin{equation*}
       \begin{split}
\text{(2nd term)} &=-\int_{\partial D} \left[ ({\bf n}\cdot\nabla_x) H(x,y) \right] \varphi(R^{-1}|x-y|) (\tau \cdot {\bf n} (y)) dS \\
&\quad -\int_D[({\bf n}\cdot\nabla_x)H(x,y) ] (\tau \cdot\nabla_y) \varphi(R^{-1}|x-y|)dy.
	\end{split}
\end{equation*}
Then Lemma \ref{EstimateH} gives
\begin{gather*}
\left| \int_D \left[ ({\bf n}\cdot\nabla_x)(\tau \cdot\nabla_y) H(x,y) \right] \varphi(R^{-1}|x-y|)dy \right| \lesssim R.
\end{gather*}
Dividing this by $R$ and integrating the result in $R$, we have an acceptable bound.
Thus we arrive at the desired bound of $({\bf n}\cdot\nabla_x)(\tau \cdot\nabla_x)u_R(x)$.
The proof is complete.
\end{proof}

\begin{proposition}\label{GoodL1NormGreenFunctionLem}
Assume the domain $D$ is $C^{2,\alpha}$, and the Green function $G$ satisfies \eqref{QuantEG0}. Then \eqref{AssumptionExtensionEF} holds. 
\end{proposition}
\begin{proof}Assume $L>1$.
We will consider the contribution of various dyadic shells around $x$. For small scales $0\le R\le 10 b(x)$ 
we can use a crude estimate
\begin{equation}\label{EstimPosNablaGSmallScales}
\begin{split}
\Vert [\nabla_xb(x)\cdot\nabla_xG(x,y)]_+\varphi(R^{-1}|x-y|)\Vert_{L^1_y}&\lesssim \Vert \nabla_xG(x,y)\Vert_{L^1(\{R/4\leq\vert x-y\vert\le 4R\})}\lesssim R.
\end{split}
\end{equation}
Let us next treat large scales $ 1/10 \le R \leq L$. Using the Hopf maximum principle, we have 
\begin{equation}\label{PreciseHopf0}
\begin{split}
[\nabla_xb(x)\cdot\nabla_xG(x,y)]_+\le[\nabla_xb(x)\cdot[\nabla_xG(x,y)-\nabla_xG(\pi_b(x),y)]]_+.
\end{split}
\end{equation}
From this, \eqref{QuantEG0}, and the mean value theorem, it is seen that
\begin{equation}\label{EstimPosNablaGLargeScales}
\begin{split}
\Vert [\nabla_xb(x)\cdot\nabla_xG(x,y)]_+\varphi(R^{-1}|x-y|)\Vert_{L^1_y}&\lesssim b(x)
\end{split}
\end{equation}
and that, when $R\ge L$,
\begin{equation*}
\Vert [\nabla_xb(x)\cdot\nabla_xG(x,y)]_+\varphi_{\ge L}(x-y)\Vert_{L^\infty_y}\lesssim b(x)L^{-3}.
\end{equation*}

It remains to consider middle scales $10 b(x) \le R \le 1/10$. We approximate the Green function by $G_W$:
\begin{equation*}
\begin{split}
	[\nabla_{x}b(x) \cdot \nabla_xG(x,y)]_{+}&=[\nabla_x b(x)\cdot  \nabla_{x}H(x,y)+ \nabla_xb(x)\cdot \nabla_{x}G_{W}(x,y)]_{+} \\
	&\leq  [\nabla_x b(x)\cdot  \nabla_{x}H(x,y)]_{+} +[\nabla_xb(x)\cdot \nabla_{x}G_{W}(x,y)]_{+}.
\end{split}
	 \end{equation*}
The second term can be treated similarly as in the case of the half-space.  We consider the correction term $H$. By direct computations, we find that when $y=(p,w(p))\in \partial D$, 
\begin{equation*}
	\begin{split}
		\nabla_{x}b(x)\cdot \nabla_{x}H(\pi_{b}(x),y)=-\frac{w(p)}{4\pi}\left[\frac{1}{[(b(x)-w(p))^2+\vert p\vert^2]^\frac{3}{2}}+\frac{1}{[(b(x)+w(p))^2+\vert p\vert^2]^\frac{3}{2}}\right]\leq 0.
	\end{split}
\end{equation*}
And thus the above is true for $y \in D$ by {the maximum principle}. 
For the first term, using Lemma \ref{EstimateH} with $\nabla_{x} b(x)={\bf n}(x)$, we see {on the support of $\varphi(R^{-1}|x-y|)$} that
\begin{equation}\label{ESigma3}
\begin{split}
	[\nabla_x b(x)\cdot  \nabla_{x}H(x,y)]_{+}&\leq  [\nabla_xb(x)\cdot [\nabla_xH(x,y)-\nabla_{x}H(\pi_{b}(x)),y ]]_{+}\\
	&\leq b(x) \sup_{\theta\in[0,1]} | {\bf n}^i(x) {\bf n}^j(x) \partial_{x^i}\partial_{x^j} H(\theta x+(1-\theta)\pi_{b}(x)),y) | \\
	&\lesssim \frac{b(x)}{R} +  \frac{b(x)}{R^{2}} + \frac{b(x)^2}{R^3}.
\end{split}
\end{equation}
From \eqref{EstimPosNablaGSmallScales}, \eqref{EstimPosNablaGLargeScales}, and  \eqref{ESigma3}, we conclude 
\begin{equation*}
\begin{split}
&\int_{R=0}^{10L}
\Vert [\nabla_xb(x)\cdot\nabla_xG(x,y)]_+\varphi(R^{-1}|x-y|)\Vert_{L^1_y}\frac{dR}{R} 
\\
&\lesssim \int_{R=0}^{10b(x)} R \frac{dR}{R} + \int_{R=10b(x)}^{1/10} (b(x)R^2  +b(x)R + b(x)^2) \frac{dR}{R} + \int_{R=1/10}^{10L} b(x) \frac{dR}{R} 
\lesssim b(x)  \ln (L+1).
\end{split}
\end{equation*}
Thus \eqref{AssumptionExtensionEF} holds.
\end{proof}

\end{appendix}

\section{Declarations}

\subsection{Thanks}

The authors thank Hongjie Dong and Masataka Shibata for interesting discussions on analysis of the Green function. The authors would also like to express sincere gratitude to anonymous referees for valuable comments and suggestions, which have greatly improved the manuscript. The authors acknowledge help and inspiration from Gemini when working on Appendix \ref{hyper}.

\subsection{Data Availability Section}

There are no additional data.

\subsection{Funding and/or Conflicts of interests/Competing interests}

B.P. thanks the hospitality of the E. Schr\"odinger Institute in Vienna where this paper was finalized. B.P. and W.H. were partially supported by NSF grant DMS-2154162 and DMS-2452275.  B.P. and M.S. were partially supported by the Research Institute for Mathematical Sciences, an International Joint Usage/Research Center located in Kyoto University. M. S. was supported by JSPS KAKENHI Grant Numbers 21K03308.

\medskip

The authors declare no conflict of interest.

\end{document}